\documentclass[3p,10pt]{elsarticle}

\usepackage{graphics}
\usepackage{graphicx}
\usepackage{epsfig}
\usepackage{amssymb}
\usepackage{amsthm}
\usepackage{amsmath}
\usepackage{float}
\usepackage{resizegather}
\usepackage{epstopdf}
\usepackage{multirow}
\usepackage{color}
\usepackage{hhline}

\newtheorem{theorem}{Theorem}[section]
\newtheorem{corollary}[theorem]{Corollary}
\newtheorem{proposition}[theorem]{Proposition}
\newtheorem{remark}{Remark}[section]

\newtheorem{lemma}{Lemma}[section]

\newtheorem*{weakproblem}{Weak problem}

\DeclareMathOperator{\dive}{\mathrm{div}}
\DeclareMathOperator{\gra}{\mathbf{grad}}
\DeclareMathOperator{\grae}{\mathrm{grad}}

\newcommand{\g}{\mathbf{g}}
\newcommand{\q}{\mathbf{q}}

\newcommand{\Wur}{\boldsymbol{\mathcal{W}}}
\newcommand{\btau}{\boldsymbol{\tau}}
\newcommand{\press}{p}
\newcommand{\Flux}{\mathbf{\mathcal{F}}}

\journal{Journal of Computational Physics}

\allowdisplaybreaks

\begin{document}

\begin{frontmatter}

\title{A staggered semi-implicit hybrid FV/FE projection method for weakly compressible flows}

\author[mati,imat,itmati]{A. Berm\'udez}
\ead{alfredo.bermudez@usc.es}

\author[LAM]{S. Busto\corref{cor1}}
\ead{saray.busto@unitn.it}

\author[LAM]{M. Dumbser}
\ead{michael.dumbser@unitn.it}
\cortext[cor1]{Corresponding author}

\author[mati,imat,itmati]{J.L. Ferr\'in}
\ead{joseluis.ferrin@usc.es}

\author[upm]{L. Saavedra} 
\ead{laura.saavedra@upm.es}

\author[mati,imat,itmati]{M.E. V\'azquez-Cend\'on}
\ead{elena.vazquez.cendon@usc.es}

\address[LAM]{Laboratory of Applied Mathematics, DICAM, University of Trento, via Sommarive 14, IT-38050 Trento, Italy}

\address[mati]{Departamento de Matem\'atica Aplicada, Universidade de Santiago de Compostela, Facultad de Matem\'aticas. 15782 Santiago de Compostela, Spain}

\address[imat]{Instituto de Matem\'aticas, Universidade de Santiago de Compostela, Facultad de Matem\'aticas. 15782 Santiago de Compostela, Spain}

\address[itmati]{ITMATI, Campus Sur 15782 Santiago de Compostela, Spain}

\address[upm]{Departamento de Matemática Aplicada a la Ingeniería Aeroespacial, Universidad Polit\'ecnica de Madrid, Madrid, Spain}

\begin{abstract}

In this article we present a novel staggered semi-implicit hybrid finite-volume/finite-element (FV/FE) method for the resolution of weakly compressible flows in two and three space dimensions. 
The pressure-based methodology introduced in \cite{BFSV14,BFTVC17} for viscous incompressible flows is extended here to solve the compressible Navier-Stokes equations. Instead of considering the classical system including the energy conservation equation, we replace it by the pressure evolution equation written in non-conservative form. 
To ease the discretization of complex spatial domains, face-type unstructured staggered meshes are considered. A projection method allows the decoupling of the computation of the density and linear momentum variables from the pressure. Then, an explicit finite volume scheme is used for the resolution of the transport diffusion equations on the dual mesh, whereas the pressure system is solved implicitly by using continuous finite elements defined on the primal simplex mesh. Consequently, the CFL stability condition depends only on the flow velocity, avoiding the severe time restrictions that might be imposed by the sound velocity in the weakly compressible regime. High order of accuracy in space and time of the transport diffusion stage is attained using a local ADER (LADER) methodology. Moreover, also the CVC Kolgan-type second order in space and first order in time scheme is considered. To prevent spurious oscillations in the presence of shocks, an ENO-based reconstruction, the minmod limiter or the Barth-Jespersen limiter are employed. To show the validity and robustness of our novel staggered semi-implicit hybrid FV/FE 
scheme, several benchmarks are analysed, showing a good agreement with available exact solutions and numerical reference data from low Mach numbers, up to Mach numbers of the order of unity.
\end{abstract}

\begin{keyword}
weakly compressible flows \sep projection method \sep finite volume method \sep finite element method \sep staggered semi-implicit schemes \sep ADER methodology.
\end{keyword}

\end{frontmatter}

\section{Introduction}
Weakly compressible flows are of great interest for the comprehension of numerous natural phenomena and industrial processes,  such as some geophysical and biological applications, heat exchangers, combustion furnaces or solar energy collectors. Therefore, during the last decades, the research community has made a great effort  in order to develop efficient numerical algorithms for their solution. Two main families of weakly compressible flow solvers can be identified depending on the approach used to compute the pressure field (see \cite{KBW04} for an extended overview).

On the one hand, density-based solvers initially focus on the computation of the density from which the pressure is recovered using the equation of state (EOS). These solvers have been traditionally used for the simulation of compressible flows \cite{EMPRS9,God59,HLL83,LW60,LV02,Munz94,OS82,Roe81,TSS94}. However, the small velocities compared with the sound speed, present in low Mach number flows, introduce a severe restriction on the time step. Even fully implicit algorithms, which avoid the dependency of the time step on the Mach number, may still produce wrong physical results due to excessive numerical dissipation of density-based Godunov-type schemes in the low Mach number limit. Weakly compressible solvers should be  consistent with the features of the incompressible regime in which the pressure becomes a purely hydrodynamic variable or, from the mathematical point of view, a Lagrange multiplier associated to the incompressibility condition (see \cite{MDZ03}). {\color{black} The low Mach number limit, however, poses not only efficiency problems of explicit density-based  Godunov-type schemes, due to the CFL number based on the sound speed, but also the scaling of the numerical dissipation matrix with the Mach number is wrong in the low Mach number limit. To overcome the latter problem, a general approach providing promising results is the so-called preconditioniong of the numerical  dissipation  matrix, see e.g.  \cite{TFL93,GV99,Turk99,MeisterMach,PLK06,Rie11,Mor16,STAJG18,CYX18,MDAB18} for a non-exhaustive list of references. } 

On the other hand, pressure-based solvers directly compute the pressure from a Poisson-type equation derived from the mass and momentum conservation equations. These methods are widely used in the resolution of incompressible Navier-Stokes equations (see, e.g., \cite{HW65,Cho67,Cho68,Pat80,BCG89,CC92,Kan86,Cas14,BFSV14,TD14,TD15,TD16}).  
When applying them to solve weakly compressible flows the major issues arise from considering a variable density field. The system of equations to be solved must be adjusted by including the time derivative term on the mass conservation equation and by adapting also the momentum equation to account for spatial density variation. In addition, a state equation is needed. Several finite volume, finite difference and finite element methods initially developed for incompressible flows have been extended to the weakly compressible and low Mach regimes (see \cite{CG84,Klein95,MRKG03,Dol08,Knik11,CordierDegond,MA16,VCLMPO17}). 
Furthermore, some of the algorithms have been developed to solve all Mach number flows, which proves the potential of the extension of this family of schemes (see \cite{DLP93,PM05,DC13,NRKCC16,XDW17,TD17,RussoAllMach,ABIR19,SIMHD,AbateIolloPuppo}). Note that the methods proposed in \cite{DC13,TD17,SIMHD} make use of the novel flux-vector splitting approach forwarded in \cite{TV12}. 

{\color{black} A key point concerning the development of accurate all-speed flow solvers is the asymptotic preserving property, which ensures the correct behaviour of the method in the incompressible limit. Further details can be found in \cite{CordierDegond,DeT11,AbateIolloPuppo} or in \cite{DLDV18,RussoAllMach,TPK20}, where implicit-explicit schemes (IMEX) are employed. 
	The main idea behind these schemes is the decoupling of the computation of acoustic and material waves. The system related to the material waves is then solved using an explicit scheme whereas the acoustic waves are treated implicitly. Consequently, within the explicit scheme, the eigenvalues of the Jacobian matrix of the flux no longer depend on the sound velocity and therefore the time step restriction, imposed by the CFL condition, is less demanding than for fully explicit schemes. Different splitting strategies have been presented in the literature, including the splitting of the Jacobian matrix in fast and slow waves, \cite{CNPT10}, multiple pressure variables, \cite{PM05},  hyperbolic  splitting, \cite{CordierDegond}, or the already mentioned flux splitting techniques, \cite{TV12,NBALM12}. For a more detailed review on semi-implicit schemes applied to low and all Mach number flows we refer to \cite{PM05,RussoAllMach} and references therein.
}

The methodology presented in this manuscript also belongs to this last family of pressure-based methods and has been extended from the hybrid FV/FE method presented in \cite{BFSV14,BBFSTVC17,BFTVC17,Bus18,BSRVC19} for the resolution of incompressible flows {\color{black} in 3D}.  The main goal of this paper is the development of an efficient methodology for solving weakly compressible flows in two and three space dimensions. To this end, we consider a modified version of the compressible Navier-Stokes equations where the total energy conservation equation has been substituted by the pressure equation written in non-conservative form. Meanwhile, the mass and linear momentum conservation equations are employed. Within this document, we take the equation of state (EOS) of ideal gases, but the above system of equations also applies to a general EOS.  {\color{black} The splitting proposed may be seen as an extension of classical projection methods used in the incompressible regime, see \cite{Guer06}, to weakly compressible flows. First, the transport-diffusion part of the mass, momentum and pressure equation is solved using an explicit scheme. Then, the pressure system is solved implicitly and the intermediate values of the linear momentum are updated with the pressure correction.}

One important feature of the developed algorithm is that it does not only rely on one specific methodology to solve the partial differential equations involved in the compressible model, but combines both finite volume (FV) and finite element (FE) methods in order to provide a more efficient and accurate solver.
Finite volume methods have proven to be highly valuable in the resolution of advection equations so they are used in order to solve transport-diffusion equations (see \cite{LV02,Toro,VC15} and references therein). On the other hand, the capabilities of classical continuous finite elements are exploited for the resolution of the resulting Poisson-type problem for the pressure (see, for instance, \cite{RT83,Cia02}). 

Most of pressure-based solvers extended to solve weakly compressible flows consider structured meshes on two-dimensional geometries. This results in an important restriction on the range of physical phenomena that can be simulated. The numerical method presented in this paper has been developed for general unstructured meshes in order to deal with complex geometries
and to enlarge its applicability. Therefore, from the spatial discretization point of view, we consider a staggered unstructured mesh of the face-type. These kinds of meshes have already been used in the finite volume and discontinuous  Galerkin frameworks with great success (see \cite{BDDV98,THD09,DHCPT10,BFSV14,TD14,TD16,BFTVC17,TD17}). The use of staggered meshes avoids checker-board phenomena that must be corrected when collocated grids are used (see \cite{Tyl16} for a particular example in low Mach number flows). Moreover, the design of the dual mesh allows for an easy implementation of boundary conditions. The nodes of the dual mesh, in which the conservative variables and the density are computed, belong to the boundary as well as the vertex of the primal mesh used to approximate the pressure. The use of these staggered meshes together with a simple specific way of passing the information between them leads to a stable scheme.

To achieve a high order numerical scheme a local ADER method is applied on the finite volume framework. {\color{black} The ADER methodology, Arbitrary high order DErivative Riemann problem, has been first put forward for linear hyperbolic PDE by Toro et al. in \cite{TMN01,Mill99}.} It is a fully discrete one-step approach that relies on non-linear reconstructions and the approximate solution of the generalised Riemann problem, up to any order of accuracy. The resulting schemes are arbitrarily accurate in both space and time in the sense that they have no theoretical accuracy barrier. 
{\color{black} Further developments and applications of different families of ADER methods include, for example, the extension to nonlinear systems of hyperbolic conservation laws  \cite{TT02royal,TT05,TT06,Tak06,Zah08,AboiyarIske,MT14}, non conservative products \cite{DCPT09}, the extension to the DG framework \cite{dumbser_jsc}
	and the use of a general space-time finite element predictor to avoid the cumbersome Cauchy-Kovalevskaya procedure and to allow also the discretization of stiff source terms, see \cite{DET08,DBTM08,Dum10,HD11,BD14,DPRZ16,BCDGP20}. } 
The development of the high order hybrid FV-FE method to be employed is based on the analysis of ADER methodology for the scalar non-linear advection-diffusion-reaction equation presented in \cite{BTVC16}. {\color{black}  Furthermore, in \cite{BFTVC17},  it has been adapted to profit from the benefits of considering the hybrid FV-FE formulation. As a result, the stencil, and thus the computational cost, is reduced with respect to a classical ADER-FV scheme.} 
{\color{black} The use of LADER in the incompressible flow regime involves only the reconstruction and time evolution of linear momentum. To solve the weakly compressible flow model we also extend this procedure to the density and pressure unknowns by using the mass and pressure equations within the Cauchy-Kovalevskaya procedure.}
Let us remark that the use of high order schemes in the presence of high discontinuities {\color{black} requires a limiter} which avoids spurious oscillations that may arise {\color{black} in} the presence of shock waves. 
Many different types of limiters can be found in the bibliography such as ENO and WENO limiters, \cite{HEOC87,SO88,Shu98}, slope and moment limiting, \cite{Kri07,Toro}, artificial 
viscosity approaches, \cite{HH02,PP06,MC11}, or a posteriori limiters \cite{CDL11,DCL12,DLC13,LDD14,DL16}. 
In this paper we consider two different strategies: an ENO reconstruction to be employed within the LADER procedure; or the use of more restrictive limiters like the minmod limiter of Roe, \cite{Roe85,Toro}, and the Barth-Jespersen limiter, \cite{BJ89}. 
{\color{black} Besides, one particular feature of the system of equations to be solved with respect to the incompressible model is the presence of a non-conservative term in the pressure equation.} {\color{black} For its discretization within the finite volume framework, in this paper we employ a path conservative approach, see  \cite{Par06,CFFP09,DCPT09,DPRZ16,GCD18}.}

{\color{black}
	To obtain the pressure system in the incompressible case the mass and momentum conservation equations were used. For the weakly compressible case we make use of the momentum and pressure equations. The weakly compressible formulation used in this paper requires not only the stiffness matrix arising in the incompressible case, but also a mass matrix. Moreover, new terms appear on the right hand side. To properly approximate them we include a pre-projection stage in which data coming from the transport diffusion stage is adequately interpolated between meshes and used to compute the sound velocity needed in the pressure system.
}

The rest of the paper is organized as follows. In Section \ref{sec:math}, the compressible Navier-stokes equations are recalled and the energy conservation equation is substituted by the non-conservative form of the pressure equation. Further details on its derivation are included in \ref{sec:app_equations}. In Section \ref{sec:ndfv}, the numerical scheme is presented. The unstructured staggered grid is described in detail and the overall algorithm is introduced. 
The extension of the LADER methodology to achieve a second order in space and time scheme for the compressible model is carefully detailed. Special attention is paid to the coupling between transport and projection stages. 
Finally, Section \ref{sec:numericalresults} is devoted to the careful testing of the new method by solving a large set of different numerical benchmark problems. A modified version of the Taylor-Green vortex including gravity effects is used, to verify the order of accuracy of the scheme. Furthermore, some Riemann problems are simulated to study the behaviour of the new method in the presence of weak discontinuities. Further benchmarks include explosion problems in two and three dimensions, natural convection tests and very low Mach number flows. Overall, all numerical test problems indicate that the method is suitable for the simulation of weakly compressible flows, from low Mach numbers up to Mach numbers of the order of unity. The paper closes with some concluding remarks and an outlook to future research in Section \ref{sec:conclusions}.

\section{Governing equations} \label{sec:math}
The compressible Navier-Stokes equations are classically presented in terms of density, $\rho$, momentum density, $\rho \mathbf{u}$, and {\color{black} specific total energy, $E$}, as

\begin{gather}
	\frac{\partial \rho}{\partial t} + \dive \left( \rho\mathbf{u}\right)=0,\label{eq:mass1}\\
	\frac{\partial \rho\mathbf{u}}{\partial t} + \dive \left( \rho\mathbf{u}\otimes \mathbf{u}\right)  + \grae \press - \dive \btau = \rho \g,\label{eq:momentum1}\\
	\frac{\partial {\color{black} \rho E}}{\partial t} + \dive \left[ \mathbf{u}\left({\color{black} \rho E} + \press \right) \right]   - \dive \left(\btau\mathbf{u}\right) +\dive \q = \rho \g \cdot \mathbf{u},\label{eq:totenerg1}\\
	\btau = \mu \left(\gra \mathbf{u} + \gra \mathbf{u}^{T} \right) -\frac{2}{3} \mu \left( \dive \mathbf{u}  \right) \mathbf{I}, \label{eq:stresstensor}\\
	\q = -\lambda \gra \theta, \label{eq:heatflux}
\end{gather}
with $\press$ and $\theta$ the pressure and the temperature, respectively, $\btau$ the viscous part of the Cauchy stress tensor, $\mu$ the dynamic viscosity, $\g$ the gravity vector, $\q$ the heat flux, and $\lambda$ the thermal conductivity coefficient. For smooth solutions, an equivalent formulation can be obtained by replacing the total energy conservation law \eqref{eq:totenerg1} with a time evolution equation for the pressure. In the following we will denote the linear momentum with $\mathbf{w}_{\mathbf{u}} := \rho\mathbf{u}$. 

\begin{gather}
	\frac{\partial \rho}{\partial t} + \dive \left( \mathbf{w}_{\mathbf{u}}\right)=0,\label{eq:mass2}\\
	\frac{\partial \mathbf{w}_{\mathbf{u}}}{\partial t} + \dive \left( \frac{1}{\rho} \mathbf{w}_{\mathbf{u}}\otimes \mathbf{w}_{\mathbf{u}}\right)  + \grae \press - \dive \btau = \rho \g,\label{eq:momentum2}\\
	\frac{\partial \press}{\partial t}  + \mathbf{u}\cdot \grae \press - c^2 \mathbf{u}\cdot\grae  \rho + c^2 \dive \left(\rho\mathbf{u}\right) + \left(\gamma-1\right)\left(\mathrm{div} \mathbf{q} - \btau \cdot \gra \mathbf{u} \right) = 0,\label{eq:pressure2}
\end{gather} 
where $c^2 = (\partial p / \partial \rho )_s$ is the square of the isentropic sound speed $c$. For the ideal gas equation of state (EOS) we have $c=\sqrt{\frac{\gamma \press}{\rho}}$, with $\gamma$ the adiabatic index.  
Furthermore, one has the well-known relation 
\begin{equation}
	\press = \rho R \theta
\end{equation}
with $R$ being the specific gas constant, $R=\mathcal{R} \displaystyle{\sum_{l=1}^{N_{e}} \frac{\mathrm{y}_{l}}{\mathcal{M}_{l}}}$, $\mathcal{R}$ the universal gas constant (8.314 J/molK), $\mathcal{M}_{l}$ the molar mass of the $l$-th species, $\mathrm{y}_{l}$ its mass fraction and $N_{e}$ the number of species of the mixture. Although the method is presented here for ideal gases only, for the sake of simplicity, any other EOS could in principle be also employed. To change the EOS in the numerical method, it is enough to introduce the proper expressions for the computation of the isentropic sound speed $c$ and the temperature $\theta$ in terms of the density and the pressure. It is precisely for this reason that we have chosen to use the pressure evolution equation rather than the total energy conservation law, since it allows a very simple and straightforward extension of the algorithm to general equations of state. 
{\color{black} Using formal asymptotic analysis it can be seen that in the incompressible limit, i.e. for $c \to \infty$ and thus for $M \to 0$, from the pressure equation  \eqref{eq:pressure2} and the momentum conservation equation \eqref{eq:momentum2} one obtains  
	\begin{equation}
		\dive \mathbf{u} \to 0,  
	\end{equation} 
	which is the well-known divergence-free condition of the velocity field for incompressible flows. Rigorous mathematical proofs concerning the asymptotic limit 
	and numerical methods based on these asymptotic results can be found in   \cite{KlaMaj,KlaMaj82,Klein95,Klein2001,MRKG03,MDR07}. }   
Further details on the derivation of system \eqref{eq:pressure2} from \eqref{eq:totenerg1} are included in \mbox{\ref{sec:app_equations}}.

The term $\boldsymbol{\tau}\cdot \gra \mathbf{u}$ in the pressure equation \eqref{eq:pressure2}, corresponds to the energy dissipation related to the gas viscosity, i.e., due to the shear between the molecules of the gas.  
Usually, for low Mach number and high Reynolds number flows, its magnitude is small with respect to the remaining terms in \eqref{eq:pressure2}, so it can be neglected in many numerical experiments.

\section{Numerical discretization} \label{sec:ndfv}
The numerical discretization of the complete system, \eqref{eq:mass2}-\eqref{eq:pressure2}, is performed by extending
the projection method put forward in \cite{BFSV14} and \cite{BFTVC17}.
The proposed methodology decouples the computation of the linear momentum and the pressure. Regarding the numerical scheme, finite volume methods are used for the approximation of the {\color{black} transport-diffusion} equations, whereas finite elements are employed to solve the pressure system.

To derive our numerical method, we consider the following semi-discretization of the governing PDE system, where first time is discretized, while all spatial operators are still kept continuous. Only later, the spatial discretization will be performed with appropriate discrete operators and methods, namely FV methods for the transport-diffusion terms and continuous FE method for the resulting semi-discrete pressure equation. 

We therefore start by considering a two-stage in time discretization scheme:
in order to get the solution at time $t^{n+1}$, we use the previously obtained approximations
$\mathbf{W}_{\mathbf{u}}^n$ of the linear momentum density $\mathbf{w}_{\mathbf{u}}(x,y,z,t^n)$,
$\mathbf{U}^n$ of  velocity $\mathbf{u}(x,y,z,t^n)$,
$\rho^{n}$ of  density  $\rho\left(x,y,z,t^n\right)$, $\theta^{n}$ of temperature $\theta \left(x,y,z,t^n\right) $
and $\press^n$ of  pressure $\press(x,y,z,t^n)$, and
compute $\mathbf{W}_{\mathbf{u}}^{n+1}$, $\rho^{n+1}$, $\theta^{n+1}$, and $\press^{n+1}$ from the following system
of equations:

\begin{align}
	& \frac{1}{\Delta t} \left(\rho^{n+1}-\rho^{n}\right) + \dive \Flux^{\rho}\left(\mathbf{W}^n\right)= 0, \label{eq:disc_mass}\\
	&\frac{1}{\Delta t}\left( \widetilde{\mathbf{W}}_{\mathbf{u}}-\mathbf{W}_{\mathbf{u}}^{n}\right) +\dive \Flux^{\mathbf{W}_{\mathbf{u}}}\left(\mathbf{W}^n\right)-\dive  \btau^n = \rho^{n}\g,
	\label{eq:disc_mom_tilde}\\
	&\frac{1}{\Delta t}\left( \mathbf{W}_{\mathbf{u}}^{n+1}-\widetilde{\mathbf{W}}_{\mathbf{u}}\right) +\grae \press^{n+1}=0,\label{eq:disc_mom_n1}\\
	&\frac{1}{\Delta t}\left( \widetilde{\press}_{\press}-\press^{n}\right) + \mathbf{U}^{n} \cdot \grae \press^{n} + \left(\gamma-1\right) \mathrm{div} \mathbf{q}^{n} = 0,\label{eq:disc_pres_tilde}\\
	& \frac{1}{\Delta t} \widetilde{\press}_{\rho} - c^{2}\mathbf{U}^{n} \cdot \grae \rho^{n}  = 0, \label{eq:disc_pres_tilder}\\
	&\frac{1}{\Delta t}\left( \press^{n+1}- \widetilde{\press}\,\right) + c^{2} \dive \mathbf{W}_{\mathbf{u}}^{n+1} - \left(\gamma-1\right) {\color{black} \boldsymbol{\tau}^{n}}\cdot \gra {\color{black} \mathbf{u}^{n}} = 0, \quad \widetilde{\press} := \widetilde{\press}_{\rho} + \widetilde{\press}_{\press}, \label{eq:disc_pres_n1}
\end{align}
where $\mathbf{W}^n=\left(\rho, \mathbf{W}_{\mathbf{u}}^{T}, p \right)^{T}$ is the vector of  unknowns, and $\Flux^{\rho}\left(\mathbf{W}^n\right):=\mathbf{W}_{\mathbf{u}}^{n}$ and $\Flux^{\mathbf{W}_{\mathbf{u}}}\left(\mathbf{W}^n\right):= \frac{1}{\rho^{n}} \mathbf{W}^n_{\mathbf{u}}\otimes \mathbf{W}^n_{\mathbf{u}}$ are the convective fluxes related to the mass and momentum conservation equations, respectively. 
Let us notice that adding \eqref{eq:disc_mom_tilde} and \eqref{eq:disc_mom_n1} (respectively, \eqref{eq:disc_pres_tilde}, \eqref{eq:disc_pres_tilder} and \eqref{eq:disc_pres_n1}) we get a time discretization of \eqref{eq:momentum2} (respectively, of \eqref{eq:pressure2}).

Accordingly, we propose an algorithm involving four stages:
\begin{itemize}
	\item \textit{Transport-difusion stage:} 
	equations \eqref{eq:disc_mass}, \eqref{eq:disc_mom_tilde}, and \eqref{eq:disc_pres_tilde} are solved through a finite volume method providing the value of the density at the new time, $\rho^{n+1}$, and intermediate approximations for linear momentum and pressure, namely $\widetilde{\mathbf{W}}_{\mathbf{u}}$ and $\widetilde{\press}_{p}$ {\color{black} on the dual mesh}.
	
	\item \textit{Pre-projection stage:} 
	the value of $\widetilde{\press}_{\rho}$ is computed solving \eqref{eq:disc_pres_tilder} in a finite volume fashion{\color{black}, i.e. 
		\eqref{eq:disc_pres_tilder} is integrated on the control volumes related to the primal mesh and the average value of $\widetilde{p}_{\rho}$ on each cell is approximated.
	}
	Then, the intermediate linear momentum and pressure are interpolated from the dual mesh to the primal mesh.
	
	\item  \textit{Projection stage:} a finite element method is applied to system \eqref{eq:disc_mom_n1}, \eqref{eq:disc_pres_n1} 
	in order to determine the pressure at the new time instant, $\press^{n+1}$, {\color{black} on the vertices of the primal mesh.}
	
	\item  \textit{Post-projection stage:}
	the new pressure is employed in equation \eqref{eq:disc_mom_n1} to update the linear momentum, $\mathbf{W}_{\mathbf{u}}^{n+1}$. Moreover, the temperature at the new time instant, $\theta^{n+1}$, is recovered from the new pressure and density by using the EOS.
\end{itemize}
Before providing a better description of each stage, we will introduce the space discretization.

\subsection{Unstructured staggered grid}
\begin{figure}[h]
	\begin{minipage}{0.6\linewidth}
		\begin{figure}[H]
			\begin{center}
				\begin{picture}(100,90)
				\put(0,0){\makebox(100,90){
						\vspace{-2cm}
						\includegraphics[width=6.2cm]{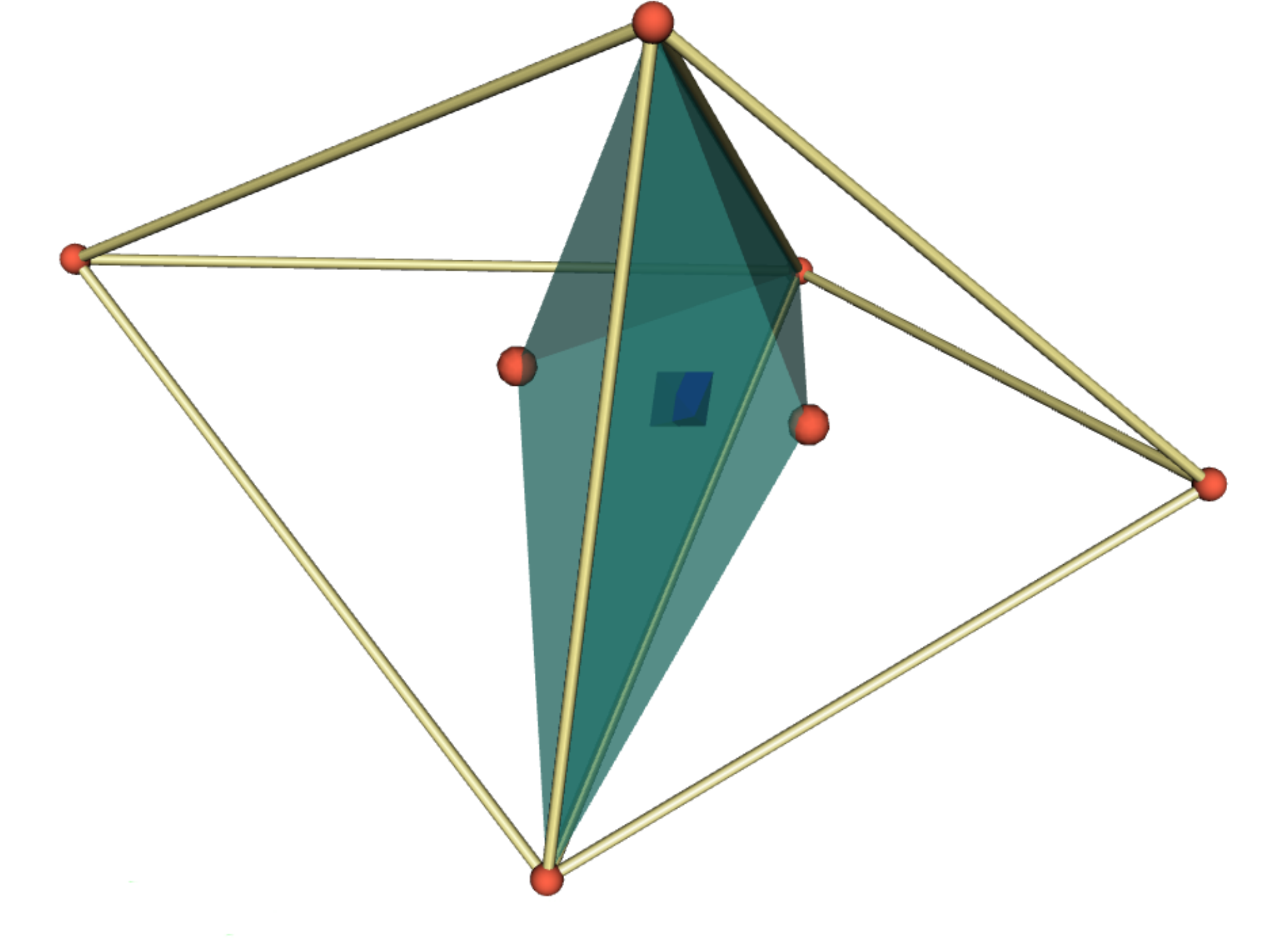}}}
				\put(33,-55){\normalsize{$V_{2}$}}
				\put(75,48){\normalsize{$V_{3}$}}
				\put(47,90){\normalsize{$V_{1}$}}
				\put(-40,52){\normalsize{$V_{4}$}}
				\put(140,12){\normalsize{$V_{4}'$}}
				\put(22,30){\normalsize{$B$}}
				\put(84,20){\normalsize{$B'$}}
				\put(56,33){\normalsize{$N_i$}}
				\end{picture}
				\vspace{1.5cm}
			\end{center}
		\end{figure}
	\end{minipage}
	\begin{minipage}{0.3\linewidth}
		\begin{figure}[H]
			\begin{center}
				\begin{picture}(100,90)
				\put(0,0){\makebox(100,90){
						\vspace{-2cm}
						\includegraphics[width=3.6cm]{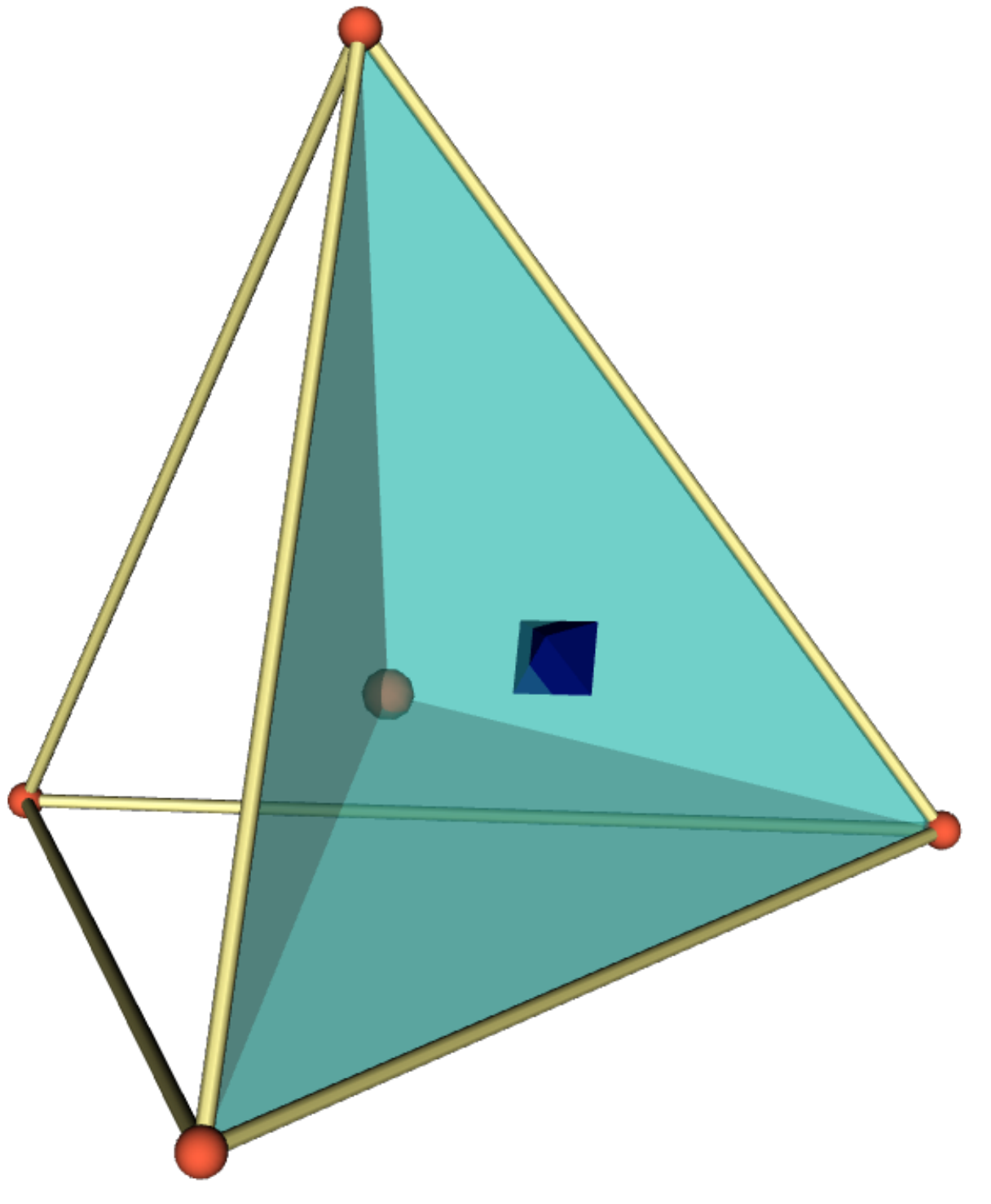}}}
				\put(12,-55){\normalsize{$V_{2}$}}
				\put(104,-13){\normalsize{$V_{3}$}}
				\put(30,87){\normalsize{$V_{1}$}}
				\put(-15,-12){\normalsize{$V_{4}$}}
				\put(30,0){\normalsize{$B$}}
				\put(58,15){\normalsize{$N_i$}}
				\end{picture}
				\vspace{1.5cm}
			\end{center}
		\end{figure}
	\end{minipage}
	\caption{Interior (left) and boundary (right) finite volumes of the face-type in 3D.  The vertex of the primal mesh are denoted by $V_{i}$, $B$ and $B'$ correspond to the barycentres of the primal tetrahedra and $N_{i}$ are the nodes of the dual mesh.
	}\label{vf3d}
\end{figure}

The computational domain is discretized in space using face-based staggered unstructured meshes, as adopted in \cite{BDDV98}, \cite{BFSV14} and \cite{BFTVC17}.
We start considering a triangular (2D) or tetrahedral (3D) grid of elements 
$\{ T_k, \, k=1,\dots ,nel\}$, to be called the \textit{primal mesh}. Then, as illustrated in Figure \ref{vf3d} for the 3D case, the nodes $\left\lbrace N_i\right\rbrace$ of the dual mesh are defined as the barycentres of the faces of the elements of the primal mesh. Then each finite volume of the dual mesh is defined as the polyhedron determined by the vertices of the face, {\color{black} $V_{j}$}, and the barycentres of the two tetrahedra sharing the face, {\color{black} $B,B^{\prime}$}. If the node is on the boundary of the domain, then the associated finite volume is the {\color{black} tetrahedron} defined by the three vertices of the face and the barycentre of the primal {\color{black} tetrahedron} including this face. 

The employed notation is as follows:
\begin{itemize}
	\item Each interior node $N_i$ has as neighbouring nodes the set ${\cal K}_{i}$
	consisting of the barycentres of the edges/faces of the two primal elements to which it belongs.
	\item Each finite volume, also to be called cell, is denoted by $C_{i}$. We denote by $\Gamma_{i}$ its boundary and $\boldsymbol{\widetilde{\eta}}_{i}$ its outward unit normal.
	\item The edge/face $\Gamma_{ij}$ is the face between cells $C_{i}$ and $C_{j}$. $N_{ij}$ is the barycentre of the $\Gamma_{ij}$.
	\item The boundary of $C_{i}$ is denoted by
	$\Gamma_i =\displaystyle \bigcup_{N_j \in {\cal K}_{i}} \Gamma_{ij}$ .
	\item  $\left|C_{i}\right|$ is the area/volume of $C_{i}$.
	\item $\widetilde{\boldsymbol{\eta}}_{ij} $ represents the outward unit normal vector to $\Gamma_{ij}$.
	We define $\boldsymbol{\eta_{ij}}:=\boldsymbol{\widetilde{\eta}_{ij}} ||\boldsymbol{\eta_{ij}} ||$, where,
	$||\boldsymbol{\eta_{ij}} ||$ {\color{black}  denotes} the length/area of $\Gamma_{ij}$. 
\end{itemize}

{\color{black}
	Therefore, the main location of the variables involved in the algorithm is as follows. The conservative variables, $\mathbf{W}^{n}_{i}$, as well as the intermediate approximations $\widetilde{\mathbf{W}}_{i}$, $\widetilde{p}_{p\, i}$, and the temperature, $\theta^{n}_{i}$, are computed on the dual cell $C_{i}$.
	$\widetilde{p}_{\rho}$ is approximated at each primal element $T_{k}$ and the pressure, $\press^{n}_{j}$, is obtained at each primal vertex $V_{j}$. 
}

\subsection{Transport-diffusion stage}
Within the {\color{black}  transport-diffusion} stage a finite volume method is applied in order to compute the density at the new time step, $\rho^{n+1}$, and to provide a first approximation of the linear momentum density, $\widetilde{\mathbf{W}}_{\mathbf{u}}$, and the pressure, $\widetilde{\press}_{\press}$.  
Integrating equations \eqref{eq:disc_mass}, \eqref{eq:disc_mom_tilde} and \eqref{eq:disc_pres_tilde} on $C_i$ and applying Gauss' theorem we get

\begin{gather}
	\frac{\left|C_i\right|}{\Delta t} \left(\rho^{n+1}_{i}-\rho^{n}_{i}\right) + \int_{\Gamma_i} \Flux^{\rho}\left(\mathbf{W}^n\right)\boldsymbol{\widetilde{\eta}}_{i}\,dS = 0, \label{eq:discr_comp_mass}
	\\
	\frac{\left|C_i\right|}{\Delta t}\left( \widetilde{\mathbf{W}}_{\mathbf{u},\, i}^{n+1}-\mathbf{W}_{\mathbf{u},\, i}^{n}\right) +\int_{\Gamma_i}\Flux^{\mathbf{W}_{\mathbf{u}}}\left(\mathbf{W}^n\right)\boldsymbol{\widetilde{\eta}}_{i}\,dS -\int_{\Gamma_i} \btau^n \,\boldsymbol{\widetilde{\eta}}_{i}\, dS = \int_{C_i} \rho^{n}\g dV,\label{eq:discr_comp_momentum}
	\\
	\frac{\left|C_i\right|}{\Delta t}\left( \widetilde{\press}_{\press\, i}-\press^{n}_{i}\right) + \int_{C_i} \mathbf{U}^{n} \cdot \grae \press^{n} dV  + \int_{C_i} \left(\gamma-1\right) \mathrm{div} \mathbf{q}^{n}  dV = 0.
	\label{eq:discr_comp_pressure}
\end{gather}

\subsubsection{Numerical flux}\label{sec:dvectionterm_lmn}
We start by approximating the flux term in \eqref{eq:discr_comp_mass}-\eqref{eq:discr_comp_momentum}. To this end we define the global normal flux  on ${\Gamma}_{i}$ as

\begin{equation}
	{\mathcal{Z} }(\mathbf{W}^{n},  {\boldsymbol{\tilde{\eta}}}_{i} ): ={
		\Flux} (\mathbf{W}^{n})  {\boldsymbol{\tilde{\eta}}}_{i},\quad\mathrm{ with }\quad
	\Flux (\mathbf{W}^{n})= \left(\Flux^{\rho} (\mathbf{W}^{n}) , \Flux^{\mathbf{w}_{\mathbf{u}}} (\mathbf{W}^{n}) \right)^{T}.
\end{equation}
Next, we split
$\Gamma_{i}$ into the cell interfaces $\Gamma_{ij}$, namely

\begin{equation}
	\int_{ \Gamma_{i} } \mathbf{\mathcal{F}}(
	\mathbf{W}^{n})  \boldsymbol{ \widetilde{\eta}}_{i} \,\mathrm{dS}
	=\displaystyle \sum_{N_j \in {\cal K}_{i}} \displaystyle
	\int_{\Gamma_{ij}} {
		\mathcal{Z} }(\mathbf{W}^{n}, {\boldsymbol{\tilde{\eta}}}_{ij} )  \,\mathrm{dS}.
\end{equation}
Then, in order to get a stable discretization, the integral on $\Gamma_{ij}$ is approximated
by an upwind scheme using a numerical flux function $\boldsymbol{\phi}=\left(\phi_{\rho},\boldsymbol{\phi}_{\mathbf{u}}\right)^{T}$.
We use the simple Rusanov scheme (see \cite{Rus62}),

\begin{gather}
	\boldsymbol{\phi} \left(
	\mathbf{W}_{i}^{n}\, ,\mathbf{W}_{j}^{n}\, , \boldsymbol{
		\eta}_{ij}                   \right)
	= \, \displaystyle\frac{1}{2} ({\cal Z}
	(\mathbf{W}_{i}^n\, ,\boldsymbol{\eta}_{ij})+{\cal Z}(\mathbf{W}_{j}^n \, , \boldsymbol{\eta}_{ij}))-\frac{1}{2}
	\alpha^{n}_{RS,\, ij} \left( \Wur_{j}^n-\Wur_{i}^n\right) \label{eq:flux_lmnt}
\end{gather}
with

\begin{equation}\alpha^{n}_{RS,\, ij}=\alpha_{RS}(\mathbf{W}_{i}^n,\mathbf{W}_{j}^n,\boldsymbol{\eta}_{ij}):=\max \left\lbrace 2\left|\mathbf{U}_i^{n}\cdot\boldsymbol{\eta}_{ij}\right|,2\left|\mathbf{U}_j^{n}\cdot\boldsymbol{\eta}_{ij}\right|\right\rbrace  \label{eq:alphars_decoupled1_comp}\end{equation}
the maximum signal speed on the edge and $\Wur=\left(\rho,\mathbf{W}_{\mathbf{u}}^{T}\right)^{T}$ .
Therefore, equations \eqref{eq:discr_comp_mass}-\eqref{eq:discr_comp_momentum} can be rewritten as

\begin{gather}
	\rho_{i}^{n+1}-\rho_{i}^{n} +\frac{\Delta t}{\left|C_i\right|}  \sum_{N_{j}\in\mathcal{K}_{i}}
	\phi_{\rho}\left( \mathbf{W}_{i}^{n},\mathbf{W}_{j}^{n},\boldsymbol{\eta}_{ij}\right) = 0,
	\\
	\frac{1}{\Delta t}\left( \widetilde{\mathbf{W}}_{\mathbf{u},\, i}-\mathbf{W}_{\mathbf{u},\, i}^{n}\right) +\frac{1}{\left|C_i\right|}  \sum_{N_{j}\in\mathcal{K}_{i}}
	\boldsymbol{\phi}_{\mathbf{u}}\left( \mathbf{W}_{i}^{n},\mathbf{W}_{j}^{n},\boldsymbol{\eta}_{ij}\right)
	-\frac{1}{\left|C_i\right|} \sum_{N_{j}\in\mathcal{K}_{i}} \boldsymbol{\varphi}_{\mathbf{u}}\left( \mathbf{U}_{i}^{n},\mathbf{U}_{j}^{n},\boldsymbol{\eta}_{ij}\right) = \frac{1}{\left|C_i\right|} \int_{C_i}\rho^n \g  dV, \label{eq:discretized_momentum_comp}
\end{gather}
where $\boldsymbol{\varphi}_{\mathbf{u}}$ denotes a diffusion flux function corresponding to the viscous stress tensor term to be detailed in Section \ref{sec:viscousterms_lmn}.
Some of the tests to be presented in Section \ref{sec:numericalresults}
account for an extra artificial viscosity term, $c_{\alpha}\in\mathbb{R}$, on the Rusanov flux,

\begin{equation}
	\alpha^{n}_{RS,\, ij} := \max \left\lbrace 2\left|\mathbf{U}_i^{n}\cdot\boldsymbol{\eta}_{ij}\right|,2\left|\mathbf{U}_j^{n}\cdot\boldsymbol{\eta}_{ij}\right|\right\rbrace  + c_{\alpha} \left\|\boldsymbol{\eta}_{ij} \right\|,\label{eq:artvisccoef}
\end{equation}
which is used to improve the stability of the scheme related to the mass conservation equation when the fluid has large density variations with respect to the magnitude of the velocity field, 
{\color{black} \cite{Hug95,HMJ00,GPP11}}.

\subsubsection{CVC Kolgan-type scheme}\label{sec:cvc}
By using the flux function \eqref{eq:flux_lmnt} we would obtain a first order scheme both in space and time. Second order in space can be reached by extending the CVC Kolgan-type scheme presented in \cite{CV12} and \cite{BFTVC17}. 
The   linear momentum density and mass density variables, $\mathbf{W}^n_{\mathbf{u},\,i}, \, \mathbf{W}^n_{\mathbf{u},\,j},\, \rho_{i}^{n},\, \rho_{j}^{n}$, are replaced in the upwind terms by their improved interpolations, $\mathbf{W}^n_{\mathbf{u},\,i\,L}, \, \mathbf{W}^n_{\mathbf{u},\,j\,R},\, \rho_{i\,L}^{n},\, \rho_{j\,R}^{n}$:

\begin{gather}
	\boldsymbol{\phi}_{\mathbf{u}} \left(\mathbf{W}_{i}^{n},\mathbf{W}_{j}^{n},\mathbf{W}_{i\,L}^{n},\mathbf{W}_{j\,R}^{n},\boldsymbol{\eta}_{ij}                   \right) 
	=\displaystyle\frac{1}{2} \left[ {\cal Z}
	(\mathbf{W}_{i}^n,\boldsymbol{\eta}_{ij})+{\cal Z}(\mathbf{W}_{j}^n,\boldsymbol{\eta}_{ij})\right]
	-\frac{1}{2}
	\alpha^{n}_{RS,\, ij} \left( \Wur_{j\,R}^n-\Wur_{i\,L}^n\right) \label{eq:flux_lmn_cvc}
\end{gather}
(see \cite{BFTVC17} for further details on the computation of the improved interpolations).

\subsubsection{LADER methodology}\label{sec:lader}
To achieve a second order scheme in space and time, we extend LADER techniques to solve equations \eqref{eq:discr_comp_mass}-\eqref{eq:discr_comp_momentum}. This method was first proposed in \cite{BFTVC17}
as a modification of ADER methodology (see \cite{TMN01} and \cite{Toro}) and makes its extension to solve the multidimensional problem easier. The main difference between that initial scheme and the one employed in this paper relates {\color{black}  to} the density variation. 
Both linear momentum and density are extrapolated and the mid-point rule is applied. 
We come now to detail the new computations to be performed at each step of the extended method in 2D:

\begin{description}
	\item[Step 1.] \label{step1} ENO-based reconstruction. Reconstruction of the data in terms of first degree polynomials is considered. At each finite volume we define {\color{black} four} polynomials,
	each of them at the neighbourhood of one of the boundary edges. That is, the cell is divided into {\color{black} four} sub-triangles having one of the edges of the finite volume as basis and the cell node as opposite vertex. Focusing on an edge $\Gamma_{ij}$ and on the discretization of a scalar variable, $W$, its two related reconstruction polynomials for the conservative  variable are
	
	\begin{equation}P^{i}_{ij}(N)=W_i+(N-N_{i})\left( \grae\, W\right)^{i}_{ij},\quad P^{j}_{ij}(N)=W_j+(N-N_{j})\left( \grae\, W\right)^{j}_{ij}.\end{equation}
	To avoid spurious oscillations on the solution we add a non-linearity on the scheme by applying an ENO (Essentially Non-Oscillatory) interpolation method. 
	The slopes are adaptively chosen as follows:
	
	\begin{gather*}
		\left( \grae\, W\right)^{i}_{ij}= \left\lbrace
		\begin{array}{lc}
			\left(\grae\, W \right)_{T_{ijL}}, & \textrm{if }\left| \left(\grae\, W \right)_{ T_{ijL}}\cdot \left(N_{ij}-N_{i}\right)\right| \leq \left|\left(\grae\, W \right)_{ T_{ij}}\cdot \left(N_{ij}-N_{i}\right) \right|,\\[8pt]
			\left(\grae\, W \right)_{ T_{ij}}, & \mathrm{otherwise};
		\end{array}
		\right.
	\end{gather*}
	\begin{gather*}
		\left( \grae\, W \right)^{j}_{ij}= \left\lbrace
		\begin{array}{lc}
			\left(\grae\, W\right)_{ T_{ijR}}, & \textrm{if }\left|\left( \grae\, W \right)_{ T_{ijR}}\cdot \left(N_{ij}-N_{j}\right)\right| \leq \left|\left(\grae\, W \right)_{ T_{ij}}\cdot \left(N_{ij}-N_{j}\right)\right|,\\[8pt]
			\left(\grae\, W\right)_{ T_{ij}},  & \mathrm{otherwise}.
		\end{array}
		\right.
	\end{gather*}
	The triangles $T_{ij}$, $T_{ijL}$ and $T_{ijR}$, in which the gradients of the conservative variable are computed, are constructed by connecting the barycentres of the faces of the finite elements (see Figure \ref{fig:up_g_t} for the 2D representation). The computation of these derivatives can also be seen as the gradient computation of Crouzeix-Raviart finite elements defined on the primal mesh.
	
	\begin{figure}
		\centering
		\includegraphics[width=0.3\linewidth]{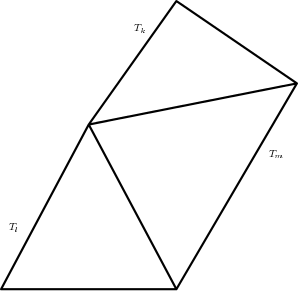}
		\includegraphics[width=0.3\linewidth]{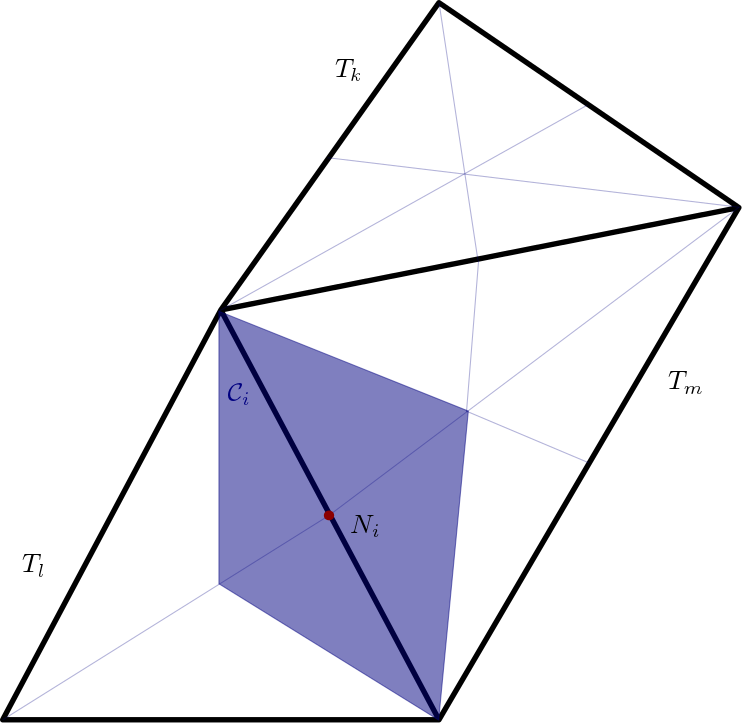}
		\includegraphics[width=0.3\linewidth]{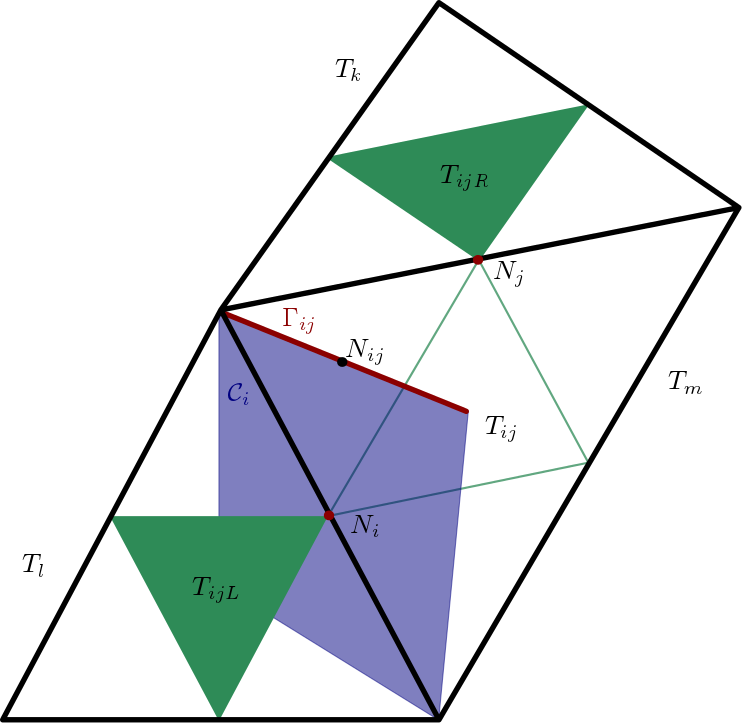}
		\caption{Construction of a dual 2D mesh and auxiliary triangles.
			Left: finite elements of the primal triangular mesh (black).
			Centre: finite volume $C_{i}$ (purple).
			Right: upwind and centred auxiliary triangles (green).}\label{fig:up_g_t}
	\end{figure}
	\item[Step 2.] Computation of the boundary extrapolated values at the barycentre of faces, $N_{ij}$:
	
	\begin{gather} W_{i\, N_{ij}} =p^{i}_{ij}(N_{ij})
		= W_{i}+(N_{ij}-N_{i})   \left( \grae\, W\right)^{i}_{ij},\\
		W_{j\, N_{ij}}=p^{j}_{ij}(N_{ij})=
		W_{j}+(N_{ij}-N_{j})   \left( \grae\, W\right)^{j}_{ij}.  \end{gather}
	
	\item[Step 3.]\label{item:step3} Computation of the variables involved in the flux term with second order of accuracy using the mid-point rule. 
	Taylor series expansion in time and the Cauchy-Kovalevskaya procedure, {\color{black}  based on the momentum conservation equation \eqref{eq:momentum2}}, are
	applied to locally approximate the conservative variables at time $\frac{\Delta t}{2}$.
	This methodology accounts for the contribution of the advection and diffusion terms to the time evolution of the flux term at the momentum conservation equation.
	The resulting evolved variables read
	
	\begin{gather}\overline{\mathbf{W}}_{\mathbf{u}\,i\, N_{ij}}=\mathbf{W}_{\mathbf{u},\,i\, N_{ij}}-\frac{\Delta t}{2\mathcal{L}_{ij}}\left[ \mathcal{Z}(\mathbf{W}_{i\, N_{ij}},\boldsymbol{\eta}_{ij}  )+\mathcal{Z}(\mathbf{W}_{j\, N_{ij}},\boldsymbol{\eta}_{ij}  )\right] \notag\\ +
		\frac{\mu \Delta t}{2\mathcal{L}_{ij}} \left[ \left( \gra\, \mathbf{W}_{\mathbf{u}}+\gra\, \mathbf{W}_{\mathbf{u}}^{T}\right)_{i\, N_{ij}}\boldsymbol{\eta}_{ij}  +  \left( \gra\, \mathbf{W}_{\mathbf{u}}+\gra\, \mathbf{W}_{\mathbf{u}}^{T}\right)_{j\,N_{ij}} \boldsymbol{\eta}_{ij} \phantom{\frac{2}{3}}
		\right. \notag\\ \left.
		-\frac{2}{3} \left( \dive \mathbf{W}_{\mathbf{u},\, i \, N_{ij}}\boldsymbol{ \eta }_{ij} +  \dive  \mathbf{W}_{\mathbf{u},\, j \, N_{ij}}\boldsymbol{ \eta }_{ij} \right)   \right]
		+{\color{black} \frac{ \Delta t}{2}}  \left((\grae\, \press\right)_{T_{ij}}
		, \label{eq:evolved_veli}\\
		\overline{ \mathbf{W}}_{\mathbf{u}\,j\, N_{ij}}= \mathbf{W}_{\mathbf{u},\,j\, N_{ij}}-\frac{\Delta t}{2\mathcal{L}_{ij}}\left[ \mathcal{Z}( \mathbf{W}_{i\, N_{ij}},\boldsymbol{\eta}_{ij}  )+\mathcal{Z}( \mathbf{W}_{j\, N_{ij}},\boldsymbol{\eta}_{ij}  )\right] \notag \\+
		\frac{\mu \Delta t}{2\mathcal{L}_{ij}} \left[ \left( \gra\, \mathbf{W}_{\mathbf{u}}+\gra\, \mathbf{W}_{\mathbf{u}}^{T}\right)_{i\,N_{ij}}\boldsymbol{\eta}_{ij}  +  \left( \gra\, \mathbf{W}_{\mathbf{u}}+\gra\, \mathbf{W}_{\mathbf{u}}^{T}\right)_{j\,N_{ij}} \boldsymbol{\eta}_{ij} \phantom{\frac{2}{3}}
		\right. \notag\\ \left.
		-\frac{2}{3} \left( \dive \mathbf{W}_{\mathbf{u},\, i \, N_{ij}}\boldsymbol{ \eta }_{ij} +  \dive \mathbf{W}_{\mathbf{u},\, j \, N_{ij}}\boldsymbol{ \eta }_{ij} \right)   \right]+{\color{black} \frac{\Delta t}{2}}  \left( (\grae\, \press\right)_{T_{ij}}.\label{eq:evolved_velj}\end{gather}
	We have denoted $\mathcal{L}_{ij}=\min \left\lbrace \frac{\left| C_{i} \right|}{\mathrm{S}(C_i)},
	\frac{\left| C_{j} \right|}{\mathrm{S}(C_j)} \right\rbrace$
	with $\mathrm{S}(C_i)$ the area of the surface of cell $C_i$ and $(\grae\, W)_{i\,N_{ij}}$, $(\grae\, W)_{j\,N_{ij}}$ the approximation of the gradients of $W$ in the neighbourhood of face $\Gamma_{ij}$ in the related finite volumes $C_{i}$ and $C_{j}$, respectively. Let us remark that these derivatives correspond to the values obtained at $T_{ijL}$ and $T_{ijR}$.	
	Regarding density, we consider the mass conservation equation, {\color{black} \eqref{eq:mass2}}, in order to apply Cauchy-Kovalevskaya procedure. Therefore, we define
	
	\begin{align}\overline{\rho}_{i\, N_{ij}}:=\rho_{i\, N_{ij}}-\frac{\Delta t}{2\mathcal{L}_{ij}}\left[ \mathcal{Z}^{\rho}(\mathbf{W}_{i\, N_{ij}},\boldsymbol{\eta}_{ij}  )+\mathcal{Z}^{\rho}(\mathbf{W}_{j\, N_{ij}},\boldsymbol{\eta}_{ij}  )\right],\label{eq:densityevoli}\\
		\overline{ \rho}_{j\, N_{ij}}:= \rho_{j\, N_{ij}}-\frac{\Delta t}{2\mathcal{L}_{ij}}\left[ \mathcal{Z}^{\rho}( \mathbf{W}_{i\, N_{ij}},\boldsymbol{\eta}_{ij}  )+\mathcal{Z}^{\rho}( \mathbf{W}_{j\, N_{ij}},\boldsymbol{\eta}_{ij}  )\right]  .\label{eq:densityevolj}\end{align}
	
	\item[Step 4.] Computation of the numerical flux using \eqref{eq:flux_lmnt}:
	\begin{gather}
		\boldsymbol{\phi} \left(\overline{\mathbf{W}}_{i\, N_{ij}}^{n},\overline{\mathbf{W}}_{j\, N_{ij}}^{n},\boldsymbol{\eta}_{ij}\right)
		= \displaystyle\frac{1}{2} \left[ {\cal Z}
		(\overline{\mathbf{W}}_{i\, N_{ij}}^n,\boldsymbol{\eta}_{ij})+{\cal Z}(\overline{\mathbf{W}}_{j\, N_{ij}}^n,\boldsymbol{\eta}_{ij})\right]
		-\frac{1}{2}
		\overline{\alpha}^{n}_{RS,\, ij} \left( \overline{\Wur}_{j\, N_{ij}}^n-\overline{\Wur}_{ i \, N_{ij}}^n\right).\label{eq:flux_lmnt_lader}
	\end{gather}
\end{description}

\begin{remark}
	In Step 1, an ENO-based reconstruction has been employed in order to introduce a non-linearity which circumvents Godunov theorem. An alternative may be the use of classical limiters like the minmod limiter of Roe, \cite{Roe85,Toro}, or the Barth-Jespersen limiter, \cite{BJ89}.
\end{remark}

\subsubsection{Viscous term}\label{sec:viscousterms_lmn}
In this section we describe the computation of the viscous term. 
First, applying Gauss' theorem we relate the volume integral of the diffusion term with
a surface integral over the boundary, $\Gamma_{i}$. Next, this integral is split into the integrals on the cell faces
$\Gamma_{ij}$. Thus, the viscous term of the momentum conservation equation reads 

\begin{gather}
	\int_{C_i}\dive \btau^n dV = \sum_{N_j\in\mathcal{K}_i}\int_{\Gamma_{ij}} \btau^n\boldsymbol{\widetilde{\eta}}_{ij}  \, \mathrm{dS} 
	= \sum_{N_j\in\mathcal{K}_i}\int_{\Gamma_{ij}}\mu\left[
	\gra\, \mathbf{U}^{n}+\left( \gra\,\mathbf{U}^{n}\right)^T  - \frac{2}{3} \dive  \mathbf{U}^{n} I\right]\boldsymbol{\widetilde{\eta}}_{ij}  \,\mathrm{dS},\label{eq:visc_term_int_comp}
\end{gather}
where a new divergence term
has appeared with respect to the incompressible model, \cite{BFTVC17}.
We define the corresponding numerical diffusion function as

\begin{gather}
	\varphi_{\mathbf{u}}\left( \mathbf{U}_{i}^n,\mathbf{U}_{j}^n, \boldsymbol{\eta}_{ij} \right) \approx \int_{\Gamma_{ij}}   \mu \left[ \gra\, \mathbf{U}^n +\left( \gra\, \mathbf{U}^n\right)^{T} - \frac{2}{3} \dive \mathbf{U}^n \right] \widetilde{\boldsymbol{\eta}}_{ij}\mathrm{dS}.
\end{gather}
Accounting for the dual mesh structure, we compute the spatial derivatives on the auxiliary tetrahedron $T_{ij}$ through a Galerkin approach. Hence,

\begin{gather}
	\boldsymbol{\varphi}_{\mathbf{u}}\left( \mathbf{U}_{i}^n,\mathbf{U}_{j}^n,\boldsymbol{\eta}_{ij} \right)
	=\mu \left( \gra\, \mathbf{U}^n\right)_{T_{ij}}\boldsymbol{\eta}_{ij}
	+\mu \left( \gra\, \mathbf{U}^{n}\right)_{T_{ij}}^{T}\boldsymbol{\eta}_{ij} - \frac{2}{3}\mu\left( \dive  \mathbf{U}^{n}\right)_{T_{ij}}  \boldsymbol{\eta}_{ij}. \label{eq:psi_viscterm_grad_comp}
\end{gather}

The former expressions are directly used when a first order or the CVC Kolgan-type methodologies are considered.
In order to attain a second order in space and time scheme we apply LADER methodology also to the diffusion term (see \cite{BFTVC17}). That is, we construct evolved variables $\overline{\overline{ \mathbf{W}}}_{\mathbf{u},\,i}^{n} $,  which lack from  the advection term as contribution for the half in time evolution. Next, the physical evolved variables are obtained {\color{black} by} dividing by the evolved densities already computed for the flux term,

\begin{equation}
	\overline{\overline{ \mathbf{U}}}_{i}^{n}=\frac{\overline{\overline{ \mathbf{W}}}_{\mathbf{u},\,i}^{n} }{\overline{\rho}_{i}^{n} }.
\end{equation}
Finally, the numerical diffusion function is evaluated, 

\begin{gather}
	\boldsymbol{\varphi}_{\mathbf{u}}\left(\overline{\overline{ \mathbf{U}}}_{i}^{n} ,\overline{\overline{ \mathbf{U}}}_{j}^{n} ,\boldsymbol{\eta}_{ij} \right)	. \label{eq:psi_viscterm_grad_comp_lader}
\end{gather}	

\begin{remark}
	The above evolved variables are also employed to compute the remaining terms in equations \eqref{eq:discr_comp_momentum}-\eqref{eq:discr_comp_pressure} when the second order in space and time scheme is derived. {\color{black} Up to now only the evolved values of linear momentum and density were needed. However,  the pressure appears in several terms of \eqref{eq:discr_comp_momentum}-\eqref{eq:discr_comp_pressure} and its evolved value will also be needed. Employing a Taylor series expansion in time and using equation \eqref{eq:pressure2} within the Cauchy-Kovalevskaya procedure yield}
	
	\begin{equation}
		\overline{\press}^{n}_{i} =  \press^{n}_{i} 
		- \rho^{n}_{i} {\color{black}\left( c^{n}_{i}\right)^2}  \dive  \mathbf{u}^{n} 
		- \mathbf{u}^{n}_{i}\cdot \grae \press^{n} 
		- \left(\gamma-1\right) \dive q^{n}
		+ \left(\gamma-1\right) \boldsymbol{\tau}^{n}\cdot \gra \mathbf{u}^{n},
	\end{equation}
	where the spatial derivatives were computed on the primal elements using the Galerkin approach and then interpolated on the dual mesh.
\end{remark}

\subsubsection{Gravity term}
The gravity term can be integrated directly per finite volume by assuming a constant value for the density on $C_{i}$,

\begin{equation}
	\int_{C_i} \rho^{n}\g  dV = \left|C_i\right| \rho^{n}_{i} \g.
\end{equation}

\subsubsection{Non conservative product}
Following \cite{Par06,CFFP09,DCPT09,DPRZ16,GCD18}, a path conservative scheme is employed to  approximate 
the non conservative product:

\begin{equation}
	\int_{C_i} \mathbf{U}^{n} \cdot \grae \press^{n} dV = \displaystyle \sum_{N_j \in {\cal K}_{i}} \displaystyle
	\frac{\mathbf{W}_{\mathbf{u}\,ij}^{n} }{2\rho_{ij}^{n}} \cdot \boldsymbol{\tilde{\eta}}_{ij} \left( \press_{j}^{n}- \press_{i}^{n}  \right). 
\end{equation}

\subsubsection{Heat flux term}
Assuming the average value of the temperature at each dual element is known, the heat flux term on \eqref{eq:discr_comp_pressure} can be approximated as

\begin{equation}
	\int_{C_i} \left( \gamma -1\right)\mathrm{div} \mathbf{q}^{n} dV = 
	\left( \gamma -1\right)\mathbf{q}_{T_{ij}}^{n}  \cdot {\boldsymbol{\tilde{\eta}}}_{ij}=
	-\left( \gamma -1\right)\lambda \left(\grae\, \theta^{n} \right)_{T_{ij}} \cdot {\boldsymbol{\tilde{\eta}}}_{ij}.
\end{equation}

\subsection{Pre-projection stage: density term in the pressure equation}\label{sec:densinpress}
Let us remark that the final objective of solving the pressure equations \eqref{eq:disc_pres_tilde}-\eqref{eq:disc_pres_tilder} is to obtain a value for $\widetilde{\press}$ that will be replaced in \eqref{eq:disc_pres_n1}. 
Therefore, for the density derivative term, instead of computing the integral at each finite volume and then passing the information to the primal mesh, we propose to directly compute the value at each primal element in a finite volume fashion, i.e., by considering the sum of the integrals on its boundary:

\begin{equation}
	\widetilde{\press}_{\rho}= -\frac{\Delta t}{\left|T_k\right|} \int_{T_k} {\color{black} \left( c^{n}\right)^{2}} \mathbf{U}^{n} \cdot \grae \rho^{n} dV = 
	-\frac{\Delta t}{\left|T_k\right|} {\color{black}\left( c^{n}_{k}\right)^2} \, \mathbf{U}^{n}_{k}\cdot \sum_{N_{i}\in T_{k}} \int_{\Gamma_{ki}}  \rho^{n}_{i}  \boldsymbol{\tilde{\eta}}_{ki} \,dV. \label{eq:densitytermapprox}
\end{equation}
In \eqref{eq:densitytermapprox} the sound speed is computed from the value of the density at each dual element and the pressure in the vertex of the primal grid:

\begin{equation}
	c^{2}_{k} = \sum_{i\in \mathcal{K}_{k}} \frac{\gamma \press_{ki}}{\rho_{i}} \frac{\left|T_{ki}\right|}{\left| T_{k}\right|}, \quad \press_{ki}= \frac{1}{{\color{black}\ell}}\sum_{m=1}^{{\color{black}\ell}} \press_{km}
\end{equation}
with $\ell$ the number of local vertex {\color{black} per primal element face}, i.e., {\color{black} $\ell=2$} in 2D and {\color{black} $\ell=3$} in 3D,  {\color{black} and $T_{ki}$ the half dual element with basis the face between the primal elements $k$ and $i, i\in \mathcal{K}_{k}$, and opposite vertex the barycentre of the primal element $T_{k}$}.
Similarly, the velocity has been interpolated from the dual to the primal mesh as
\begin{equation}
	\mathbf{U}_{k} = \sum_{i\in \mathcal{K}_{k}} \mathbf{U}_{i} \frac{\left|T_{ki}\right|}{\left| T_{k}\right|}. \label{eq:interpolation_dual2primal}
\end{equation}

\subsection{Projection stage}
Within the projection stage, equations \eqref{eq:disc_mom_n1}, \eqref{eq:disc_pres_n1} are solved using a finite element method. To obtain the weak formulation of the pressure system, we multiply equation \eqref{eq:disc_mom_n1} by the gradient of a test function $z\in V_{0}$, $V_{0}:=\left\lbrace z\in H^{1}\left(\Omega\right): \int_{\Omega} z dV = 0\right\rbrace$ and  integrate in $\Omega$:
\begin{eqnarray}
	\int_{\Omega} \grae \press^{n+1} \cdot \grae z dV =
	\frac{1}{\Delta t} \int_{\Omega} \widetilde{\mathbf{W}}_{\mathbf{u}} \cdot \grae z dV
	-\frac{1}{\Delta t} \int_{\Omega} \mathbf{W}_{\mathbf{u}}^{n+1}\cdot \grae z dV. \label{eq:varform}
\end{eqnarray} 
Besides, applying a Green's formula in \eqref{eq:disc_pres_n1}, we get
\begin{equation}
	\int_{\Omega} \mathbf{W}_{\mathbf{u}}^{n+1}\cdot \grae z dV = \int_{\Gamma} \mathbf{W}_{\mathbf{u}}^{n+1} \cdot \boldsymbol{\eta} z dS + \frac{1}{c^{2}\Delta t} \int_{\Omega} \left(\press^{n+1}-\widetilde{\press}\right) z dV
	-\frac{1}{c^{2}}\int_{\Omega}  \left(\gamma-1\right) \boldsymbol{\tau}^{n}\cdot \gra \mathbf{U}^{n}\, z dV
	. \label{eq:press_green}
\end{equation} 
Replacing \eqref{eq:press_green} in the variational formulation, \eqref{eq:varform}, we obtain the following weak problem:

\begin{weakproblem}
	Find  $\press^{n+1}\in V_{0}$ satisfying
	\begin{gather}
		\Delta t^{2}\int_{\Omega} \grae \press^{n+1} \cdot \grae z dV + \frac{1}{c^{2}} \int_{\Omega} \press^{n+1} z dV= 
		\Delta t \int_{\Omega} \widetilde{\mathbf{W}}_{\mathbf{u}} \cdot \grae z dV
		-\Delta t\int_{\Gamma} \mathbf{W}_{\mathbf{u}}^{n+1} \cdot \boldsymbol{\eta} z dS \notag\\
		+ \frac{1}{c^{2}} \int_{\Omega}\widetilde{\press} z dV
		-\frac{\Delta t}{c^{2}}\int_{\Omega}  \left(\gamma-1\right) \boldsymbol{\tau}^{n}\cdot \gra \mathbf{U}^{n}\, z dV
		\label{eq:weak_problem}
	\end{gather} 
	for all $z\in V_{0}$.
\end{weakproblem}
Since the contribution of the density term on the pressure equation has already been computed at each primal element whereas the remaining terms were approximated on the dual elements, the final contribution to $\widetilde{\press}$ is made of two different parts,

\begin{equation}
	\frac{1}{c^{2}} \int_{\Omega}\widetilde{\press} z dV= \frac{1}{c^{2}} \int_{\Omega}\widetilde{\press}_{\press} z dV +  \frac{1}{c^{2}} \int_{\Omega}\widetilde{\press}_{\rho} z dV.
\end{equation}
Interpolation from the dual mesh to the primal mesh of $\widetilde{\press}_{\press}$ and $\widetilde{\mathbf{W}}_{\mathbf{u}}$
is done within the pre-projection stage according to \eqref{eq:interpolation_dual2primal}.
The discretization of \eqref{eq:weak_problem} is performed using classical  $\mathbb{P}_{1}$  finite elements and the final system is solved using an optimized conjugate gradient method.

An alternative formulation to \eqref{eq:weak_problem} consists in applying Green's formula to the first term of the right hand side,

\begin{equation}
	\int_{\Omega} \widetilde{\mathbf{W}}_{\mathbf{u}} \cdot \grae z dV = - \int_{\Omega}\dive \widetilde{\mathbf{W}}_{\mathbf{u}} \cdot z dV + \int_{\Gamma} \widetilde{\mathbf{W}}_{\mathbf{u}} \cdot \boldsymbol{\eta} z dS.
\end{equation}
To approximate the integral related to the divergence of the linear momentum we proceed similarly to the computation of the density term in the pressure equation, Section \ref{sec:densinpress},

\begin{equation}
	\int_{T_k}  \dive \widetilde{\mathbf{W}}_{\mathbf{u}} dV = 
	\sum_{N_{i}\in T_{k}} \int_{\Gamma_{ki}} \widetilde{\mathbf{W}}_{\mathbf{u}\, i} \cdot  \boldsymbol{\tilde{\eta}}_{ki} \,dV. \label{eq:wtilde_tnt}
\end{equation}
Thus, we obtain a constant approximation of the linear momentum divergence at each primal element. 

The last term to be computed is the viscosity term in \eqref{eq:weak_problem}. Using the Galerkin approach we can get a constant approximation of the gradients of the velocity at each primal element. Then, the viscous stress tensor is computed and multiplied by the velocity gradient. Multiplication by the test function and integration on each simplex element provides the contribution of the term to the right hand side of equation \eqref{eq:weak_problem}.
The final resulting system is solved using a conjugate gradient method.

\subsection{Post-projection stage}

Finally, at the post-projection stage, $\mathbf{W}^{n+1}_{\mathbf{u}}$ is updated with the pressure contribution,

\begin{equation}
	\mathbf{W}_{\mathbf{u}}^{n+1}=\widetilde{\mathbf{W}}_{\mathbf{u}} -\Delta t\grae \press^{n+1}.\label{eq:velocity_update_comp}
\end{equation}
The pressure gradient involved in the previous equation is computed using the classical finite element gradients for $\mathbb{P}_{1}$ which provides a value for each finite element. Then, it is transferred into the dual mesh by considering the weighted average of the contributions of the two halves of each cell.

As a needed post-process, when considering a non-zero heat flux, we approximate the temperature using the state equation,

\begin{equation}
	\theta^{n+1}_{i} = \frac{\press^{n+1}_{i}}{\rho_{i}^{n+1} R}.
\end{equation}

\subsection{Boundary conditions}
Definition of boundary conditions can be reduced to diverse combinations of the following types:
\begin{itemize}
	\item Periodic boundary conditions. They are {\color{black} built} on the assumption that a periodic mesh is provided. Dual elements on periodic boundaries are constructed by joining the two boundary volumes that share the common face, so that a new dual element of the interior type is generated. Regarding the FE computations, the vertices on the boundaries are merged resulting in a reduction on the size of the pressure system.
	
	\item Dirichlet boundary conditions. 
	They can be imposed strongly or weakly. In the first case the velocity and the density/temperature at the boundary are overwritten using their exact values. In the second case, the exact velocity and the density/temperature on the boundary are imposed to construct the fluxes and the gradients involved in the viscous terms. 
	Dirichlet boundary conditions for the linear momentum correspond to Neumann boundary conditions for the pressure field. 
	
	\item Adiabatic wall. It corresponds to Dirichlet boundary conditions but for the computation of the heat flux which is set to zero. Moreover, the density at the boundary is not imposed.
	
	\item Neumann boundary conditions. The definition of $\widetilde{\mathbf{W}}_{\mathbf{u}}$ takes into account inflow/outflow boundary conditions with no need for additional treatment of the conservative variables. Then, the pressure is imposed on the boundary nodes of the primal mesh.
\end{itemize}

\section{Numerical results}\label{sec:numericalresults}
In this section, classical benchmarks for weakly compressible flows are presented in order to assess the performance of the proposed methodology. Let us note that the international system of units (SI) is considered for all the tests. {\color{black} The time step size $\Delta t$ is computed according to the CFL condition based on the flow velocity and the kinematic viscosity as follows 
	\begin{equation}
		\Delta t = \mathrm{CFL}_{\mathbf{u}}\, \cdot \min_{i} \frac{ r_{i}^{2} }{ \left\| \mathbf{u}_i \right\| r_{i} + \left|\nu_{i}\right|_{\mathrm{max}}},
	\end{equation}
	where $r_{i}$ denotes the incircle radius of the dual volume $C_{i}$, while $\left\|\mathbf{u}_{i}\right\|$ and  $\left|\nu_{i}\right|_{\mathrm{max}}$ are the maximum eigenvalues in absolute value on the dual cell associated with the convective and the diffusive terms, respectively, which in our semi-implicit scheme are both discretized explicitly. 
	The maximum CFL number based on the sound speed is denoted by $\mathrm{CFL}_{c}$ and is defined as  
	\begin{equation}
		\mathrm{CFL}_{c} = \max_{i} \left( c_i \,  \frac{\Delta t}{r_{i}} \right).
	\end{equation}
	It is reported for all the tests in order to show the benefits of a semi-implicit scheme over a classical explicit Godunov-type finite volume method.
}
{\color{black} All 2D simulations have been run on an Intel$^{\textrm{\textregistered}}$ Core\texttrademark   i$7-4720$HQ CPU @$2.60$GHz processor, whereas an Intel$^{\textrm{\textregistered}}$  Xeon$^{\textrm{\textregistered}}$ Gold 6140M processor has been used for the 3D tests.}

\subsection{Taylor-Green vortex}\label{sec:TGV}
To check the accuracy of the numerical method, we employ the Taylor Green vortex benchmark defined in $\Omega=\left[0,2\pi\right]\times\left[0,2\pi\right]$. An exact solution of this test case with gravity source terms can be defined as
\begin{equation}
	\rho\left(x,y,t\right) = 1,\qquad
	\press \left(x,y,t\right) = \frac{\press_{0}}{\gamma-1} + \frac{1}{4} \left(\cos(2x)+\cos(2y) \right)e^{-4\mu t} , \\
	{u}_{1} \left(x,y,t\right) = \sin(x)\cos(y)e^{-2 \mu t}, \qquad
	{u}_{2} \left(x,y,t\right) = -\cos(x)\sin(y)e^{-2\mu t} + gt,
\end{equation}
where $\press_{0}=10^{5}$, $\gamma=1.4$ {\color{black} and the characteristic Mach number is $M\approx1.7 \cdot10^{-3}$}. In order to \eqref{eq:mass2}-\eqref{eq:pressure2} verify the former analytical solution we need to impose the following source terms

\begin{gather}
	f_{\rho}\left(x,y,t\right)=0,\\
	{f}_{u_1}\left(x,y,t\right)=- g t e^{-2 t \mu}\sin(x) \sin(y),\\
	{f}_{u_2}\left(x,y,t\right)=-g te^{-2 t \mu}\cos(x)\cos(y),\\
	{f}_{\press}\left(x,y,t\right)=- 2 \mu e^{-4 t \mu}\left(\cos^2(x) + \cos^2(y) - 1\right) - e^{-6 t \mu}\cos(x)\cos(y)\sin^{2}(x)  \notag\\   - e^{-4 t\mu}\cos(y)\sin(y) \left(gt - e^{-2t\mu}\cos(x)\sin(y)\right)
\end{gather}
{\color{black} that are discretized in the finite volume framework using a fourth order quadrature rule as introduced in \cite{BFSV14} and \cite{BFTVC17}.} {\color{black} The time step has been determined following the CFL condition based on the flow velocity with $\mathrm{CFL}_{\mathbf{u}}=0.5$ which corresponds to a maximum $\mathrm{CFL}$ based on the sound velocity of $\mathrm{CFL}_{{c}}= 294.1$}.
Two different test cases have been run on the four meshes described in Table \ref{TGV_mesh}.  The first of them, T1, corresponds with the classical steady state Taylor-Green vortex benchmark, with $g=0$, $\mu=0$. Meanwhile, in T2 we have set $g=-9.81$ and $\mu=0.1$ so the flow is no longer stationary. The $L_{2}$ error norms and the corresponding convergence rates are presented in Table~\ref{TGV_errors}, where 

\begin{gather}
	E(W)_{M_i} = \left\|W-W_{M_i} \right\|_{l^2(L^2(\Omega))},\quad
	o_{W_{M_i/M_j}} = \frac{\log\left( E(W)_{M_i}/E(W)_{M_j}\right) }{\log\left( h_{M_i}/h_{M_j}\right) }\end{gather}
for any scalar variable $W$ and $h_{M_{i}}$ the minimum area of the finite volumes on $M_{i}$.
We observe that the first order scheme is slightly below the expected accuracy for the pressure variable on the first test case, but the LADER methodology overcomes that issue and successfully achieves the second order of accuracy sought. Second order is also attained for the second test case, where gravity and viscous terms play an important role.
\begin{table}[H]
	\renewcommand{\arraystretch}{1.2}
	\begin{center}
		\begin{tabular}{|c||c|c|c|}\hline
			Mesh & Elements & Vertices & Dual elements \\\hline\hline
			$M_1$ & $128 $ & $81 $ & $208 $ \\
			$M_2$ & $512 $ & $289 $ & $800 $ \\
			$M_3$ & $2048 $ & $1089 $ & $3136 $ \\
			$M_4$ & $8192 $ & $4225 $ & $12416 $ \\\hline
		\end{tabular}
		\caption{Taylor-Green vortex. Mesh features. }\label{TGV_mesh}
	\end{center}
\end{table}
\begin{table}[H]
	\renewcommand{\arraystretch}{1.2}
	\begin{center}
		\begin{tabular}{|c|c||cccc|ccc|}
			\hline Test/Method & Variable & $E_{M_1}$ & $E_{M_2}$ & $E_{M_3}$  & $E_{M_4}$ & $o_{M_1/M_2}$ & $o_{M_2/M_3}$ & $o_{M_3/M_4}$\\\hline\hline
			\multirow{2}{*}{T1 / Order 1} &$\press$                       
			&$4.51E+02$ & $4.21E+02$ & $3.08E+02$ & $1.91E+02$ & $0.10$ & $0.45$ & $0.69$ \\
			&$\mathbf{w}_{\mathbf{u}}$   
			&$8.20E-02$ & $4.49E-02$ & $2.45E-02$ & $1.29E-02$ & $0.87$ & $0.88$ & $0.92$ \\\hline
			\multirow{2}{*}{T1 / LADER} &$\press$                       
			& $3.85E+01$ & $4.06E+00$ & $3.11E-01$ & $2.24E-02$ & $3.25$ & $3.71$ & $3.80$ \\
			&$\mathbf{w}_{\mathbf{u}}$   
			& $3.29E-02$ & $8.93E-03$ & $2.34E-03$ & $6.06E-04$ & $1.88$ & $1.94$ & $1.94$ \\\hline
			\multirow{2}{*}{T2 / LADER} 
			&$\press$                       
			& $1.22E+01$ & $1.66E+00$ & $3.36E-01$ & $7.54E-02$ & $2.88$ & $2.31$ & $2.16$ \\
			&$\mathbf{w}_{\mathbf{u}}$   
			& $4.62E-02$ & $9.59E-03$ & $2.22E-03$ & $5.87E-04$ & $2.27$ & $2.11$ & $1.92$ \\\hline
		\end{tabular}
		\caption{Taylor-Green vortex. Observed $L_{2}$ errors in space and time, $E_{M_i}$, and convergence rates , $o_{M_{i}/M_{i+1}}$ {\color{black} ($\mathrm{CFL}_{\mathbf{u}}=0.5$; $c_{\alpha}=1$ for T1 and $c_{\alpha}=5$ for T2)}. }\label{TGV_errors}
	\end{center}
\end{table}

{\color{black}
	The simulations have been carried out in serial on one single core of an Intel$^{\textrm{\textregistered}}$ Core$^{\textrm{\texttrademark}}$  i$7-4720$HQ CPU @$260$GHz processor.
	Table \ref{tab:TGV_times_errors} contains the computational time (wall clock time), the number of time steps needed to reach the final simulation time and the computational time per dual element and time step,
	\begin{equation}
		t_{e} = \frac{\textrm{CPU time}}{\textrm{ N. Dual elements } \cdot \textrm{N. time steps} },
	\end{equation}
	which gives an independent statement, but for the processor used, of the cost of each numerical method per element and time step.
	The computational cost of solving T1 using the weakly compressible flow solver can be compared with the one of the fully incompressible version of the same code. The results reported in Tables \ref{tab:TGV_times_errors_incompressible} and \ref{tab:TGV_errors_incompressible} show close CPU times for similar errors of the linear momentum variable. 
}

\begin{table}[H]
	\renewcommand{\arraystretch}{1.2}
	\begin{center}
		\begin{tabular}{|c|c||cccc|}
			\hhline{~~|----|} 
			\multicolumn{2}{c|}{ } & $M_1$ & $M_2$ & $M_3$  & $M_4$\\\hhline{--====}
			\multirow{3}{*}{T1 / Order 1} & CPU time (s)&
			$0.03$ & $0.05$ & $0.45$ & $3.00 $
			\\\hhline{|~|-||----|}
			& $t_{e}$ ($\mu$s) &			
			$30.0$ & $6.51$ & $7.60$ & $6.53$
			\\\hhline{|~|-||----|}
			& Time steps &
			$5$ & $9$ & $19$ & $37 $
			\\\hline
			\multirow{3}{*}{T1 / LADER} & CPU time (s)&
			$0.3$ & $0.08$ & $0.61$ & $3.70$
			\\\hhline{|~|-||----|}
			& $t_{e}$ ($\mu$s) &			
			$30.0$ & $11.0$ & $10.2$ & $8.06 $		
			\\\hhline{|~|-||----|}
			& Time steps &
			$5$ & $9$ & $19$ & $37 $
			\\\hline
			\multirow{3}{*}{T2 / LADER} & CPU time (s)&
			$0.06$ & $0.59$ & $4.81$ & $59.1$
			\\\hhline{|~|-||----|}
			& $t_{e}$ ($\mu$s) &			
			$13.7$ & $13.0$ & $9.03$ & $8.48$
			\\\hhline{|~|-||----|}
			& Time steps &
			$22$ & $57$ & $170$ & $561$
			\\\hline
		\end{tabular}
		{\color{black} \caption{Taylor-Green vortex.
				CPU time, CPU time per element and iteration, $t_{e}$, and number of time steps for the weakly compressible Navier-Stokes solver.
			}\label{tab:TGV_times_errors}}
	\end{center}
\end{table}

\begin{table}[H]
	\renewcommand{\arraystretch}{1.2}
	\begin{center}
		\begin{tabular}{|c|c||cccc|}
			\hhline{~~|----|} 
			\multicolumn{2}{c|}{ } & $M_1$ & $M_2$ & $M_3$  & $M_4$\\\hhline{--====}
			\multirow{3}{*}{T1 / Order 1} & CPU time (s)&
			$0.02$ & $0.08$ & $0.45$ & $2.89$
			\\\hhline{|~|-||----|}
			& $t_{e}$ ($\mu$s) &			
			$15.0$ & $9.77$ & $7.60$ & $6.13$
			\\\hhline{|~|-||----|}
			& Time steps &
			$5$ & $10$ & $19$ & $38$
			\\\hline
			\multirow{3}{*}{T1 / LADER} & CPU time (s)&
			$0.03$ & $0.08$ & $0.50$ & $3.36$
			\\\hhline{|~|-||----|}
			& $t_{e}$ ($\mu$s) &			
			$30.0$ & $9.77$ & $7.97$ & $7.12$		
			\\\hhline{|~|-||----|}
			& Time steps &
			$5$ & $10$ & $20$ & $38$
			\\\hline
		\end{tabular}
		{\color{black} \caption{Taylor-Green vortex.
				CPU time, CPU time per element and iteration, $t_{e}$, and number of time steps for the incompressible Navier-Stokes solver presented in \cite{BFTVC17}.
			}\label{tab:TGV_times_errors_incompressible}}
	\end{center}
\end{table}

\begin{table}[H]
	\renewcommand{\arraystretch}{1.2}
	\begin{center}
		\begin{tabular}{|c|c||cccc|ccc|}
			\hline Test/Method & Variable & $E_{M_1}$ & $E_{M_2}$ & $E_{M_3}$  & $E_{M_4}$ & $o_{M_1/M_2}$ & $o_{M_2/M_3}$ & $o_{M_3/M_4}$\\\hline\hline
			\multirow{2}{*}{T1 / Order 1} &$\press$                       
			&$1.12E-01$ & $4.41E-02$ & $2.27E-02$ & $1.20E-02$ & $1.34$ & $0.96$ & $0.93 $ \\
			&$\mathbf{w}_{\mathbf{u}}$   
			&$6.71E-02$ & $4.28E-02$ & $2.76E-02$ & $1.70E-02$ & $0.65$ & $0.63$ & $0.70 $ \\\hline
			\multirow{2}{*}{T1 / LADER} &$\press$                       
			& $1.26E-01$ & $3.68E-02$ & $9.48E-03$ & $2.39E-03$ & $1.78$ & $1.96$ & $1.99 $ \\
			&$\mathbf{w}_{\mathbf{u}}$   
			& $2.99E-02$ & $8.30E-03$ & $2.22E-03$ & $5.76E-04$ & $1.85$ & $1.90$ & $1.95 $ \\\hline
		\end{tabular}
		{\color{black}\caption{Taylor-Green vortex. Observed $L_{2}$ errors in space and time, $E_{M_i}$, and convergence rates , $o_{M_{i}/M_{i+1}}$ for the incompressible Navier-Stokes solver presented in \cite{BFTVC17} ($\mathrm{CFL}_{\mathbf{u}}=0.5$). }\label{tab:TGV_errors_incompressible}}
	\end{center}
\end{table}

{\color{black} As a showcase of the possible speed up in terms of CPU time of the proposed semi-implicit method with respect to explicit algorithms, we finally also consider a fully explicit density-based Godunov-type finite volume scheme solving the compressible Navier-Stokes system \eqref{eq:mass1}-\eqref{eq:totenerg1}. Table \ref{tab:TGV_times_errors_explicit} reports the errors and time consumption for a modified version of T1 in which we have set $p_{0}=10^{4}$ so that the Mach number is greater than in the previous test cases, $M=0.01$, and thus favourable for the explicit algorithm. We observe that the number of time steps for the explicit simulation is much larger than for the semi-implicit method, as expected, since the CFL stability condition now also depends on the sound speed in the medium. Moreover, from the obtained results we can conclude that the CPU time per element and time step is not greatly increased due to the resolution of the pressure system in the semi-implicit scheme. Consequently, the proposed hybrid FV/FE method appears to be computationally more efficient than the fully explicit scheme in the low Mach number limit, as expected. It also avoids the issues related to the wrong scaling of the numerical dissipation in terms of the Mach number that usually arises with classical explicit density-based Godunov-type finite volume schemes without preconditioning, resulting in inaccurate numerical solutions. On the other hand, the weakly compressible semi-implicit scheme proposed in this paper is not much more expensive than a fully incompressible Navier-Stokes solver. 
	For a fair comparison, all tests have been carried out on the same mesh, on the same computer and within the same code basis.
	From these results we can conclude that the proposed hybrid FV/FE method for weakly compressible flows is computationally efficient, both, compared to incompressible flow solvers and explicit density-based Godunov-type finite volume schemes.  }

\begin{table}[H]
	\renewcommand{\arraystretch}{1.2}
	\begin{center}
		\begin{tabular}{|c||ccccc|}
			\hhline{-||-----|} 
			Method & $E_{M_4}(p)$ & $E_{M_4}(\mathbf{w}_{\mathbf{u}})$ & CPU time (s) &  $t_{e}$ ($\mu$s)  & Time steps\\\hhline{======}
			Semi-implicit, incompressible &
			$2.91E-03$ & $6.64E-04$ & $2.58$ & $10.93$ & $19$ \\\hhline{|-||-----|}
			Semi-implicit, weakly compressible &			
			$2.96E-03$ & $6.77E-04 $ & $3.36 $ & $14.24 $ &  $19$		
			\\\hhline{|-||-----|}
			Explicit, compressible &
			$9.58E+02$ & $9.33E+00$ & $262.15$ & $12.09$ & $1746$
			\\\hline
		\end{tabular}
		{\color{black} \caption{Taylor-Green vortex.
				Errors, CPU time, CPU time per element and iteration, $t_{e}$, and number of time steps for test T3 run on mesh $M_{4}$ using the semi-implicit codes and the fully compressible FV Navier-Stokes solver.
			}\label{tab:TGV_times_errors_explicit}}
	\end{center}
\end{table}

{\color{black}	Let us remark that in the semi-implicit simulations we have employed a non zero artificial viscosity coefficient. 
	To better analyse the influence of the new term on the incompressible limit, T1 has been run for different values of the coefficient, $c_{\alpha} \in \left\lbrace 0, 0.1, 1, 10, 100 \right\rbrace$. 
	In Table \ref{tab:TGV_errors_calpha} it can be observed that even large values of $c_{\alpha}$ do not substantially affect the convergence of the numerical method. Moreover, setting $c_{\alpha}=10$ provides  errors of the pressure field comparable to the ones obtained with the incompressible code. 
	Nevertheless, this artificial viscosity must be taken into account on the computation of the time step leading to smaller time steps as its value increases, see Table \ref{tab:TGV_times_calpha}. 
}

\begin{table}[H]
	\renewcommand{\arraystretch}{1.2}
	\begin{center}
		\begin{tabular}{|c|c||cccc|ccc|}
			\hline$c_{\alpha}$ & Variable & $E_{M_1}$ & $E_{M_2}$ & $E_{M_3}$  & $E_{M_4}$ & $o_{M_1/M_2}$ & $o_{M_2/M_3}$ & $o_{M_3/M_4}$\\\hline\hline
			\multirow{2}{*}{$0$} &$\press$                       
			&$4.89E+01$ & $1.14E+01$ & $1.52E+00$ & $1.42E-01$ & $2.10 $ & $2.91$ & $3.42$ \\
			&$\mathbf{w}_{\mathbf{u}}$   
			&$3.09E-02$ & $8.16E-03$ & $2.29E-03$ & $6.54E-04$ & $1.92 $ & $1.84$ & $1.80$ \\\hline
			\multirow{2}{*}{$0.1$} &$\press$                       
			& $4.53E+01$ & $1.04E+01$ & $1.22E+00$ & $1.01E-01$ & $2.12 $ & $3.09$ & $3.60$ \\
			&$\mathbf{w}_{\mathbf{u}}$   
			& $3.04E-02$ & $8.27E-03$ & $2.29E-03$ & $6.42E-04$ & $1.88 $ & $1.85$ & $1.83$ \\\hline
			\multirow{2}{*}{$10$}
			&$\press$                       
			& $4.02E+00$ & $2.49E-01$ & $2.20E-02$ & $3.09E-03$ & $4.01 $ & $3.51$ & $2.83$ \\
			&$\mathbf{w}_{\mathbf{u}}$   
			& $7.35E-02$ & $1.26E-02$ & $2.60E-03$ & $6.10E-04$ & $2.54 $ & $2.28$ & $2.09$ \\\hline
			\multirow{2}{*}{$100$}
			&$\press$                       
			& $2.60E-01$ & $2.92E-02$ & $5.77E-03$ & $1.94E-03$ & $3.15 $ & $2.34$ & $1.57$ \\
			&$\mathbf{w}_{\mathbf{u}}$   
			& $4.60E-01$ & $8.02E-02$ & $1.07E-02$ & $1.44E-03$ & $2.52 $ & $2.91$ & $2.90$ \\\hline
		\end{tabular}
		{\color{black} \caption{Taylor-Green vortex. Observed $L_{2}$ errors in space and time, $E_{M_i}$, and convergence rates , $o_{M_{i}/M_{i+1}}$, for T1 using LADER with $c_{\alpha} \in \left\lbrace 0, 0.1, 10, 100 \right\rbrace$, $\mathrm{CFL}_{\mathbf{u}}=0.5$. }\label{tab:TGV_errors_calpha}}
	\end{center}
\end{table}

\begin{table}[H]
	\renewcommand{\arraystretch}{1.2}
	\begin{center}
		\begin{tabular}{|c|c||cccc|}
			\hhline{~~|----|} 
			\multicolumn{2}{c|}{ } & $M_1$ & $M_2$ & $M_3$  & $M_4$\\\hhline{--====}
			\multirow{3}{*}{$c_{\alpha}=0$} & CPU time (s)&
			$0.02$ & $	0.08$ & $	0.36$ & $	1.92$
			\\\hhline{|~|-||----|}
			& $t_{e}$ ($\mu$s) &			
			$3.76$ & $	9.77$ & $	22.9$ & $	8.15$
			\\\hhline{|~|-||----|}
			& Time steps &
			$2$ & $	10$ & $	5$ & $	19$
			\\\hline
			\multirow{3}{*}{$c_{\alpha}=0.1$} & CPU time (s)&
			$0.02$ & $	0.05$ & $	0.34$ & $	2.13$
			\\\hhline{|~|-||----|}
			& $t_{e}$ ($\mu$s) &			
			$25.0$ & $	11.7$ & $	11.0$ & $	8.56$		
			\\\hhline{|~|-||----|}
			& Time steps &
			$3$ & $	5$ & $	10$ & $	20$
			\\\hline
			\multirow{3}{*}{$c_{\alpha}=10$} & CPU time (s)&
			$0.09$ & $	0.41$ & $	2.70$ & $	18.7$
			\\\hhline{|~|-||----|}
			& $t_{e}$ ($\mu$s) &			
			$18.0$ & $	10.2$ & $	8.53$ & $	7.48$
			\\\hhline{|~|-||----|}
			& Time steps &
			$25$ & $	50$ & $	101$ & $	201$
			\\\hline
			\multirow{3}{*}{$c_{\alpha}=100$} & CPU time (s)&
			$0.47$ & $	3.16$ & $	21.5$ & $	165$
			\\\hhline{|~|-||----|}
			& $t_{e}$ ($\mu$s) &			
			$9.80$ & $	4.93$ & $	7.45$ & $	7.22$
			\\\hhline{|~|-||----|}
			& Time steps &
			$230$ & $	800$ & $	921$ & $	1841$
			\\\hline
		\end{tabular}
		{\color{black} \caption{Taylor-Green vortex.
				Comparison of CPU time, CPU time per element and iteration, $t_{e}$, and number of time steps, for T1 using LADER scheme with $c_{\alpha} \in \left\lbrace 0, 0.1, 10, 100 \right\rbrace$, $\mathrm{CFL}_{\mathbf{u}}=0.5$. 
			}\label{tab:TGV_times_calpha}}
	\end{center}
\end{table}

\subsection{Riemann problems}\label{sec:RP}
A careful study of the performance of the proposed numerical method in the presence of weak discontinuities is carried out using several Riemann problems, see \cite{Toro,TD17}, for Euler flows. We stress that the proposed algorithm is designed for
\textit{weakly compressible} flows and relies on the non-conservative pressure evolution equation, hence we cannot expect the
method to work for high Mach numbers and strong shocks. Nevertheless, the mass and momentum equations are discretized in conservative form and therefore we still expect the method to work for Mach numbers up to unity. 
We consider a rectangular domain with $x\in \left[-0.5,0.5\right]$. Aiming at decreasing the computational cost of the simulations the width of the domain depends on the size of the used elements. More precisely, two different grids, M1 and M2, with $150$ and $500$ divisions along $x$-direction, are employed. 
Periodic boundary conditions are set on $y$-direction whereas Dirichlet boundary conditions are set in the left and right boundaries. Initial conditions are of the form,

\begin{equation}
	\rho\left(x,y,0\right) = \left\lbrace \begin{array}{lc}
		\rho_{L} & \mathrm{ if } \; x \le 0,\\
		\rho_{R} & \mathrm{ if } \; x> 0;
	\end{array}\right.\qquad
	\press \left(x,y,0\right) = \left\lbrace \begin{array}{lc}
		\press_{L} & \mathrm{ if } \; x \le 0,\\
		\press_{R} & \mathrm{ if } \; x > 0;
	\end{array}\right. \qquad
	{u}_{1} \left(x,y,0\right) =\left\lbrace \begin{array}{lc}
		u_{L}  & \mathrm{ if } \; x \le 0,\\
		u_{R} & \mathrm{ if } \; x > 0;
	\end{array}\right. \qquad
	{u}_{2} \left(x,y,0\right) =0;
\end{equation}
where $\rho_{L}$, $\rho_{R}$, $\press_{L}$, $\press_{R}$, $u_{L}$, $u_{R}$ are defined in Table \ref{tab:RP_IC}.

\begin{table}[H]
	\renewcommand{\arraystretch}{1.2}
	\begin{center}
		\begin{tabular}{|c||c|c|c|c|c|c||c|}
			\hline 
			Test &  $\rho_{L}$ &  $\rho_{R}$ &  $\press_{L}$ &  $\press_{R}$ &  $u_{L}$ &  $u_{R}$ & $t_{\mathrm{end}}$ \\ \hhline{========}
			RP1 & $1 $ & $0.125 $ & $1 $ & $0.1 $ & $0 $ & $0 $ & $0.25$ \\ 	\hline 
			RP2 & $1 $ & $1 $ & $0.4 $ & $0.4 $ & $-1 $ & $1 $ & $0.15$ \\\hline 
			RP3 & $1 $ & $0.125 $ & $1 $ & $1 $ & $0.5 $ & $0 $ & $0.1$ \\ 	
			\hline 
		\end{tabular} 
	\end{center}
	\caption{Riemann problems. Initial condition and final time for each test case.}
	\label{tab:RP_IC}
\end{table}

\begin{figure}[h]
	\centering
	\includegraphics[width=0.325\linewidth]{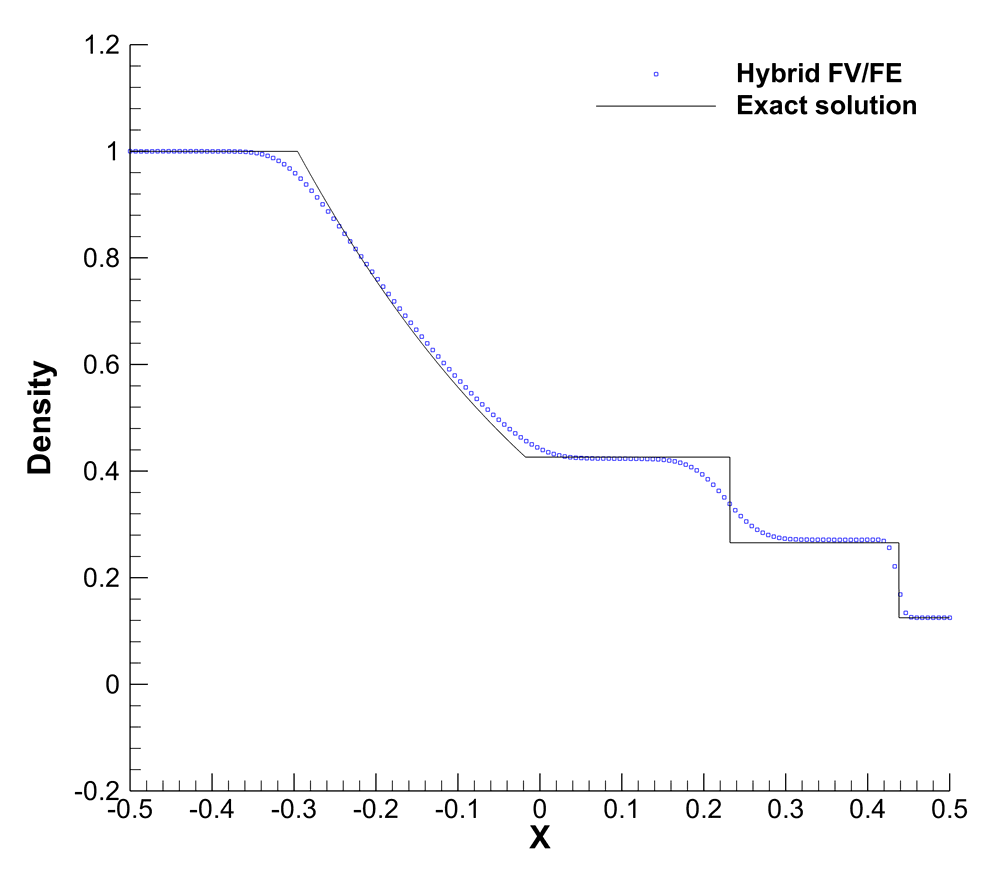}
	\includegraphics[width=0.325\linewidth]{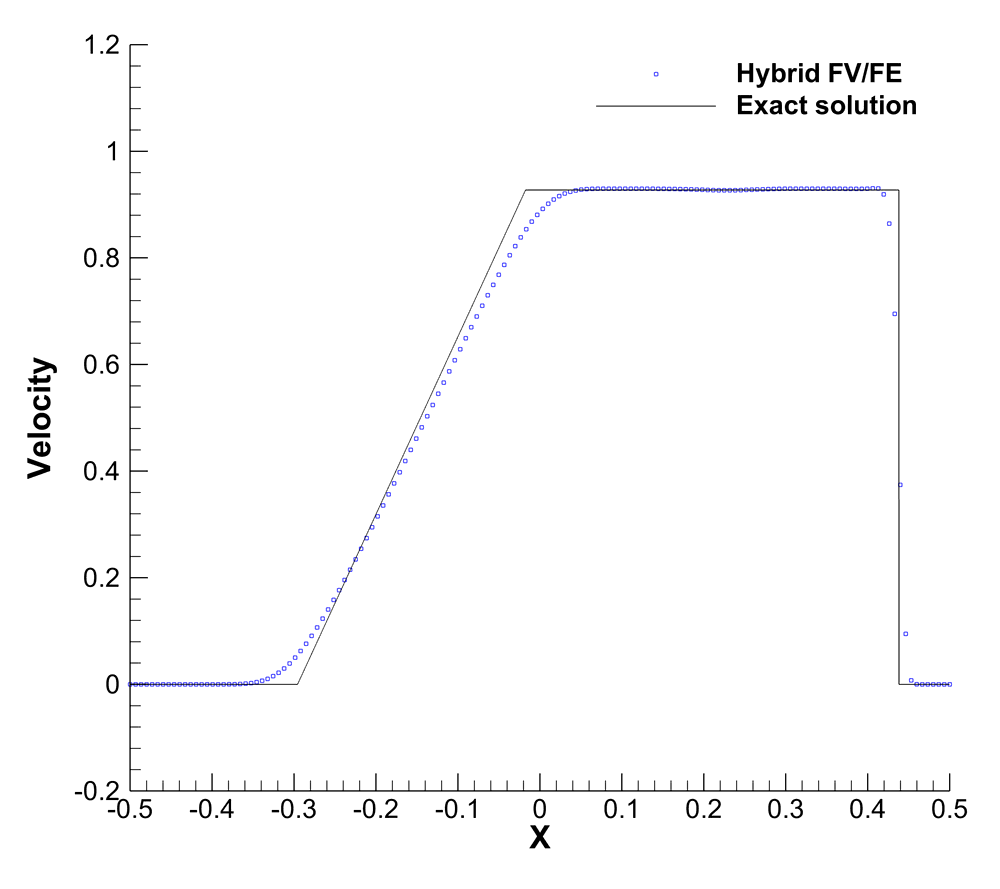}
	\includegraphics[width=0.325\linewidth]{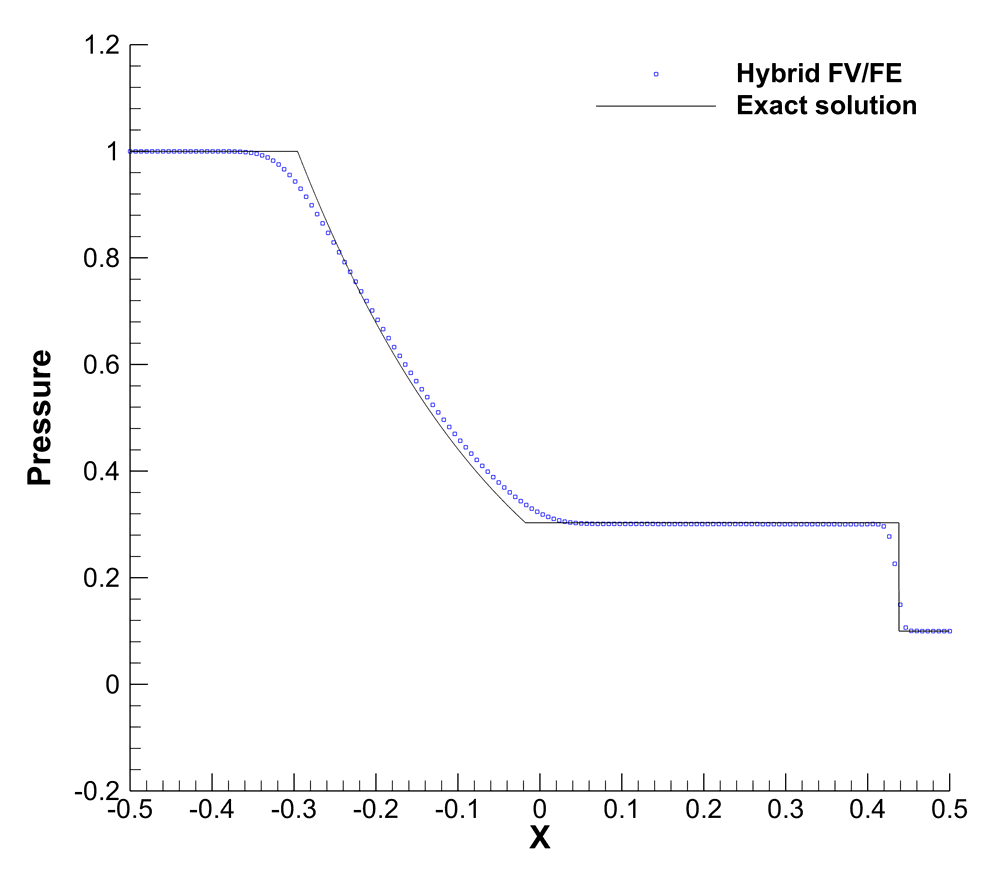}
	
	\includegraphics[width=0.325\linewidth]{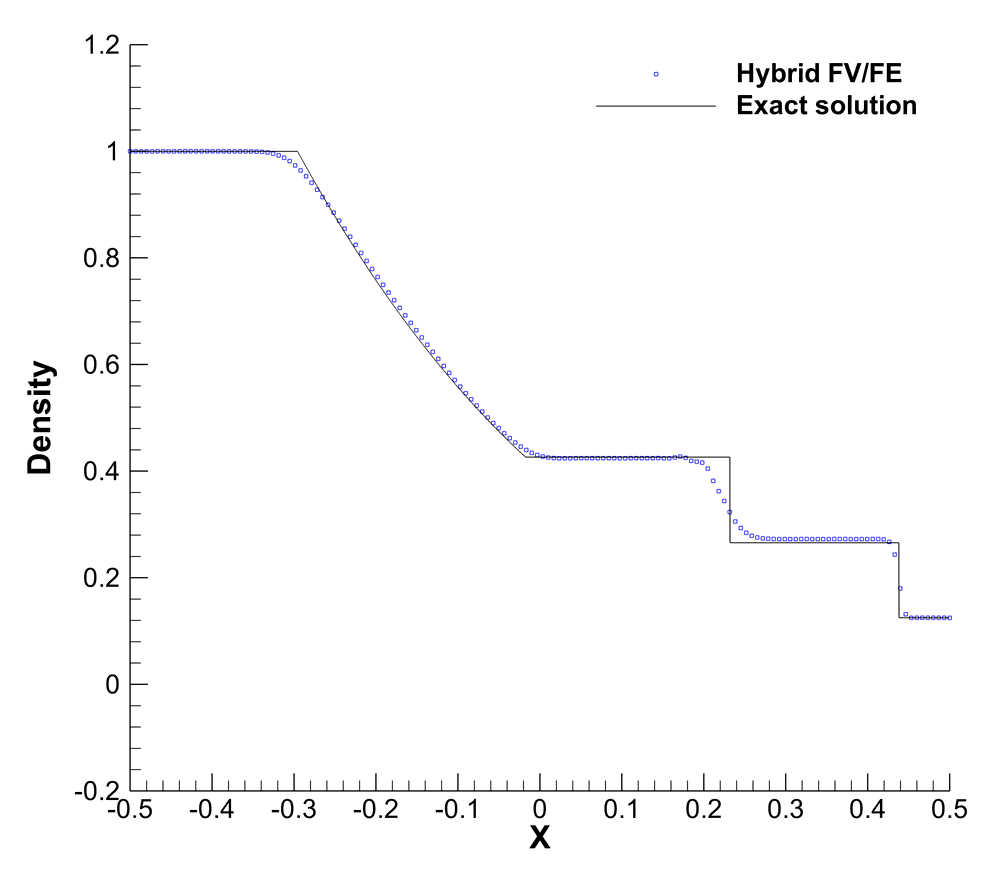}
	\includegraphics[width=0.325\linewidth]{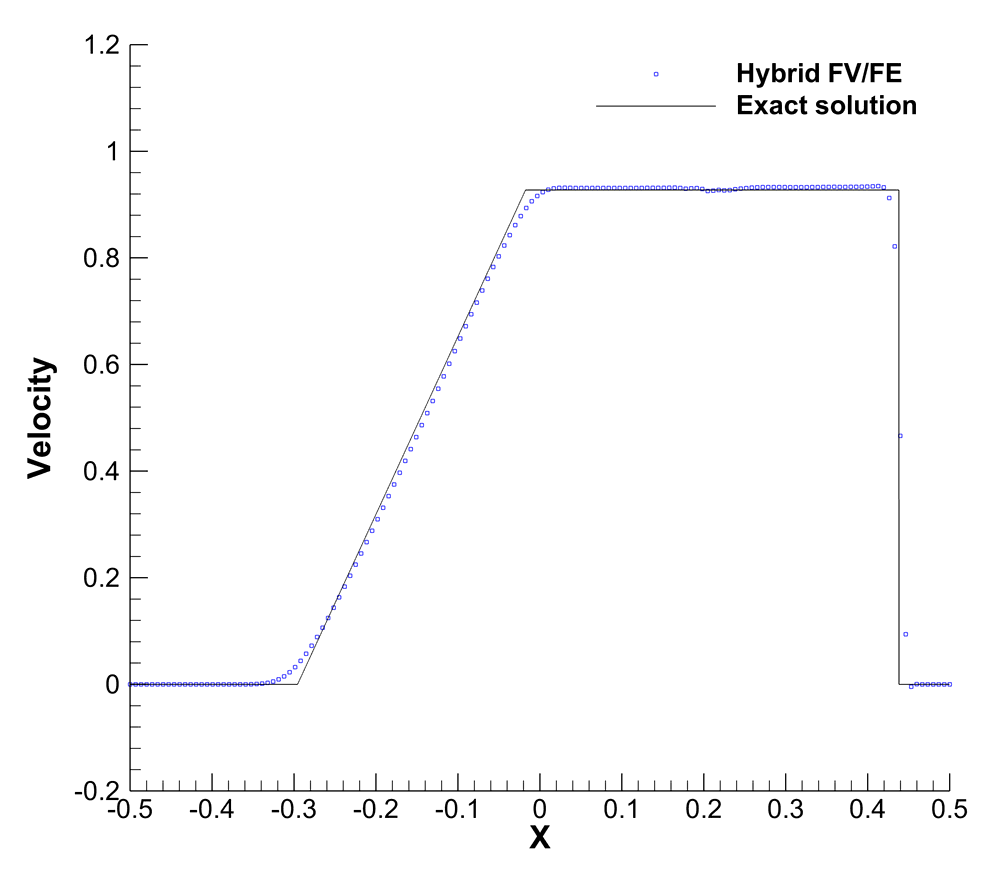}
	\includegraphics[width=0.325\linewidth]{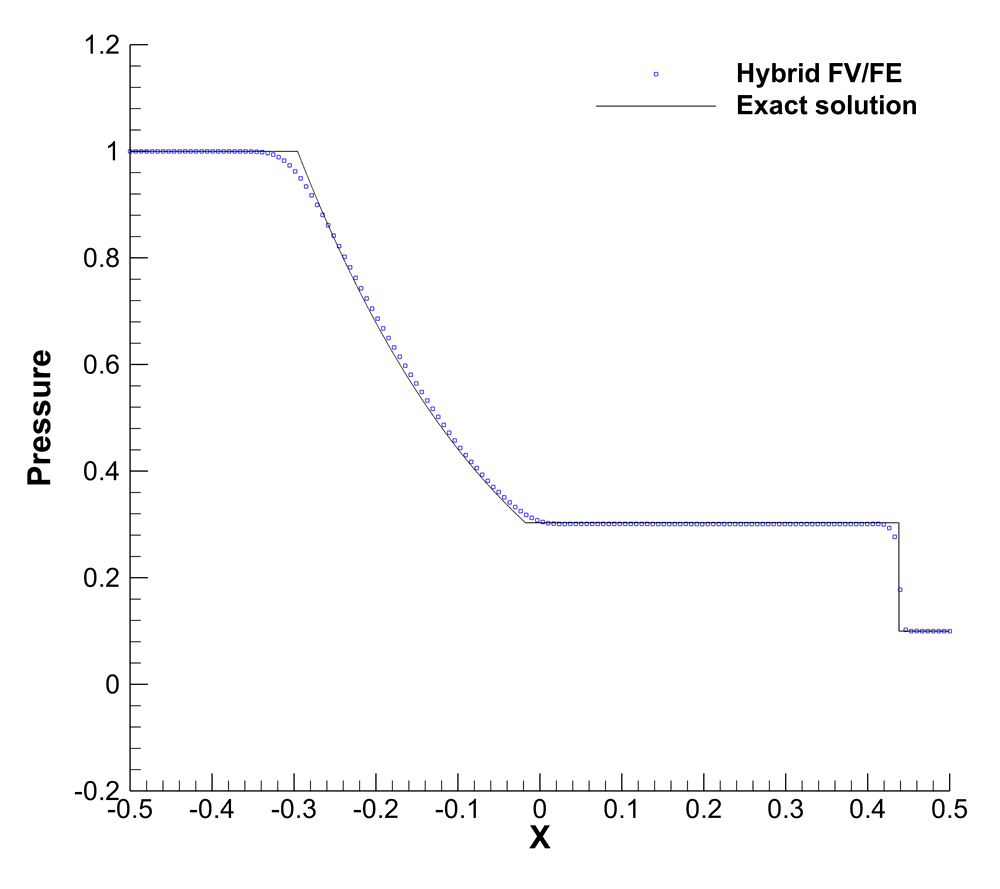}
	\caption{Riemann problem 1. Plot over line $y=0$ of the exact and numerical solutions of density, velocity and pressure at $t_{\mathrm{end}}=0.25$ using LADER-BJ method (bottom) on mesh M1, {\color{black} ($\mathrm{CFL}_{c}=0.56$, $c_{\alpha}=0.8$, $M\approx0.9$)}.}
	\label{fig:RP1_t025_c08}
\end{figure}

The first Riemann problem, RP1, corresponds to the so called Sod problem first put forward in \cite{Sod78}. The results obtained using the first order scheme and the second order LADER method with a Barth and Jespersen limiter (LADER-BJ) show a good agreement with the exact solution even for the coarse grid M1; see Figure \ref{fig:RP1_t025_c08}.

\begin{figure}[h]
	\centering
	\includegraphics[trim= 5 5 5 5,clip,width=0.325\linewidth]{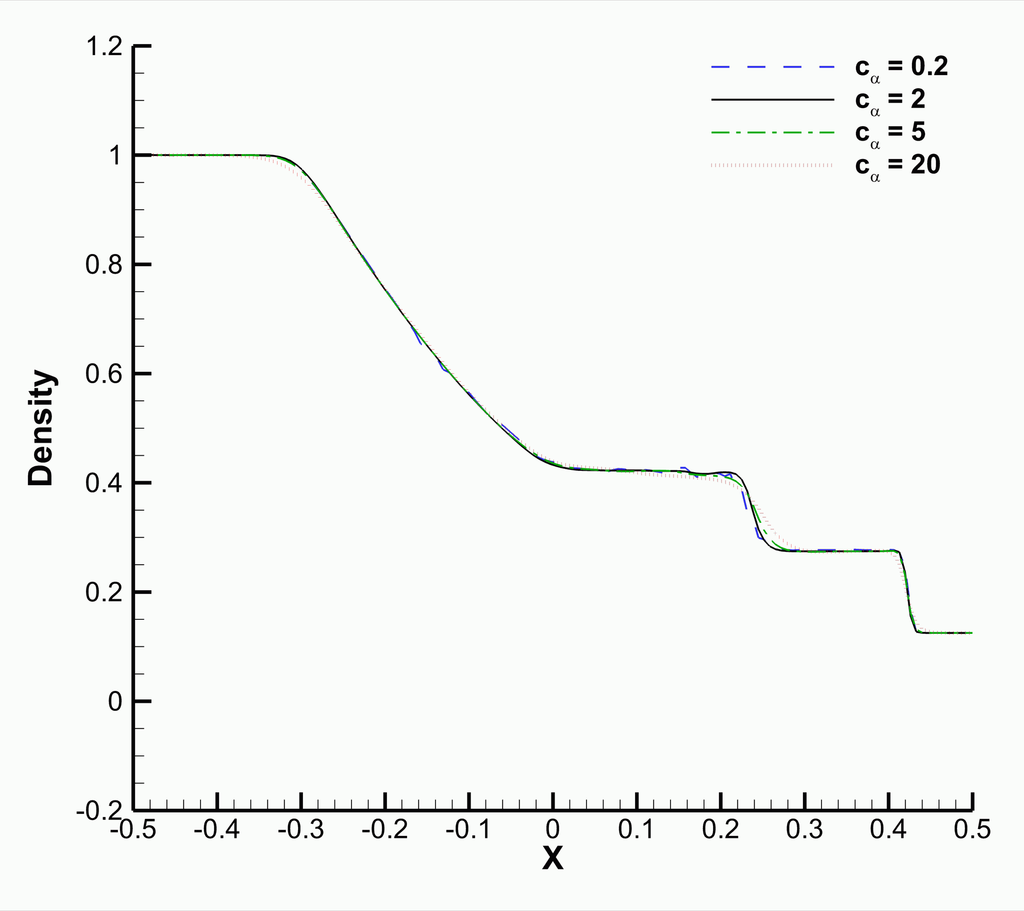}
	\includegraphics[trim= 5 5 5 5,clip,width=0.325\linewidth]{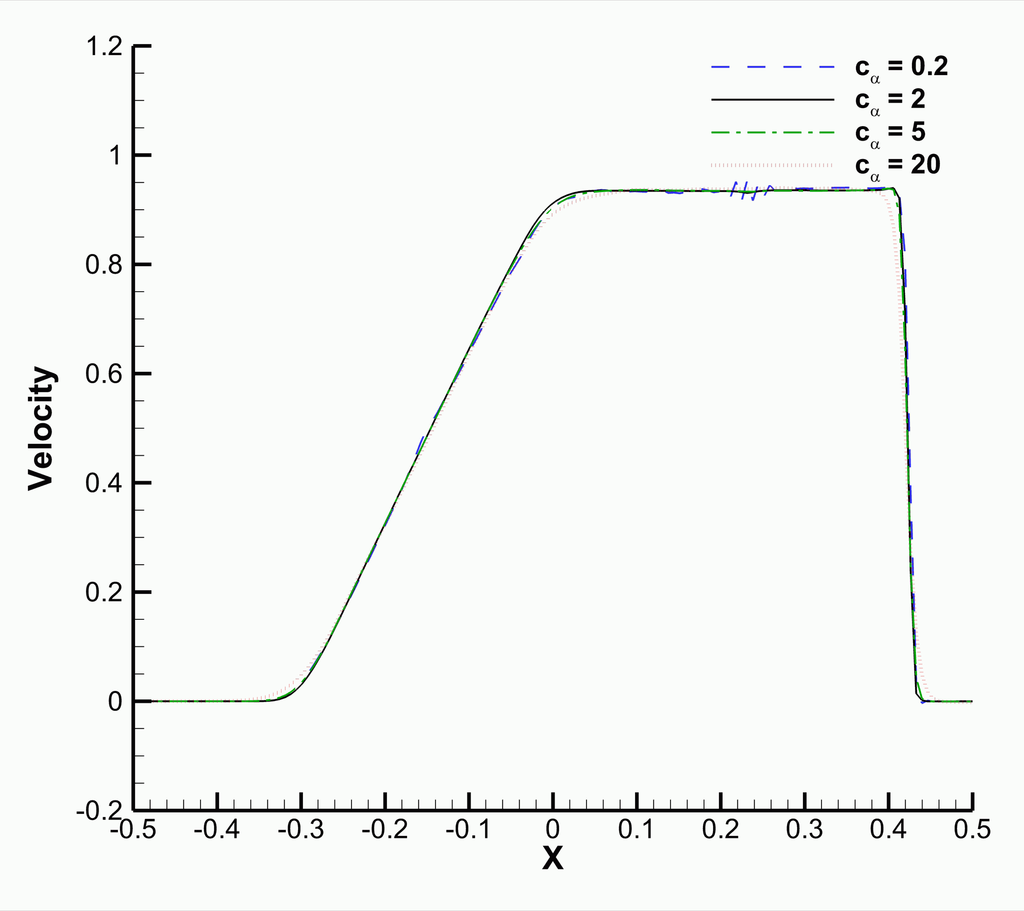}
	\includegraphics[trim= 5 5 5 5,clip,width=0.325\linewidth]{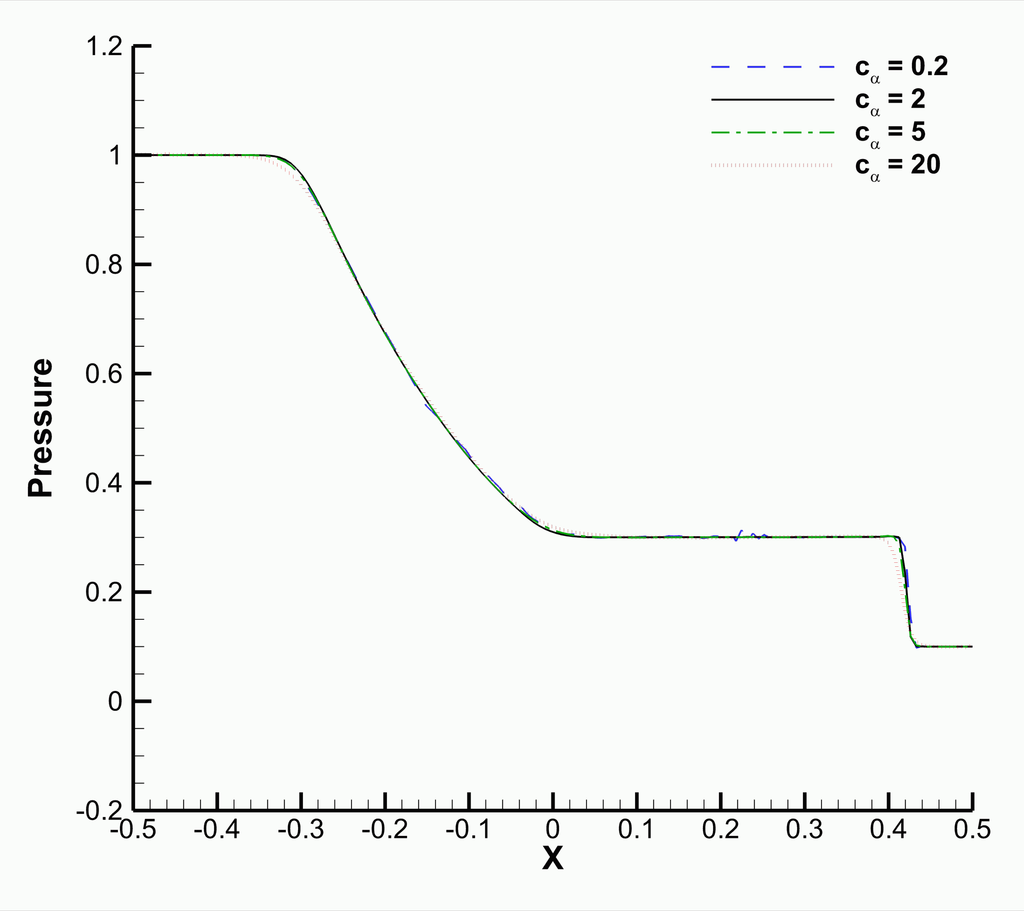}
	
	\caption{{\color{black} Riemann problem 1. Plot over line $y=0$ of the 
			numerical solutions of density, velocity and pressure at $t_{\mathrm{end}}=0.25$ for $c_{\alpha}\in\left\lbrace0.2,2,5,20\right\rbrace$ using the LADER-ENO method on mesh M1,  ($\mathrm{CFL}_{c}=0.56$, $M\approx0.9$)}.}
	\label{fig:RP1_t025_cvar}
\end{figure}

{\color{black} In \eqref{eq:artvisccoef} we have introduced an artificial viscosity coefficient, $c_{\alpha}$, to guarantee the stability of the scheme in the presence of large density variations with respect to the velocity magnitude. 
	Figure \ref{fig:RP1_t025_cvar} depicts the results obtained for RP1 using LADER-ENO for $c_{\alpha}\in\left\lbrace0.2,2,5,20\right\rbrace$. We observe that the artificial viscosity term has a diffusive effect on the solution and stabilizes the density field. Currently, the value of $c_{\alpha}$ is set manually for each test case taking into account the density variations. Nevertheless, we observe that the numerical results are rather insensitive to variations of this parameter within one order of magnitude.
}

\begin{figure}[h]
	\centering
	\includegraphics[width=0.325\linewidth]{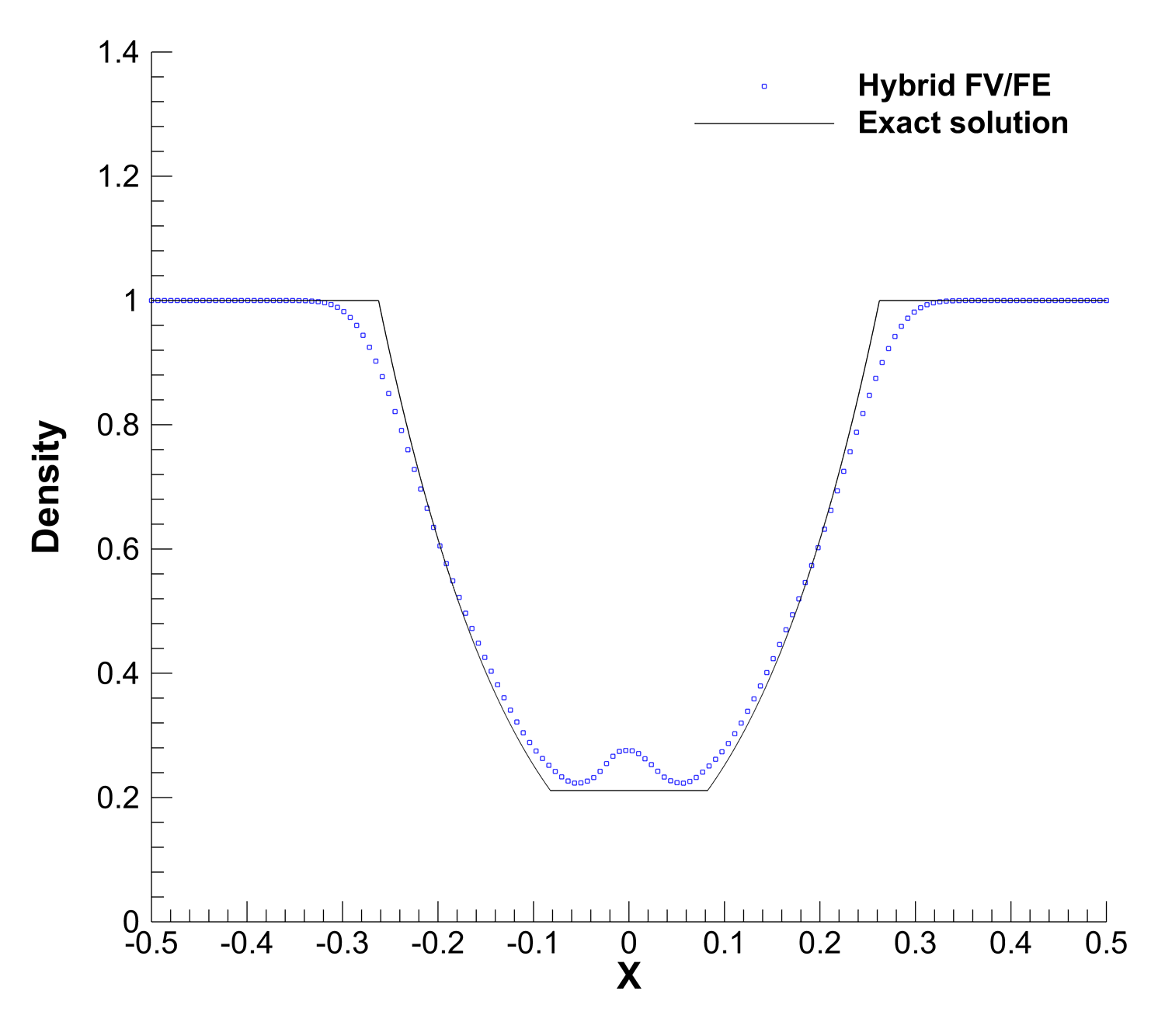}
	\includegraphics[width=0.325\linewidth]{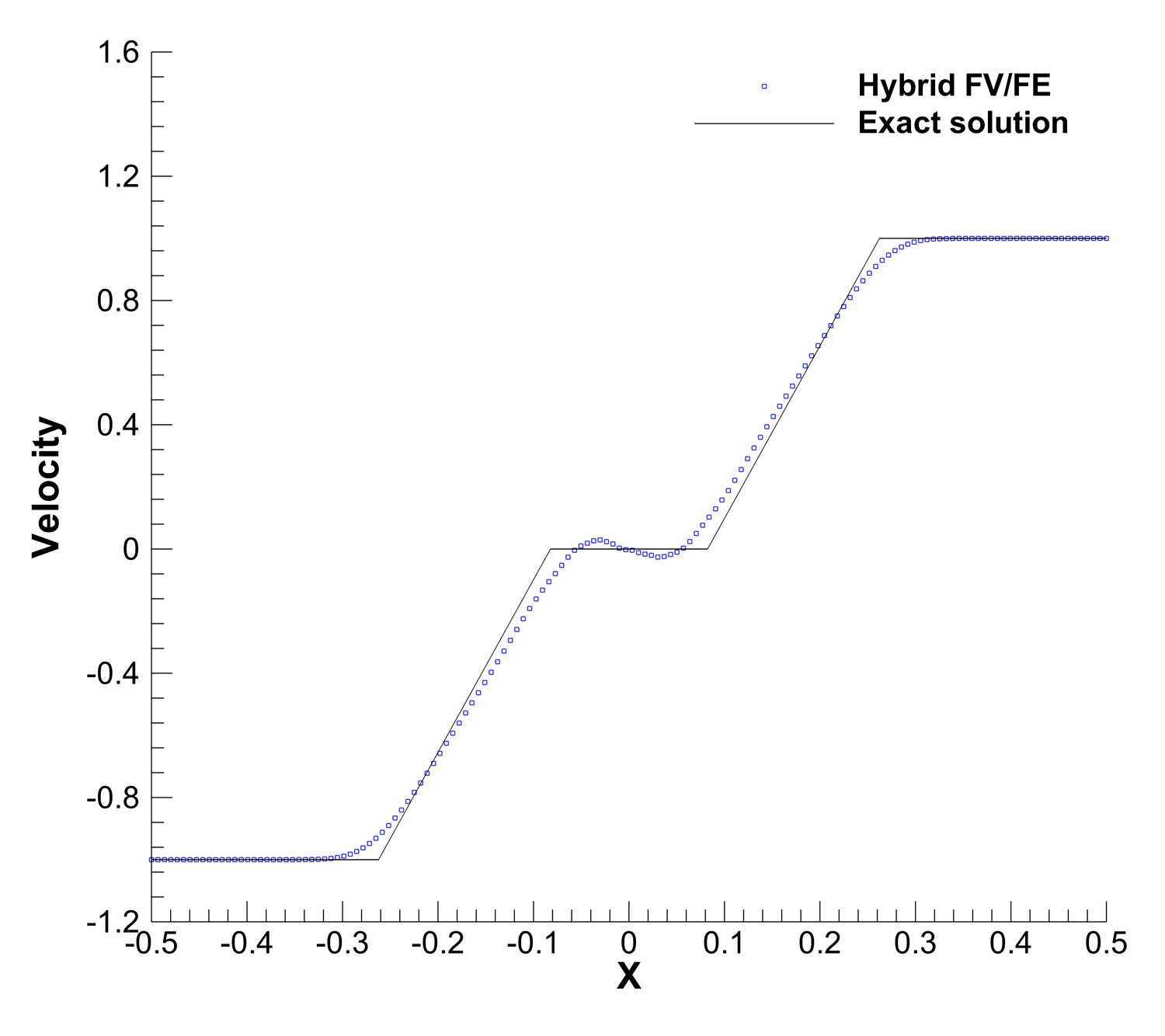}
	\includegraphics[width=0.325\linewidth]{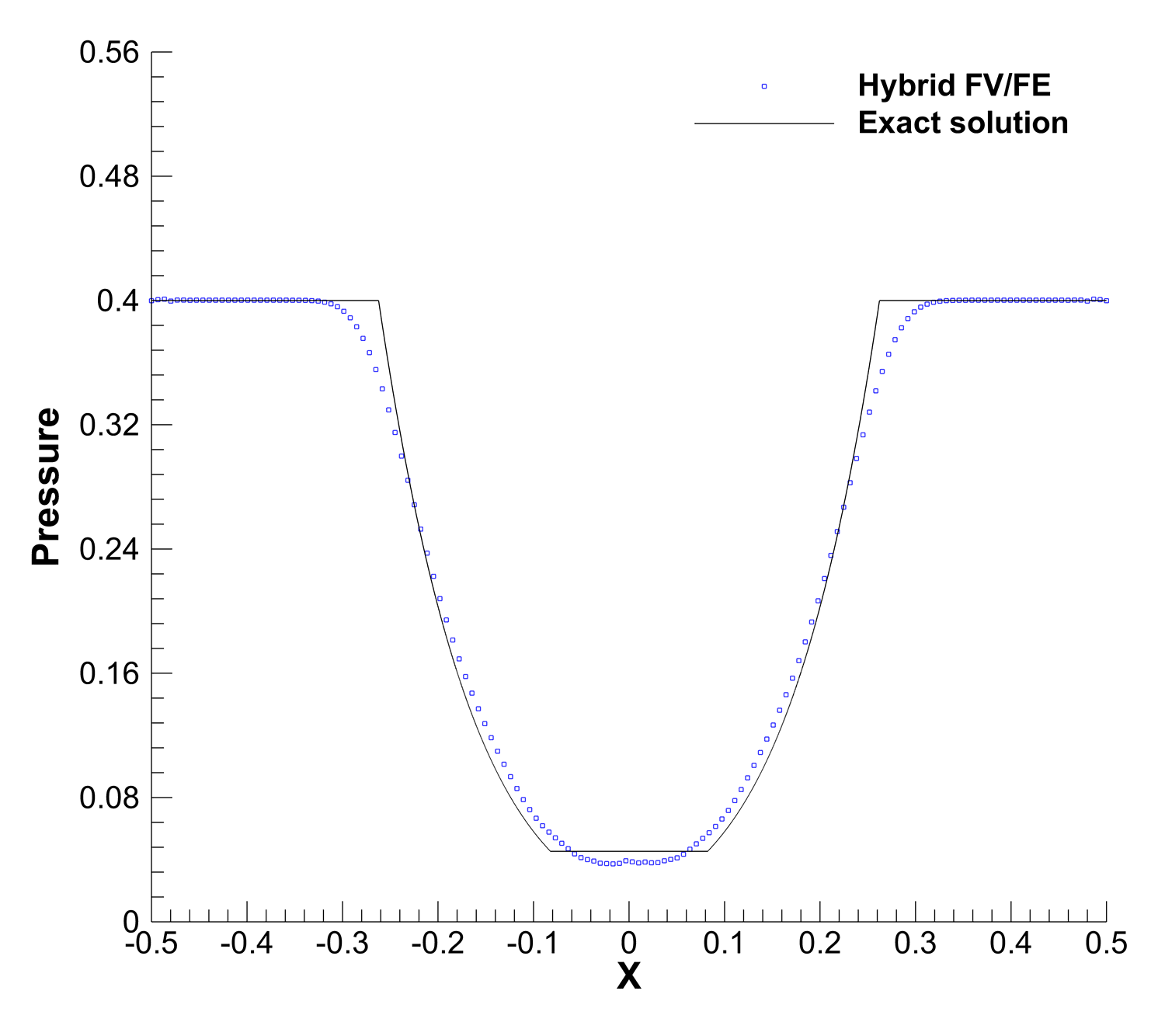}
	
	\includegraphics[width=0.325\linewidth]{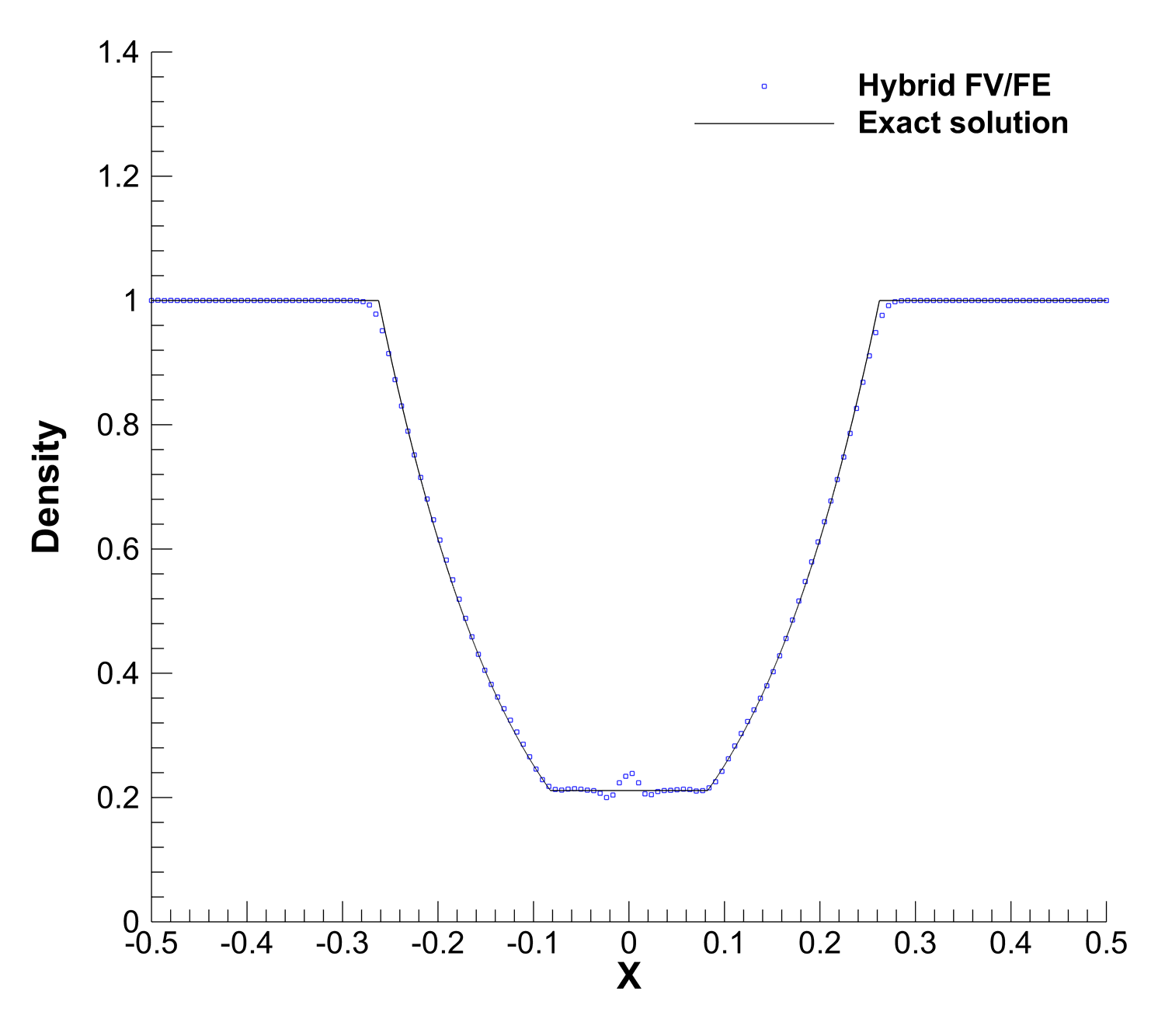}
	\includegraphics[width=0.325\linewidth]{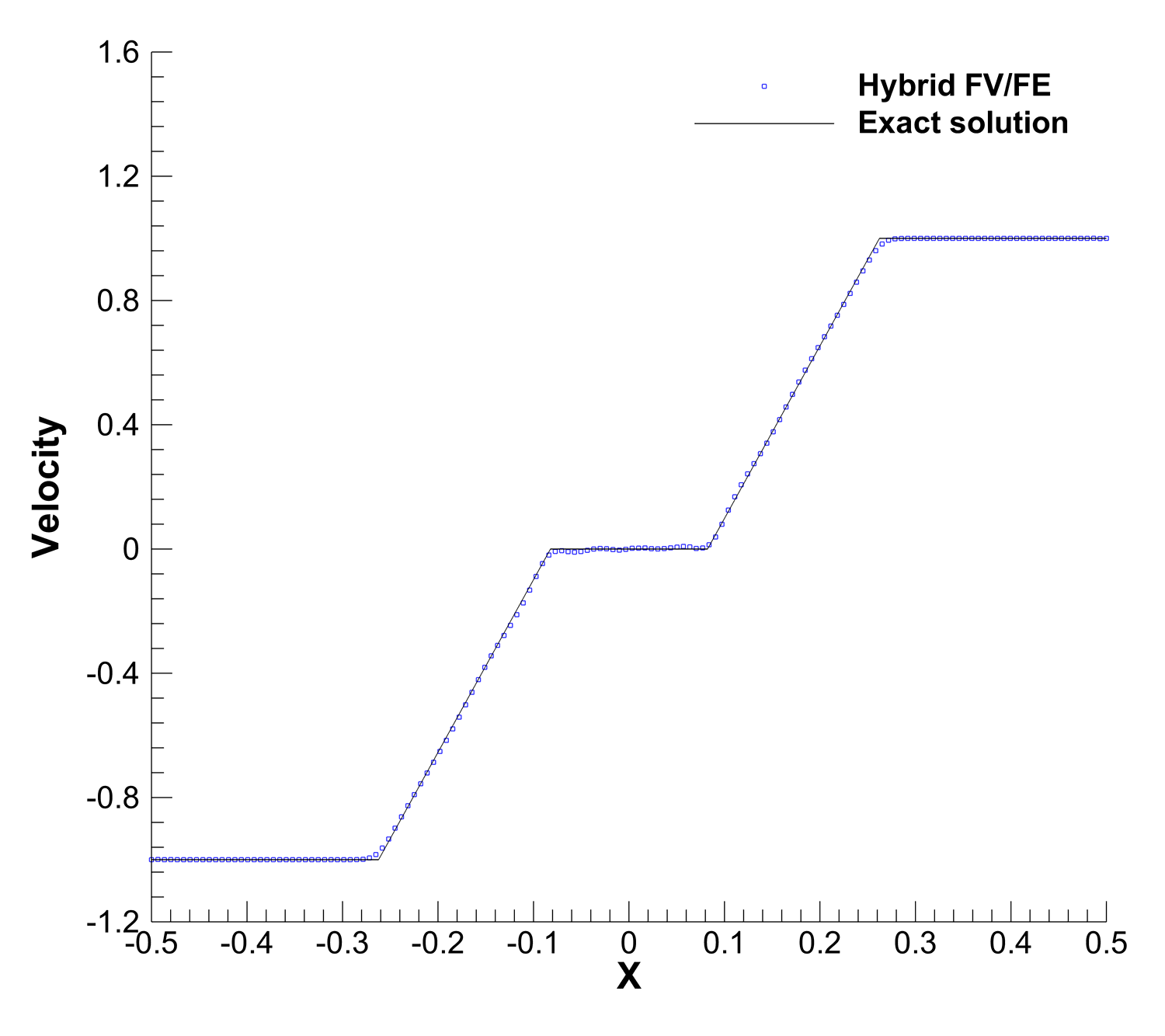}
	\includegraphics[width=0.325\linewidth]{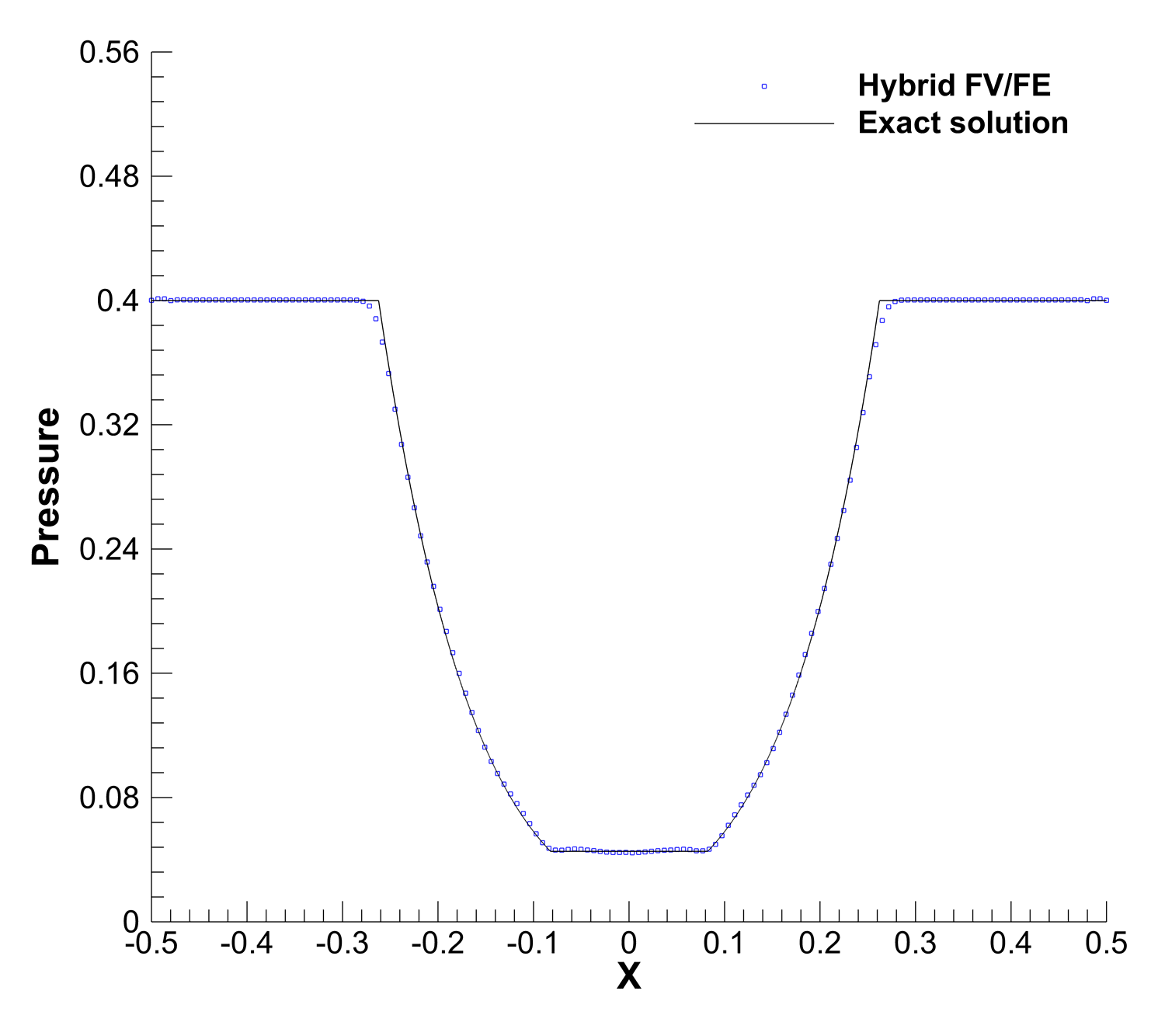}
	\caption{Riemann problem 2. Plot over line $y=0$ of the exact and numerical solutions of density, velocity and pressure at $t_{\mathrm{end}}=0.15$ using the first order scheme (top) and LADER-BJ method (bottom) on mesh M2, {\color{black} ($\mathrm{CFL}_{c}=0.71$, $c_{\alpha}=0.5$, $M\approx 0.07$)}.} 
	\label{fig:RP6_t015_c05} 
\end{figure}

The double rarefaction problem, RP2, in Figure \ref{fig:RP6_t015_c05}, presents a small spurious oscillation of the density at the origin arising from the jump in the initial velocity field. Nevertheless, the magnitude of this unphysical behaviour decreases when refining the mesh and the pressure and velocity field are non affected by it.

\begin{figure}[h]
	\centering
	\includegraphics[width=0.325\linewidth]{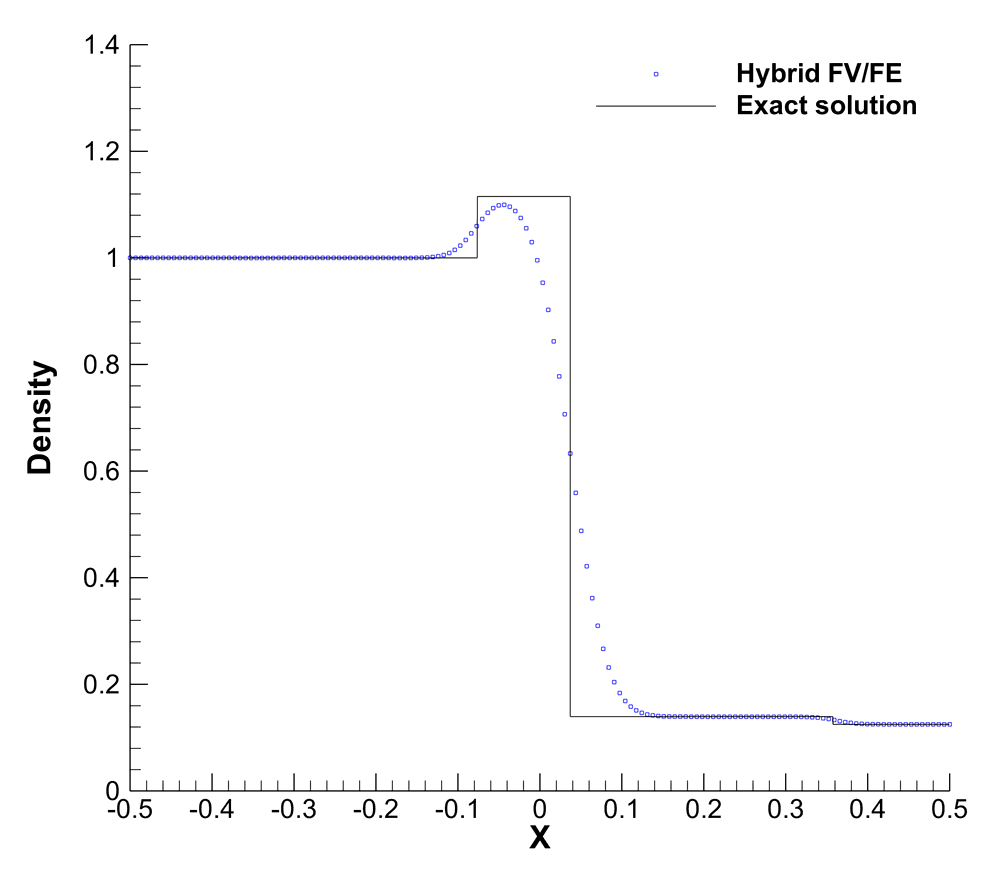}
	\includegraphics[width=0.325\linewidth]{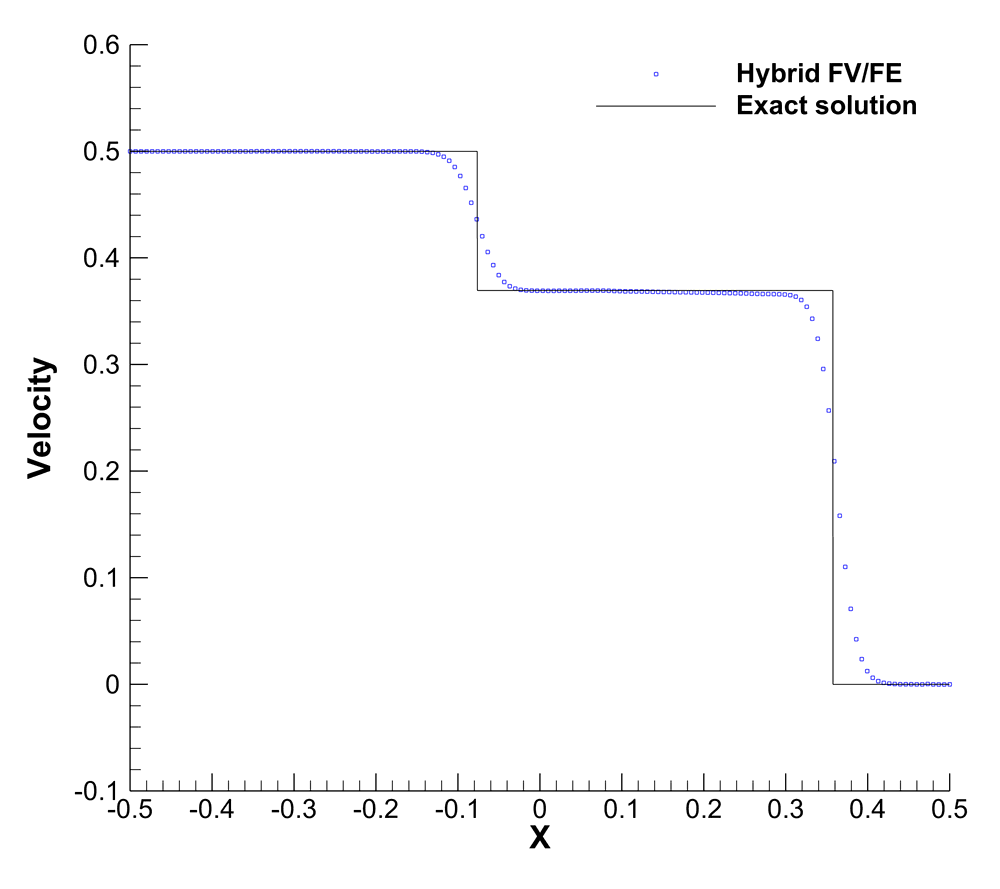}
	\includegraphics[width=0.325\linewidth]{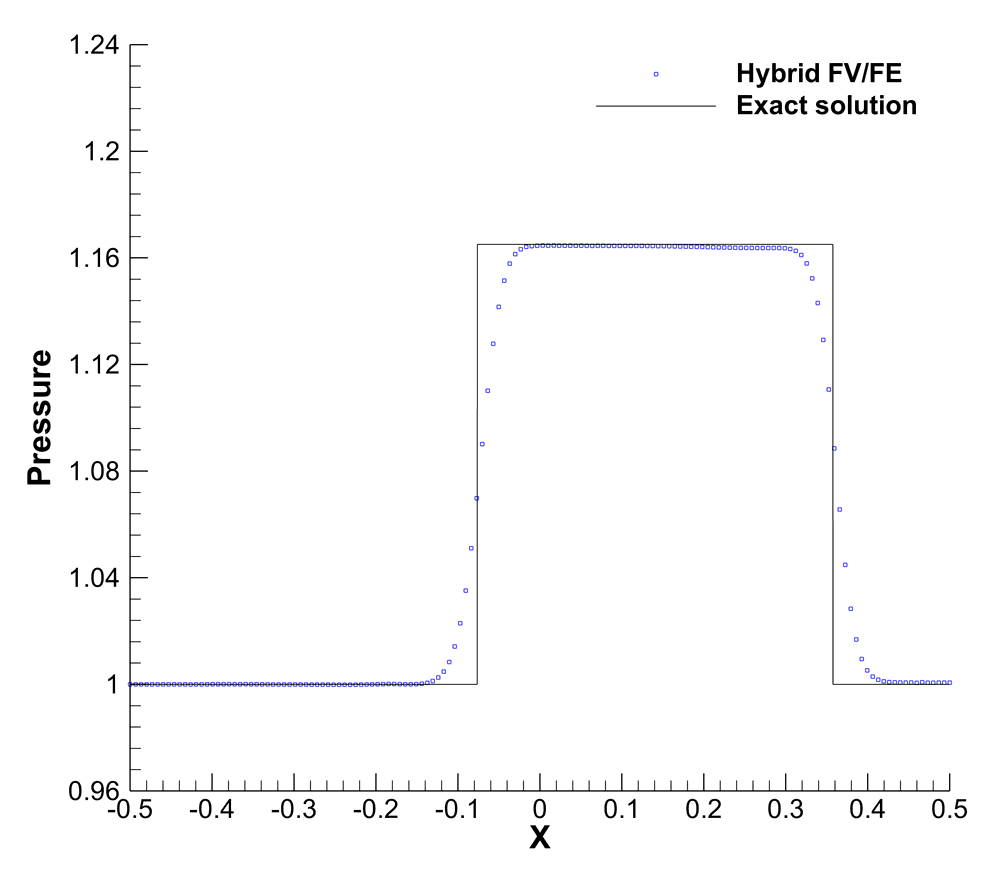}
	
	\includegraphics[width=0.325\linewidth]{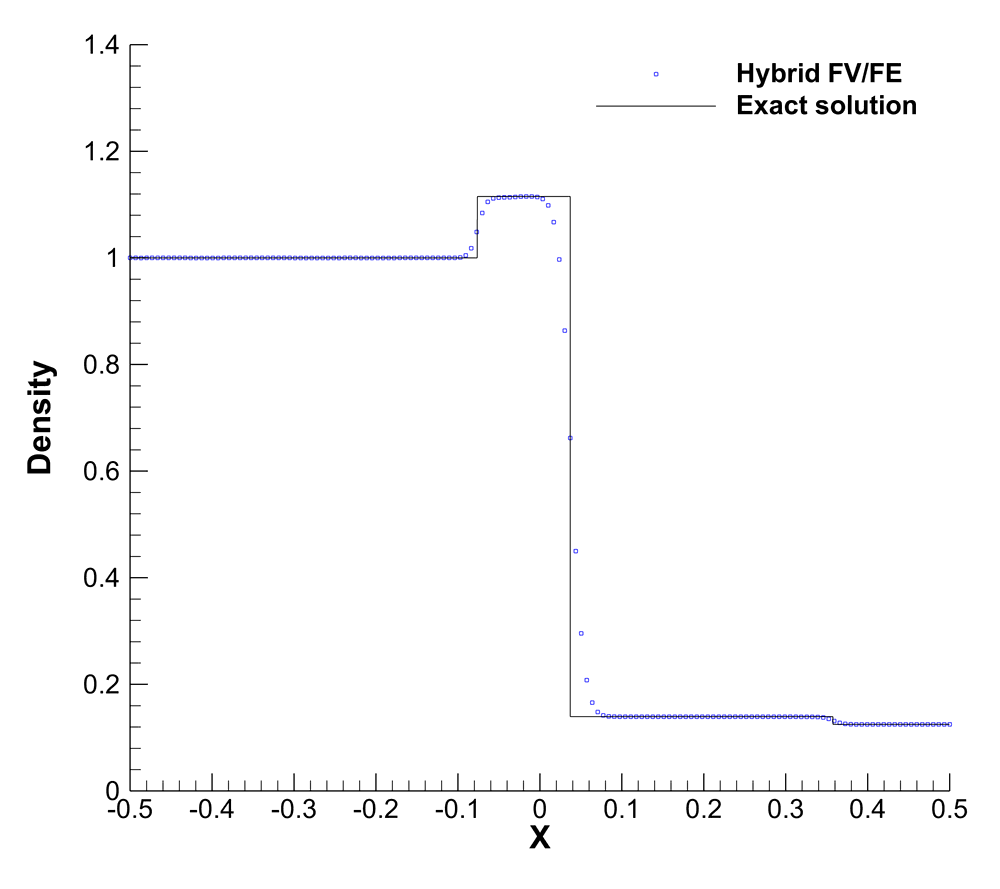}
	\includegraphics[width=0.325\linewidth]{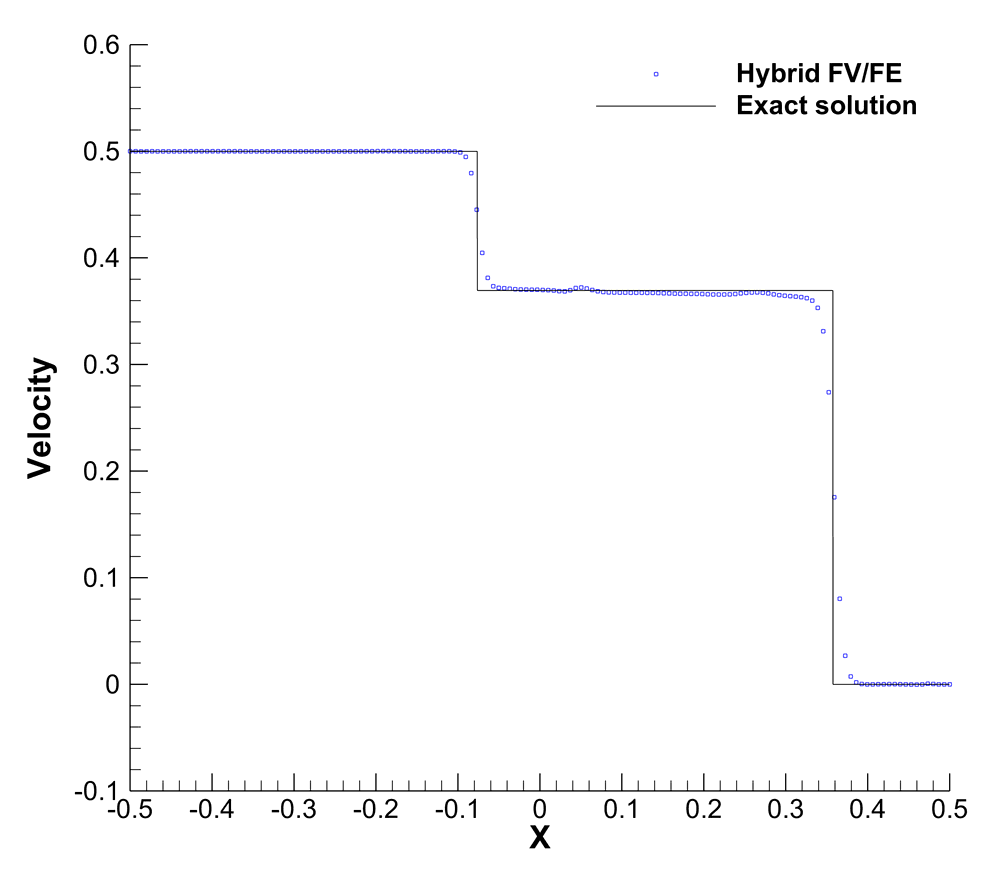}
	\includegraphics[width=0.325\linewidth]{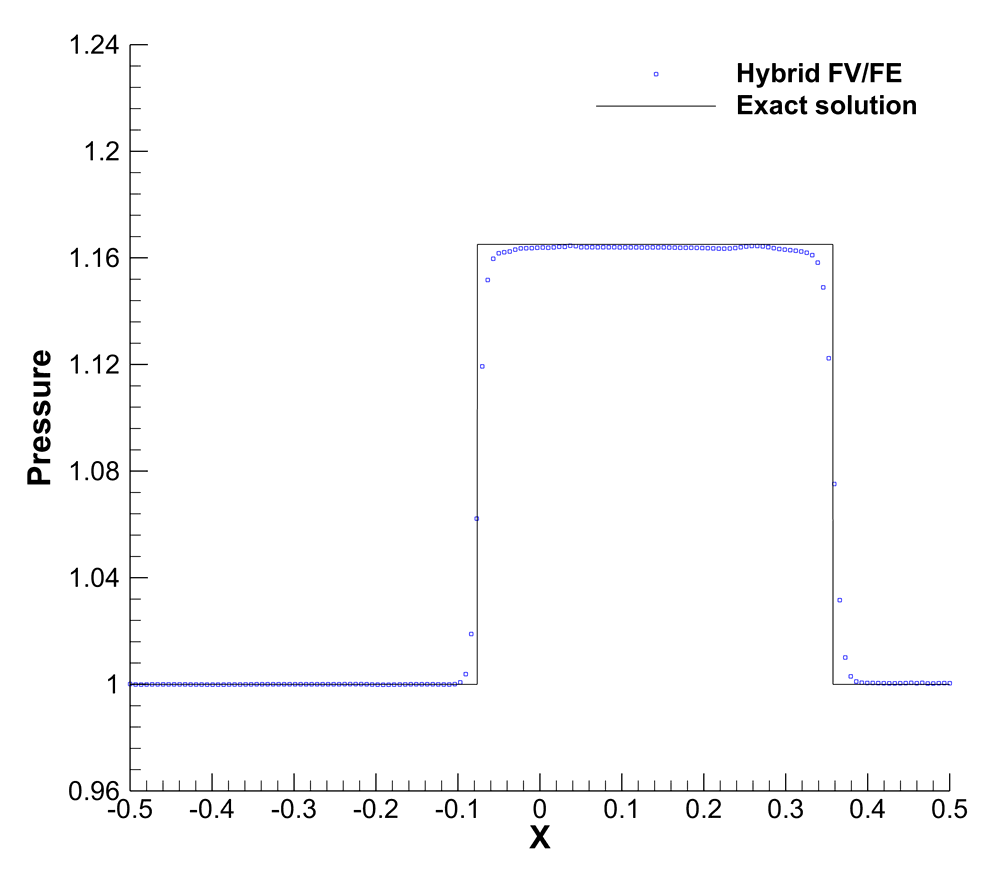}
	\caption{Riemann problem 3. Plot over line $y=0$ of the exact and numerical solutions of density, velocity and pressure at $t_{\mathrm{end}}=0.1$ using the first order scheme (top) and LADER-BJ method (bottom) on mesh M1, {\color{black} ($\mathrm{CFL}_{c}=11.9$, $c_{\alpha}=3$, $M\approx 0.42$)}.} 
	\label{fig:RP3_t01_c3}
\end{figure}

Finally, Figure \ref{fig:RP3_t01_c3} shows the results obtained for the third Riemann problem, RP3, obtained using the coarse mesh M1. Also in this test case the shape of the solution obtained using LADER-BJ agrees pretty well with the exact solution being able to capture also the small jump of the density field at $x=0.36, \, t=0.1$.

\subsection{The first problem of Stokes}\label{sec:1SP}
\begin{figure}[h]
	\centering
	\includegraphics[width=0.325\linewidth]{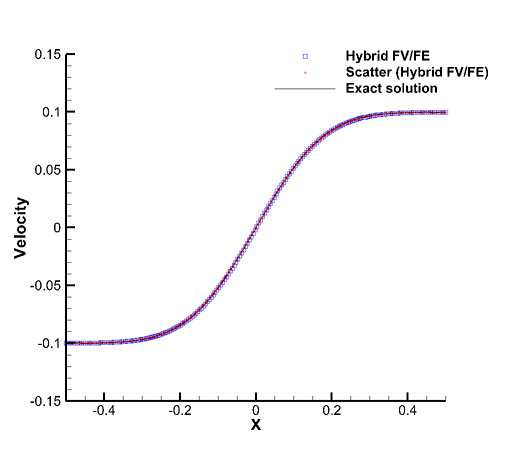}
	\includegraphics[width=0.325\linewidth]{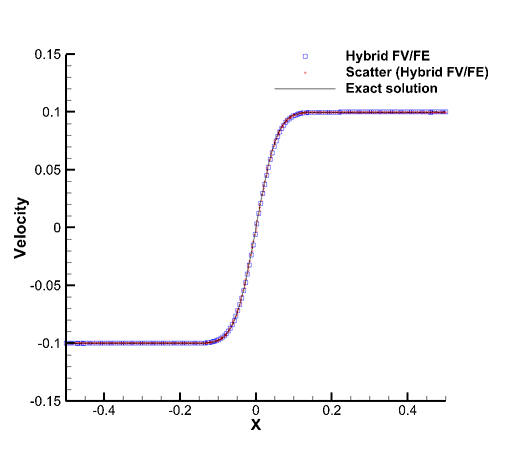}
	\includegraphics[width=0.325\linewidth]{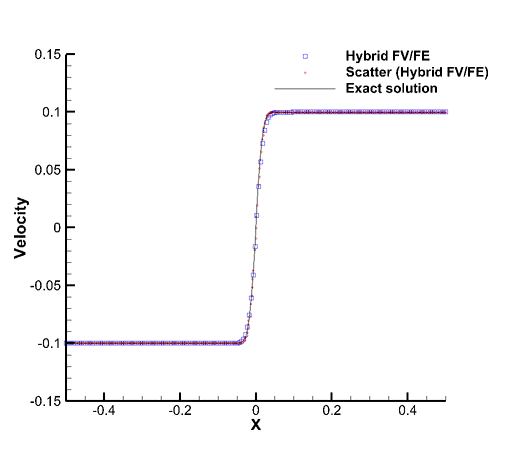}
	\caption{First Stokes problem. Comparison between the exact solution and LADER of the velocity component $u_{2}$ along the cut $y=0$ at $t_{\mathrm{end}}=1$. From left to right: $\mu=10^{-2}$, $\mu=10^{-3}$, $\mu=10^{-4}$.}
	\label{FSP_vvelocity}
\end{figure}

The first problem of Stokes is one of the few  test problems for the unsteady incompressible Navier-Stokes equations with exact analytical solution, \cite{SG16}. We consider the computational domain $\Omega=[-0.5,0.5]\times[-0.5,0.5]$ and the initial condition given by

\begin{equation}
	\rho\left(x,y,0\right) = 1,\qquad
	\press \left(x,y,0\right) = \frac{1}{\gamma}, \qquad
	{u}_{1} \left(x,y,0\right) = 0, \qquad
	{u}_{2} \left(x,y,0\right) = \left\lbrace \begin{array}{lc}
		-0.1 & \mathrm{ if } \; y \le 0,\\
		0.1 & \mathrm{ if } \; y > 0.
	\end{array}\right.
\end{equation}
with exact solution for the vertical velocity:

\begin{equation}
	{u}_{2} \left(x,y,t\right) = \frac{1}{10} \mathrm{erf}\left( \frac{x}{2\sqrt{\mu t}}\right).
\end{equation}
The fluid parameters, $\gamma= c_{\press} = 1.4$, $\lambda=0$, are set to provide a low Mach number, $M=0.1$. Besides, three different viscosities are considered $\mu\in\left\lbrace 10^{-2},10^{-3},10^{-4}\right\rbrace$. Velocity and density are imposed on the left and right boundaries, whereas periodic boundary conditions are set along the $y$-direction. {\color{black} The time step is computed verifying $\mathrm{CFL}_{c}=50$.} The numerical results, obtained at time $t_{\mathrm{end}}=1$ for a 1D cut along $y=0$, have been plotted jointly with the exact solution in \mbox{Figure \ref{FSP_vvelocity}}. Even with a quite coarse mesh of only $1000$ triangular elements, the algorithm provides an accurate solution.

\subsection{Lid driven cavity}\label{sec:LDC}
\begin{figure}[h]
	\centering
	\includegraphics[width=0.4\linewidth]{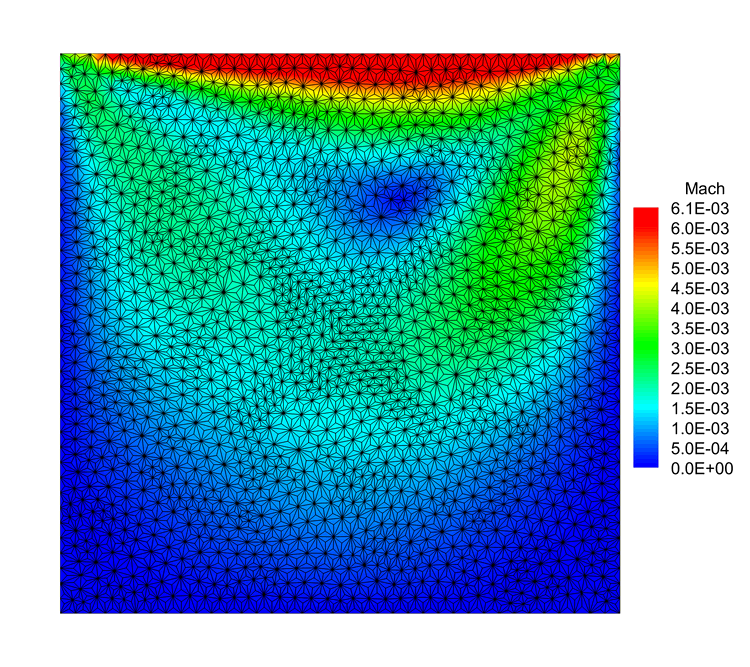}\hspace{0.1\linewidth}
	\includegraphics[width=0.4\linewidth]{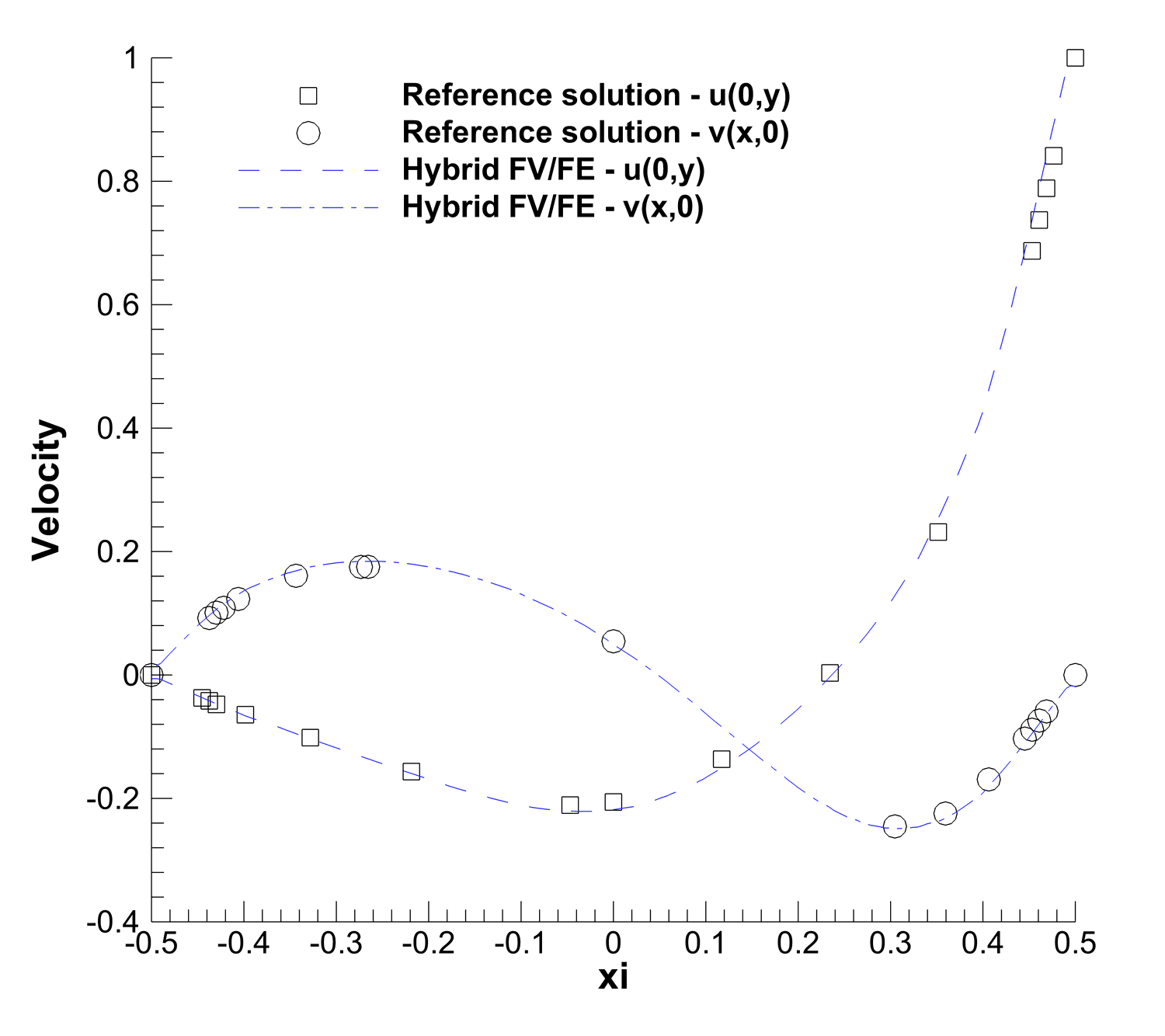}
	\caption{Lid driven cavity test ($Re=100$). Left: Mach number contour plot and dual grid (each triangle corresponds with one half of a dual element). Right: velocity profiles compared with the results given by \cite{GGS82}.}
	\label{LDC_figure}
\end{figure}

The lid driven cavity test is a classical benchmark employed to assess incompressible flow solvers since its introduction  in \cite{GGS82}. Thus, it is a good candidate to test the developed algorithm in the incompressible limit. The computational domain, $\Omega=\left[-0.5,0.5\right]\times\left[-0.5,0.5\right]$, is discretized using an unstructured mesh of $2906$ primal elements {\color{black} and the time step is computed under condition $\mathrm{CFL}_{c}=125$}.
Initially, the flow is considered at rest, the density is set to $\rho = 1$, the pressure is set to $p=10^4$ and we impose homogeneous non slip boundary conditions on the lateral and bottom boundaries and a wall boundary condition with fixed velocity field $\mathbf{u}=\left(1,0\right)^{T}$ at the top boundary. Moreover, we fixed the viscosity $\mu=10^{-2}$ so that $Re=100$. In Figure \ref{LDC_figure}, we compare the results obtained using LADER-ENO against the reference solution in \cite{GGS82}. Also the Mach number contour plot has been included, showing the main features of the flow. 
The characteristic Mach number of this flow based on the lid velocity is $M \approx 8 \cdot 10^{-3}$.

\subsection{Double shear layer}\label{sec:DSL}
\begin{figure}[h]
	\centering
	\includegraphics[width=0.4\linewidth]{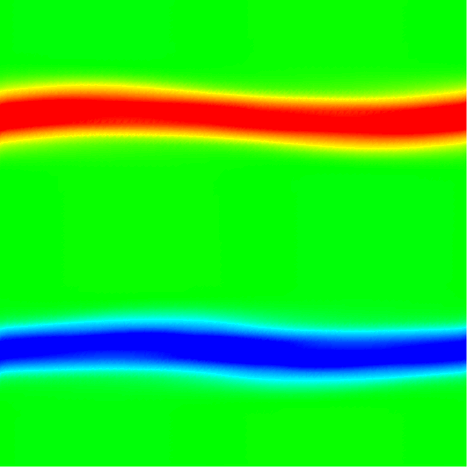}\hspace{0.05\linewidth}
	\includegraphics[width=0.4\linewidth]{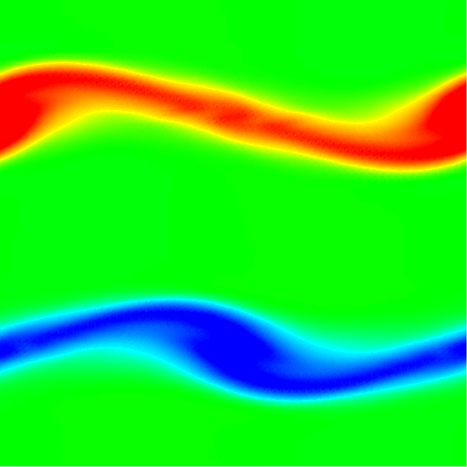}
	
	\vspace{0.05\linewidth}
	\includegraphics[width=0.4\linewidth]{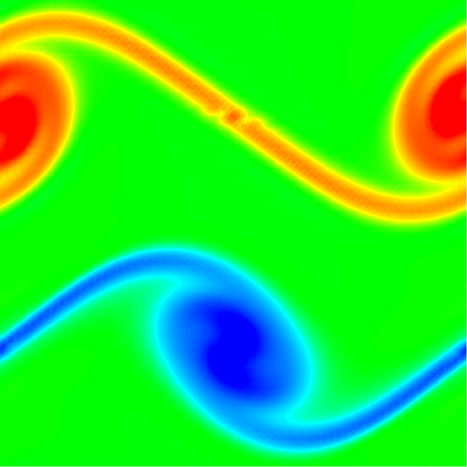}\hspace{0.05\linewidth}
	\includegraphics[width=0.4\linewidth]{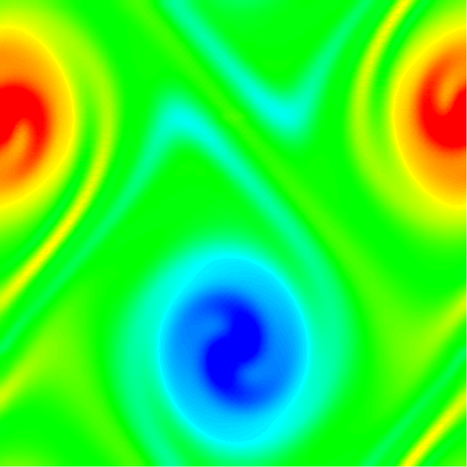}
	\caption{Double shear layer. Vorticity contours obtained using CVC scheme {\color{black} ($\mathrm{CFL}_{c}=250$)}. Time instants from top left to bottom right: $t= 0.8$, $t=1.6$, $t=2.4$, $t=3.6$.}
	\label{DSL64_cvcnl}
\end{figure}

\begin{figure}[h]
	\centering
	\includegraphics[width=0.4\linewidth]{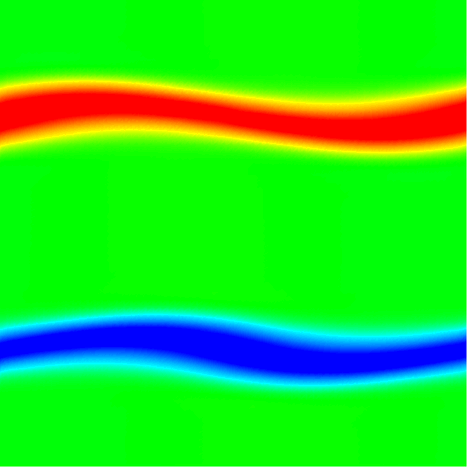}\hspace{0.05\linewidth}
	\includegraphics[width=0.4\linewidth]{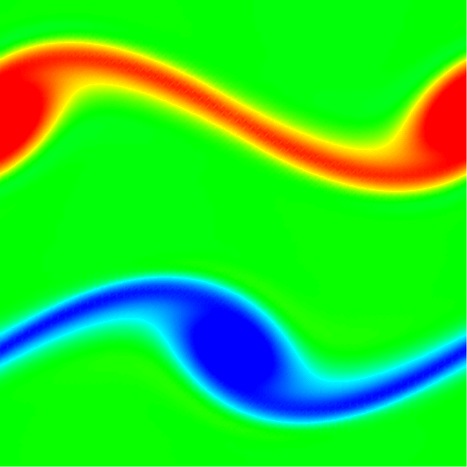}
	
	\vspace{0.05\linewidth}
	\includegraphics[width=0.4\linewidth]{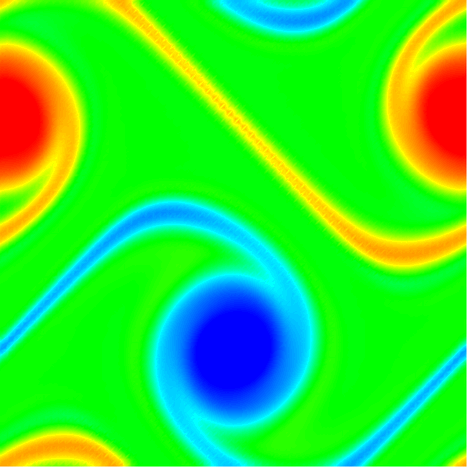}\hspace{0.05\linewidth}
	\includegraphics[width=0.4\linewidth]{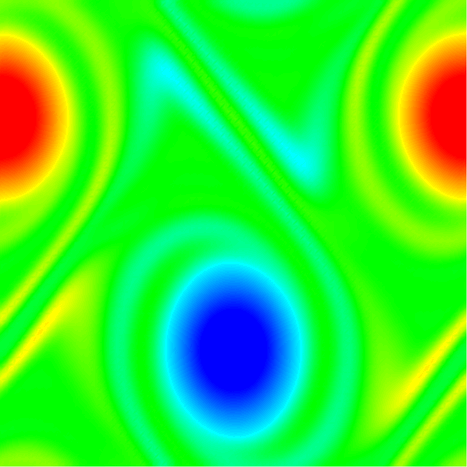}
	\caption{Double shear layer. Vorticity contours obtained using LADER scheme {\color{black} ($\mathrm{CFL}_{c}=250$)}. Time instants from top left to bottom right: $t= 0.8$, $t=1.6$, $t=2.4$, $t=3.6$}
	\label{DSL64_Lnl}
\end{figure}

The double shear layer test problem, \cite{BCG89}, is also used to assess the behaviour of the algorithm in the incompressible limit. The initial flow is characterized by a high velocity gradient which produces complex flow patterns. In particular, we consider the computational domain $\Omega=[-1,1]\times[-1,1]$ and the perturbed double shear layer profile

\begin{equation}
	\rho\left(x,y,0\right) = 1,\qquad
	\press \left(x,y,0\right) = \frac{10^{5}}{\gamma}, \qquad
	{u}_{1} \left(x,y,0\right) = \left\lbrace \begin{array}{lc}
		\tanh \left[\hat{\rho}(\hat{y}-0.25)\right] & \mathrm{ if } \; \hat{y} \le 0.5,\\
		\tanh \left[\hat{\rho}(0.75-\hat{y})\right] & \mathrm{ if } \; \hat{y} > 0.5,
	\end{array}\right.\quad
	{u}_{2} \left(x,y,0\right) = \delta \sin \left(2\pi \hat{x}\right)
\end{equation}
with $\hat{x}= \frac{x+1}{2}$ and $\hat{y}= \frac{y+1}{2}$ normalized vertical and horizontal coordinates, $\hat{\rho} = 30$ the parameter that determines the slope of the shear layer, $\mu= 2\cdot 10^{-4}$ the viscosity and $\delta =0.05$ the amplitude of the initial perturbation. The characteristic Mach number of this flow problem is of the order $M \approx 2 \cdot 10^{-3}$. Periodic boundary conditions are applied everywhere. The used primal mesh has $8192$ triangles. 
The vorticity contours for several time instants, $t\in\left\lbrace 0.8, 1.6, 2.4, 3.6\right\rbrace$, are depicted in \mbox{Figure \ref{DSL64_cvcnl}} for the CVC scheme and in \mbox{Figure \ref{DSL64_Lnl}} for LADER scheme. Overall, the flow structure seems to be well resolved, although only a second order scheme is employed (see, for instance, \cite{DPRZ16,TD15} for comparison with solutions obtained with higher order staggered semi-implicit discontinuous Galerkin schemes). {\color{black} The computational time required for each simulation is given in Table \ref{tab:DSL64_times}.}

\begin{table}[H]
	\renewcommand{\arraystretch}{1.2}
	\begin{center}
		\begin{tabular}{|c|c||c|}
			\hline
			\multirow{3}{*}{CVC} & CPU time (s)&
			$263.03$ 
			\\\hhline{|~|-||-|}
			& $t_{e}$ ($\mu$s) &			
			$7.01$ 
			\\\hhline{|~|-||-|}
			& Time steps &
			$3021$ 
			\\\hline
			\multirow{3}{*}{LADER} & CPU time (s)&
			$291.14$ 
			\\\hhline{|~|-||-|}
			& $t_{e}$ ($\mu$s) &			
			$7.46$ 		
			\\\hhline{|~|-||-|}
			& Time steps &
			$3145$ 
			\\\hline
		\end{tabular}
		{\color{black} \caption{Double shear layer.
				CPU time, CPU time per element and iteration, $t_{e}$, and number of time steps.
			}\label{tab:DSL64_times}}
	\end{center}
\end{table}

\subsection{Smooth acoustic wave}\label{sec:SAW}
\begin{figure}[h]
	\centering
	\includegraphics[width=0.45\linewidth]{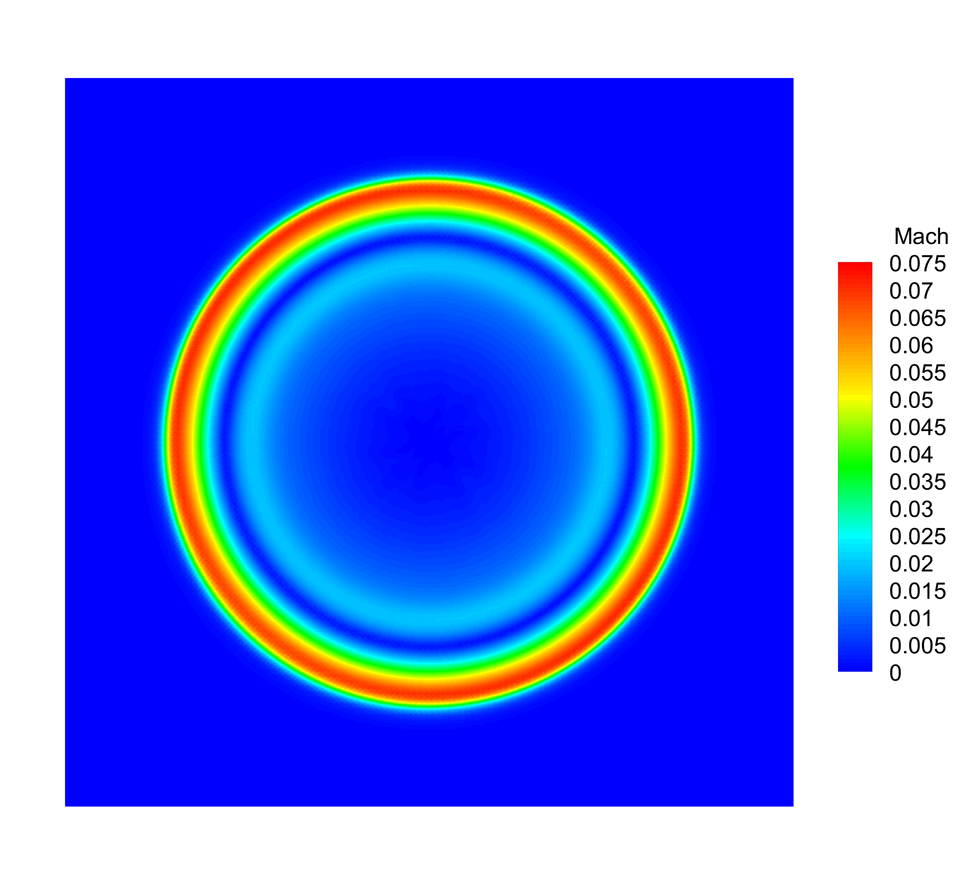}\hspace{0.05\linewidth}
	\includegraphics[width=0.45\linewidth]{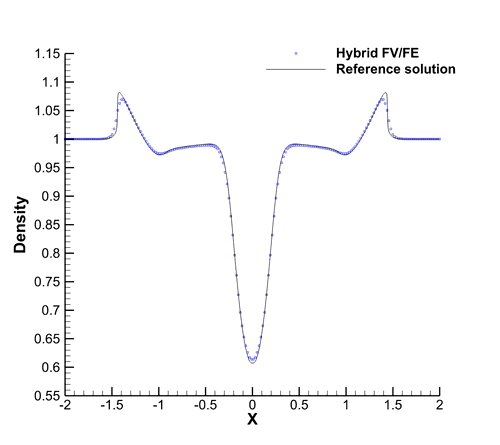}
	
	\vspace{0.05\linewidth}
	\includegraphics[width=0.45\linewidth]{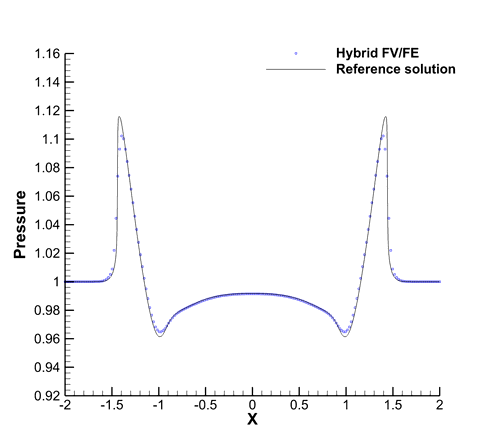}\hspace{0.05\linewidth}
	\includegraphics[width=0.45\linewidth]{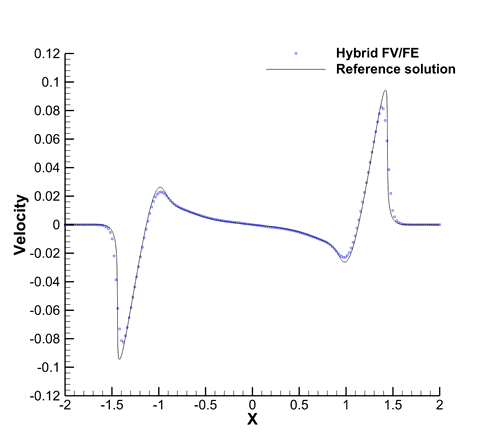}
	\caption{Smooth acoustic wave. Solution obtained at $t_{\mathrm{end}}=1$ {\color{black}($\mathrm{CFL}_{\mathbf{u}}=0.1$, $\mathrm{CFL}_{c}=1.33$)}. From left top to right bottom: Mach number contours, 1D profile (blue squares) and reference solution (black line) of density, pressure and horizontal velocity, ${u}_1$.}
	\label{SAW128}
\end{figure}

To assess the correct propagation of sound waves, which are a characterizing feature that distinguishes weakly compressible flows from incompressible ones, we consider a smooth acoustic wave problem (see \cite{TD17}). 
The two-dimensional computational domain is $\Omega=[-2,2]\times[-2,2]$  and the initial conditions are given by

\begin{gather}
	\rho\left(x,y,0\right) = 1,\qquad
	\press \left(x,y,0\right) = 1+\exp \left(-\alpha r^2\right), \qquad
	\mathbf{u} \left(x,y,0\right) = 0,
\end{gather}
where $r^2=x^2+y^2$ is the distance to the origin and we set $\alpha=40$ and $\mu=\lambda=0$. Moreover, we consider periodic boundary conditions everywhere. 
The angular symmetry of the problem allows us the computation of a reference solution by simply solving an equivalent one-dimensional PDE in the radial direction with a geometrical source term (see \cite{Toro}). More precisely, we have employed a second order TVD scheme in order to compute the reference solution on a 1D grid of $10^{4}$ elements. For the 2D case we consider a primal triangular mesh made of $32768$ elements. Figure \ref{SAW128} shows the results obtained at the final time $t_{\mathrm{end}}= 1$. We observe a good agreement with the reference solution for density, pressure and velocity variables along the 1D cut in $y=0$. The Mach number contours are also shown in Figure \ref{SAW128}. It is important to note that in this test the Mach number is quite low, but compressibility plays still an important role in this test problem. One can even clearly see the steepening of the acoustic wave front. Despite our CFL condition being based on the flow velocity, and not on the sound speed, {\color{black} taking $\mathrm{CFL}_{\mathbf{u}}=0.1$} the code still succeeds in properly capturing the position of the travelling acoustic wave. {\color{black} The solution has been obtained after $365$ time steps with a total CPU time of $134.25$s}

\subsection{2D circular explosion}\label{sec:CE}
\begin{figure}[h]
	\centering
	\includegraphics[trim= 5 0 5 0,clip,width=0.45\linewidth]{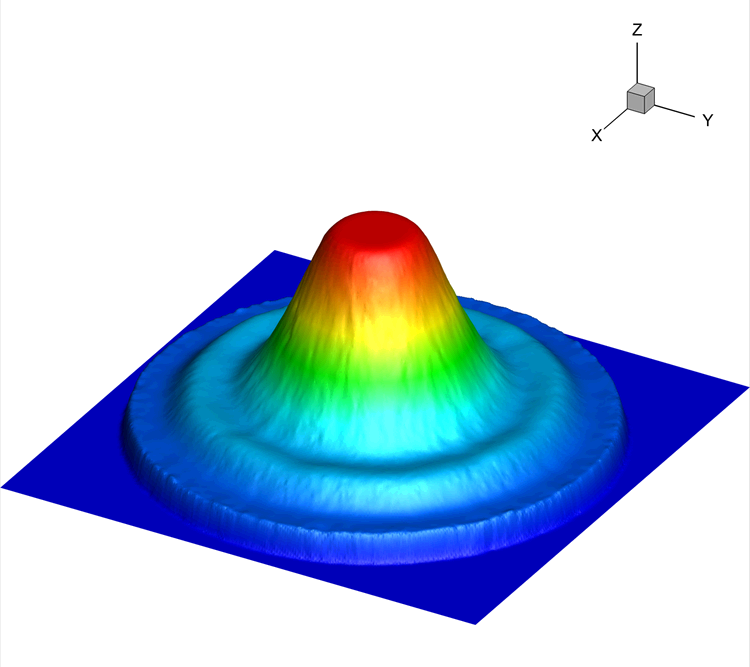}\hspace{0.05\linewidth}
	\includegraphics[trim= 5 0 5 0,clip,width=0.45\linewidth]{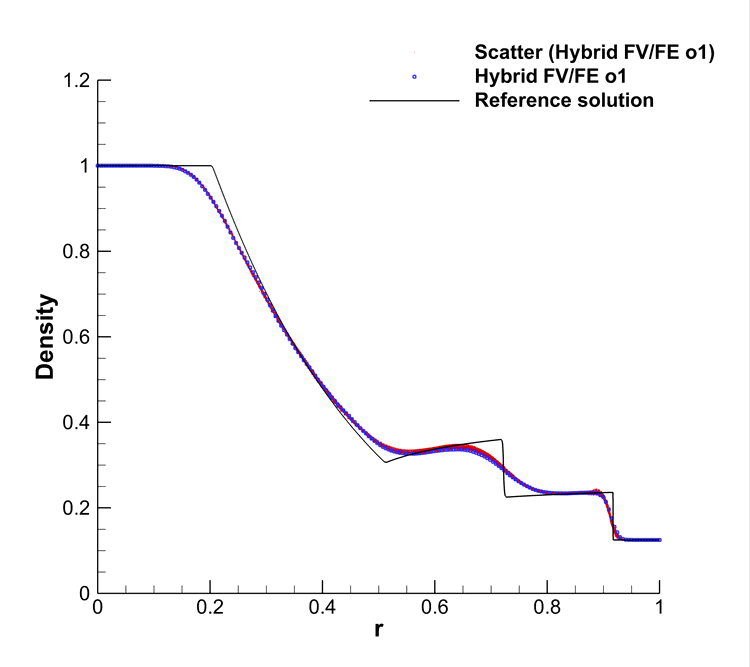}
	
	\vspace{0.05\linewidth}
	\includegraphics[trim= 5 0 5 0,clip,width=0.45\linewidth]{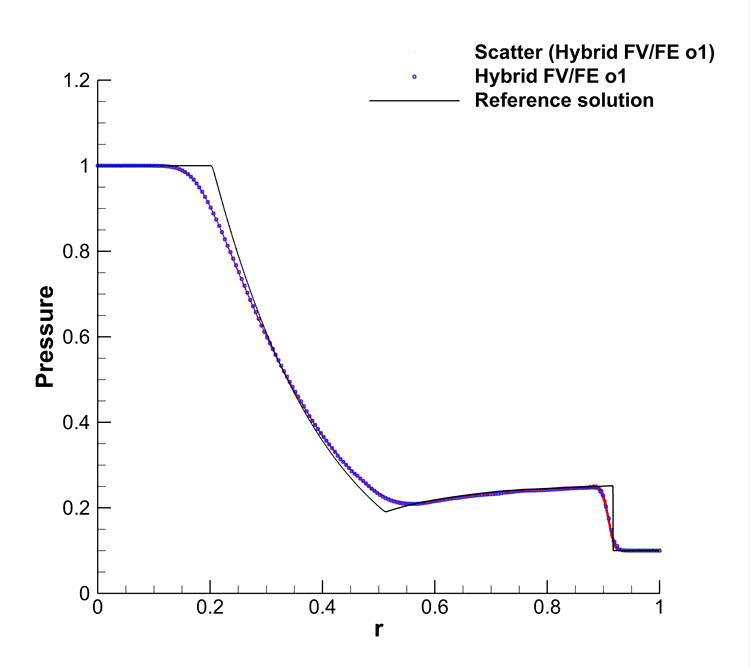}\hspace{0.05\linewidth}
	\includegraphics[trim= 5 0 5 0,clip,width=0.45\linewidth]{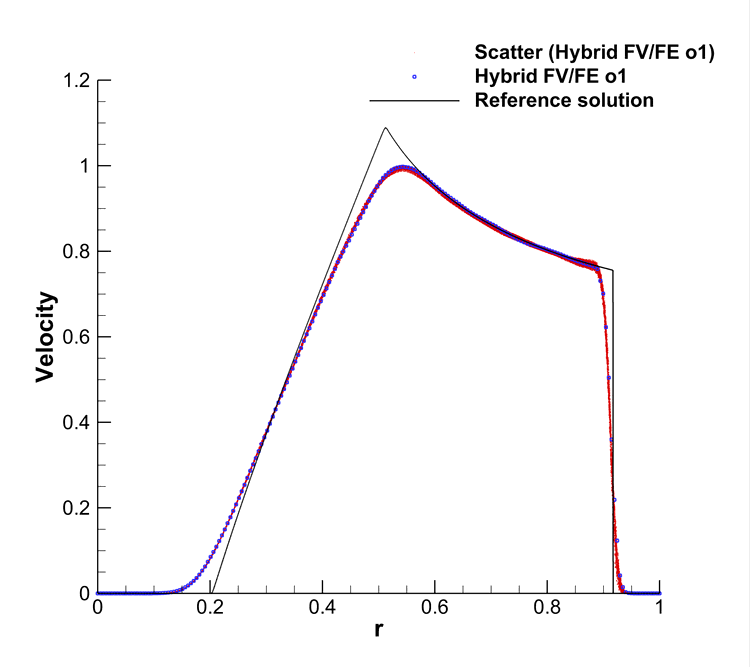}
	\caption{Circular explosion. Solution obtained at $t_{\mathrm{end}}=0.25$ using the first order scheme {\color{black} with $\mathrm{CFL}_{c}=0.83$}. From left top to right bottom: 3D plot of the density. Scatter (red dots), 1D profile (blue squares) and reference solution (black line) of density, pressure and velocity magnitude.}
	\label{CE20_o1_t025}
\end{figure}

\begin{figure}[h]
	\centering
	\includegraphics[trim= 5 0 5 0,clip,width=0.45\linewidth]{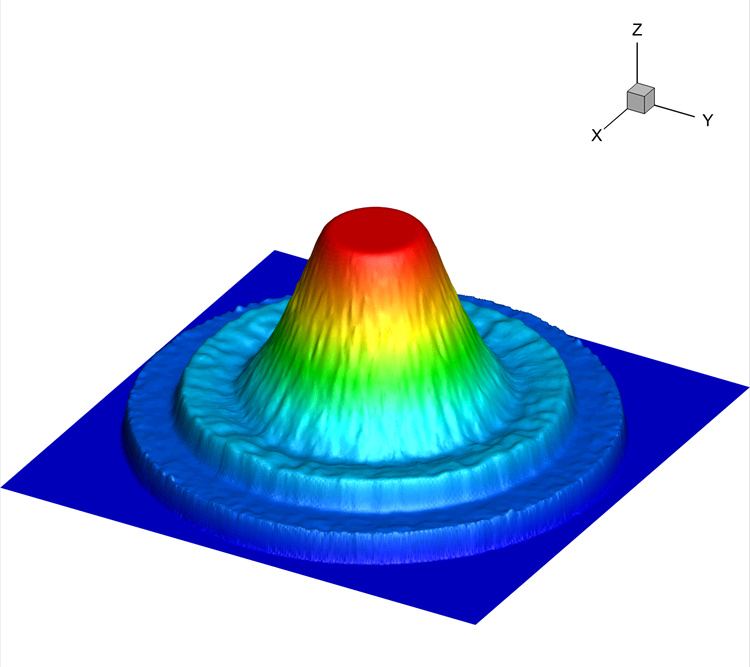}\hspace{0.05\linewidth}
	\includegraphics[trim= 5 0 5 0,clip,width=0.45\linewidth]{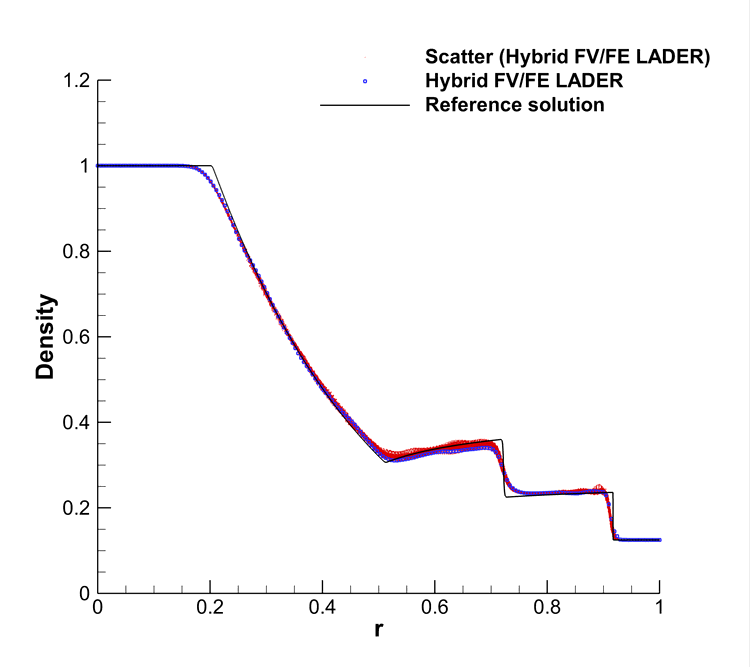}
	
	\vspace{0.05\linewidth}
	\includegraphics[trim= 5 0 5 0,clip,width=0.45\linewidth]{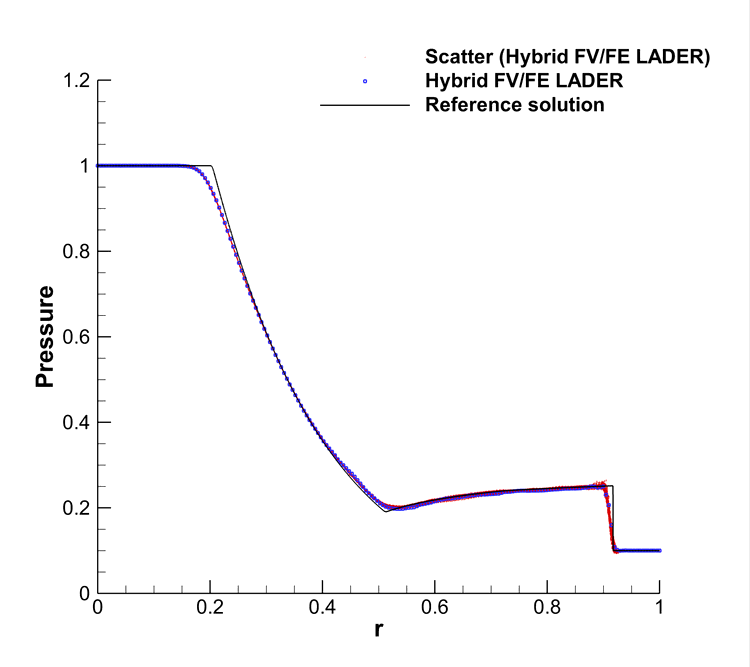}\hspace{0.05\linewidth}
	\includegraphics[trim= 5 0 5 0,clip,width=0.45\linewidth]{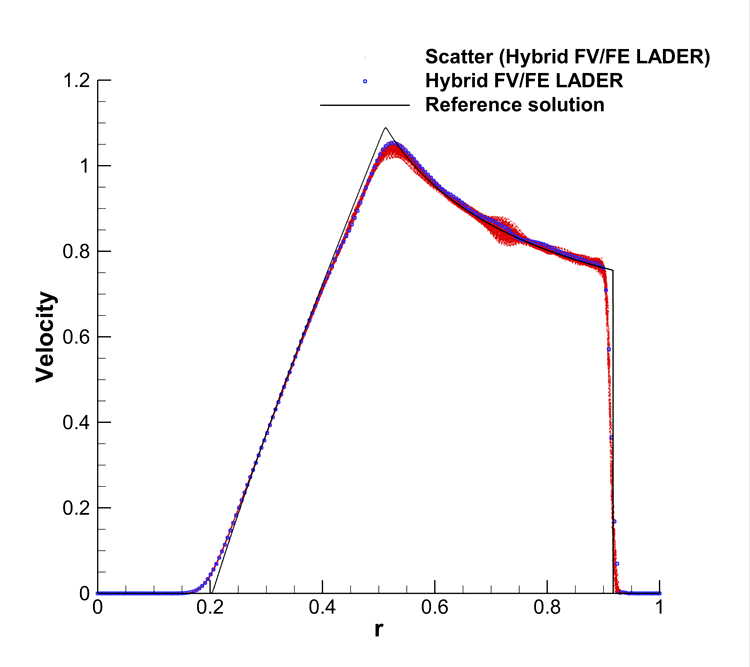}
	\caption{Circular explosion. Solution obtained at $t_{\mathrm{end}}=0.25$ using LADER-ENO {\color{black} with $\mathrm{CFL}_{c}=0.83$}. From left top to right bottom: 3D plot of the density. Scatter (red dots), 1D profile (blue squares) and reference solution (black line) of density, pressure and velocity magnitude.}
	\label{CE20_Le_t025}
\end{figure}

The seventh test consists in a two-dimensional circular explosion problem (see \cite{Toro}, \cite{TT05}, \cite{DPRZ16}) with initial condition

\begin{equation}
	\rho\left(x,y,0\right) =  \left\lbrace \begin{array}{lr}
		1 & \mathrm{ if } \; r \le 0.5,\\
		0.125 & \mathrm{ if } \; r > 0.5,
	\end{array}\right. \qquad
	\press \left(x,y,0\right) = \left\lbrace \begin{array}{lr}
		1 & \mathrm{ if } \; r \le 0.5,\\
		0.1 & \mathrm{ if } \; r > 0.5,
	\end{array}\right. \qquad
	\mathbf{u} \left(x,y,0\right) = 0,
\end{equation}
defined on the computational domain $\Omega=[-1,1]\times[-1,1]$. Moreover, we consider $\mu = \lambda = 0$ and periodic boundary conditions. The primal mesh used in the simulation consists of 
$85344$ primal elements. 
In the same way as for the smooth acoustic wave problem in Section \ref{sec:SAW}, this test case can be compared against a 1D reference solution obtained from solving, in radial direction, the compressible Euler equations with appropriate geometrical source terms (see again \cite{Toro}) on a very fine mesh of $10000$ elements. For the 2D simulations, two different schemes have been considered: a first order scheme and the second order LADER scheme using ENO reconstruction. In both cases we have taken an auxiliary artificial viscosity of $c_{\alpha}=1$. The comparison between the 2D solution along a 1D cut and the 1D reference solution, for $t_{\mathrm{end}}=0.25$, is portrayed in Figure \ref{CE20_o1_t025} for the first order scheme and in Figure \ref{CE20_Le_t025} for LADER-ENO methodology. We observe a good agreement for all flow variables. {\color{black} Let us remark that this test case has a characteristic Mach number above unity, namely $M\approx 1.2$. The CPU time of the simulations is reported in Table  \ref{tab:CE_Le_t025_times}.}

\begin{table}[H]
	\renewcommand{\arraystretch}{1.2}
	\begin{center}
		\begin{tabular}{|c|c||c|}
			\hline
			\multirow{3}{*}{Order 1} & CPU time (s)&
			$420.78$ 
			\\\hhline{|~|-||-|}
			& $t_{e}$ ($\mu$s) &			
			$5.74$ 
			\\\hhline{|~|-||-|}
			& Time steps &
			$570$ 
			\\\hline
			\multirow{3}{*}{LADER} & CPU time (s)&
			$595.45$ 
			\\\hhline{|~|-||-|}
			& $t_{e}$ ($\mu$s) &			
			$7.87$ 		
			\\\hhline{|~|-||-|}
			& Time steps &
			$589$ 
			\\\hline
		\end{tabular}
		{\color{black} \caption{Circular explosion.
				CPU time, CPU time per element and iteration, $t_{e}$, and number of time steps.
			}\label{tab:CE_Le_t025_times}}
	\end{center}
\end{table}

\subsection{Heat conduction}\label{sec:HC}
\begin{figure}[h]
	\centering
	\includegraphics[width=0.325\linewidth]{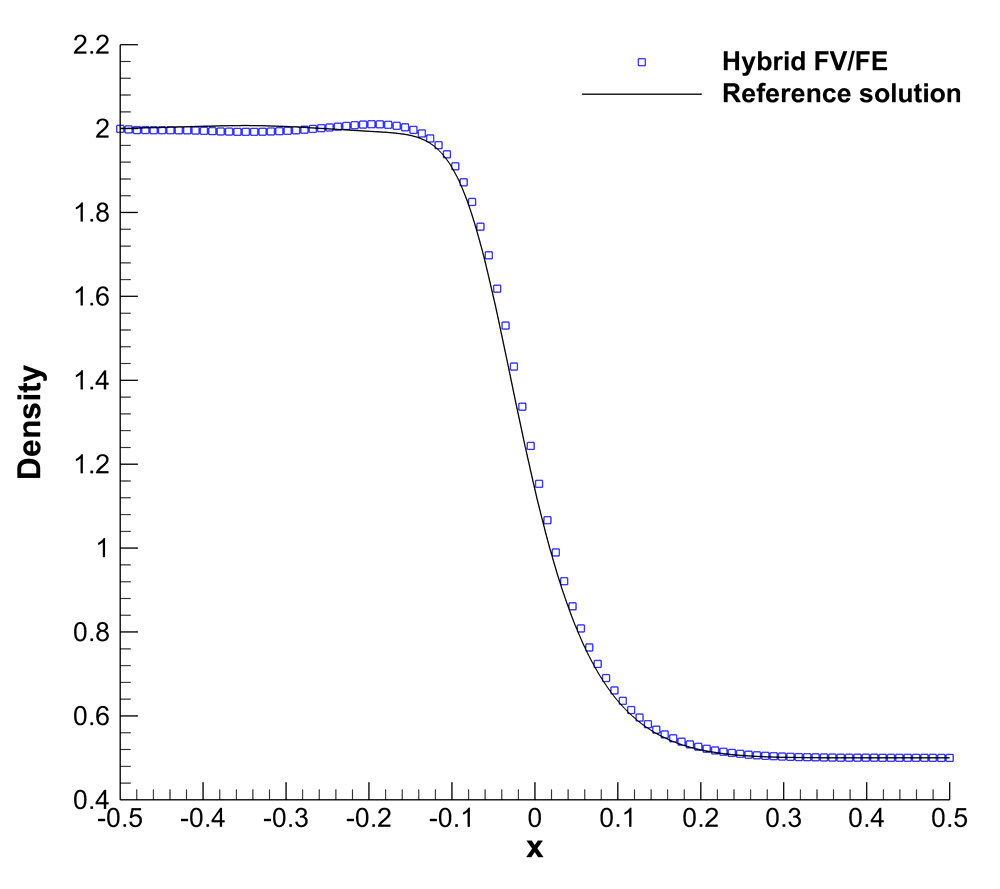}
	\includegraphics[width=0.325\linewidth]{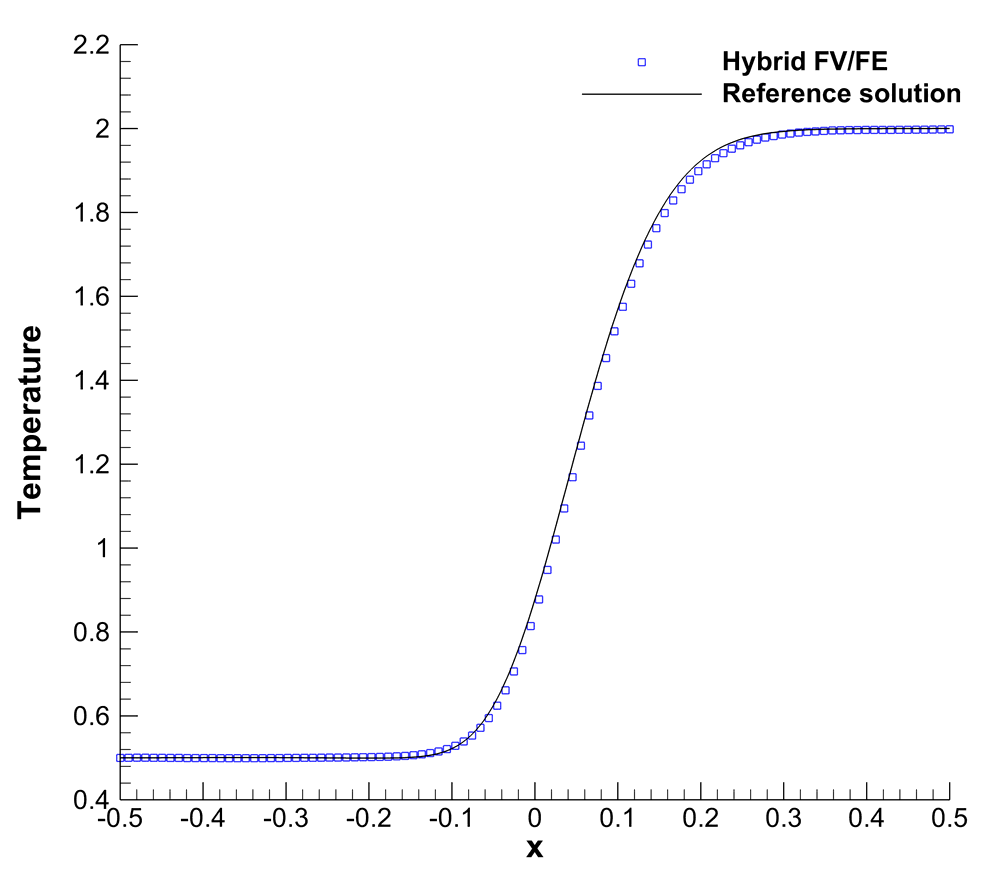}
	\includegraphics[width=0.325\linewidth]{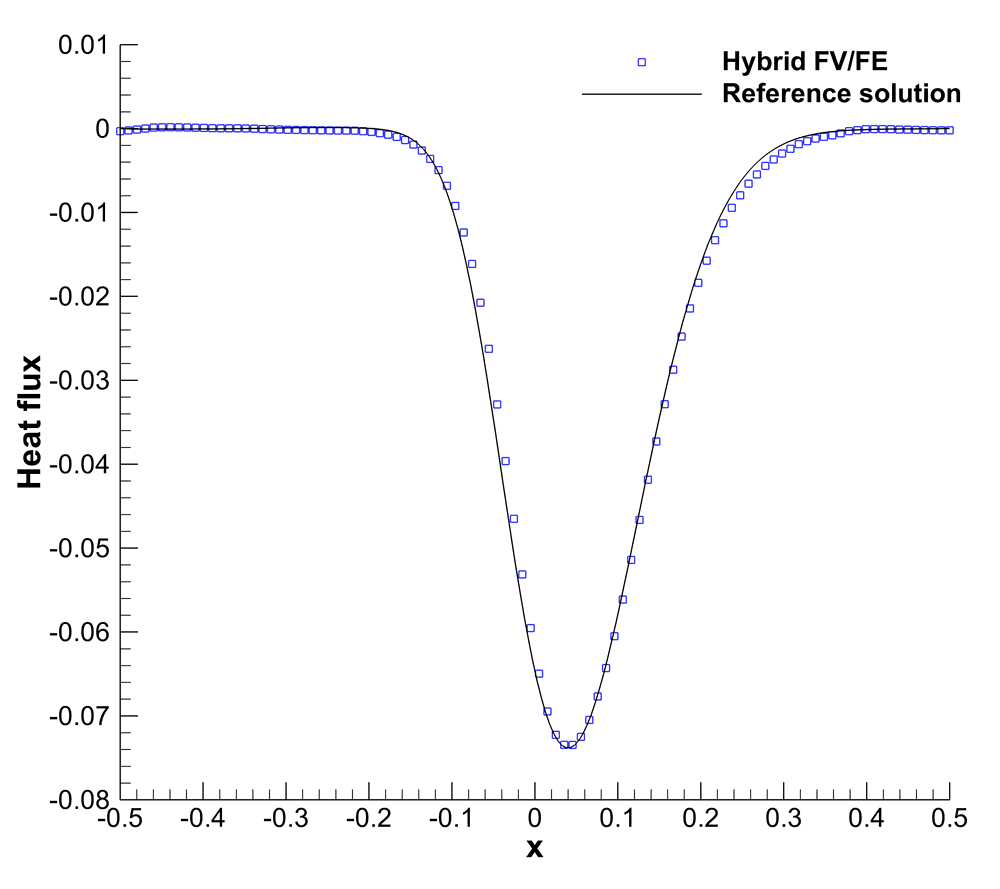}
	\caption{Heat conduction at time  $t_{\mathrm{end}}=1$ {\color{black} ($\mathrm{CFL}_{c}=277.8$)}. From left to right: density, temperature and heat flux, $q_{1}=-\kappa \partial_{x}\theta$, along the cut $y=0$.}
	\label{HC_figure}
\end{figure}

Following \cite{DPRZ16}, we define a heat conduction dominated test problem with initial condition

\begin{equation}
	\rho\left(x,y,0\right) =  \left\lbrace \begin{array}{lr}
		2 & \mathrm{ if } \; x \le 0,\\
		0.5 & \mathrm{ if } \; x > 0,
	\end{array}\right. \qquad
	\press \left(x,y,0\right) = 1 \qquad
	\mathbf{u} \left(x,y,0\right) = 0,
\end{equation} 
on $\Omega=[-0.5,0.5]\times [-0.1,0.1]$. Moreover, the fluid properties are $\gamma = 1.4$, $\mu= \lambda = 10^{-2}$, $c_{\press}=3.5$, $M\approx 1.8\cdot 10^{-2}$. Dirichlet boundary conditions for the velocity and density and Neumann boundary conditions for the pressure are set on $x$-direction, whereas in $y$-direction we consider periodic boundary conditions. 
The simulation is run on a mesh of $1000$ primal elements using the LADER scheme.
The density, temperature and heat flux obtained are portrayed in Figure \ref{HC_figure}. We also include the reference solution obtained in \cite{DPRZ16} by solving the 1D compressible Navier-Stokes equations on a very fine mesh. We observe that even the heat flux, which involves derivatives of the temperature, is in good agreement with the reference solution. 

\subsection{Rising bubble problem}\label{sec:GBP}

\begin{figure}[h]
	\centering
	\includegraphics[trim= 5 0 5 0,clip,width=0.3\linewidth]{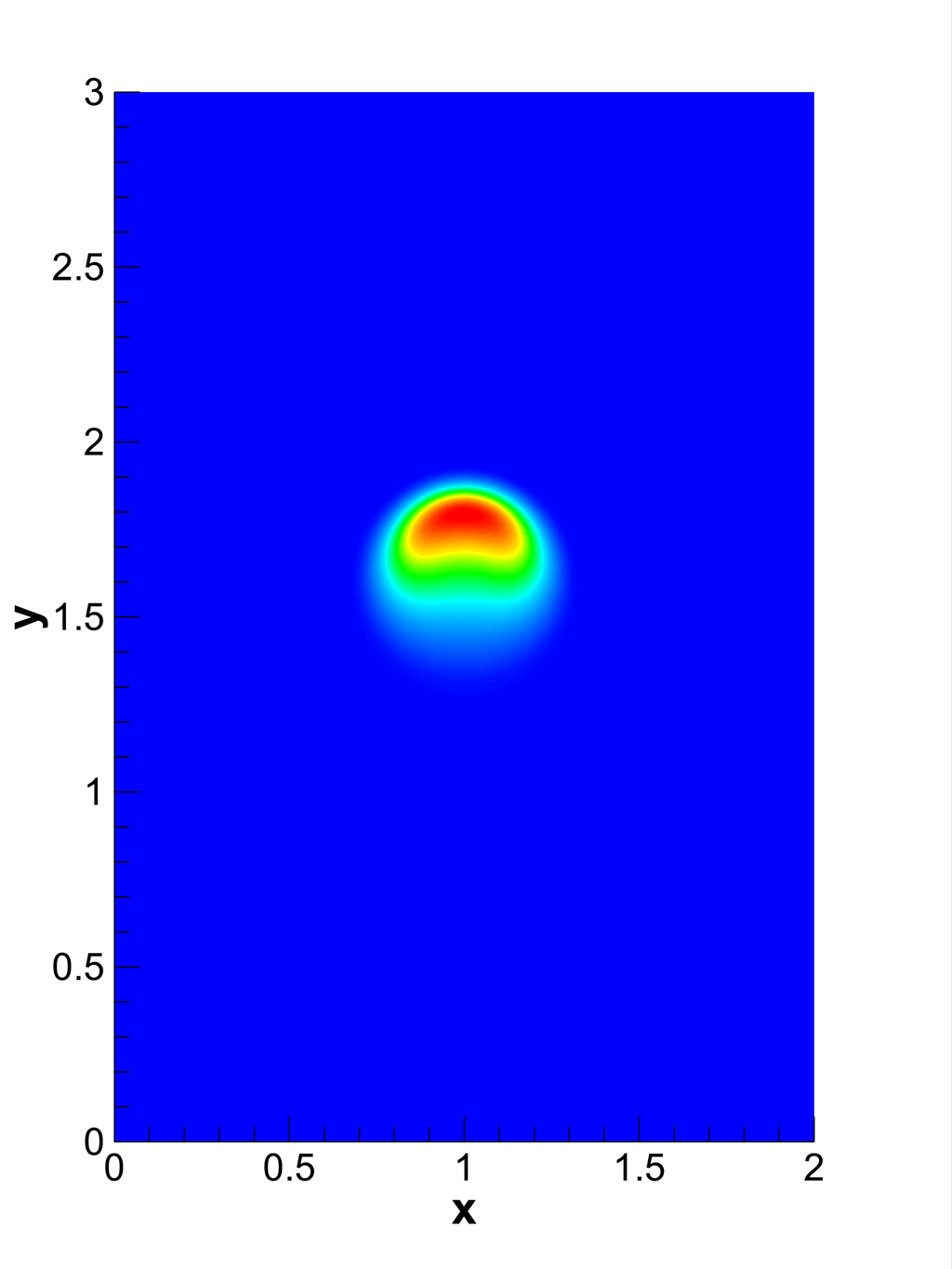}\hspace{0.05\linewidth}
	\includegraphics[trim= 5 0 5 0,clip,width=0.3\linewidth]{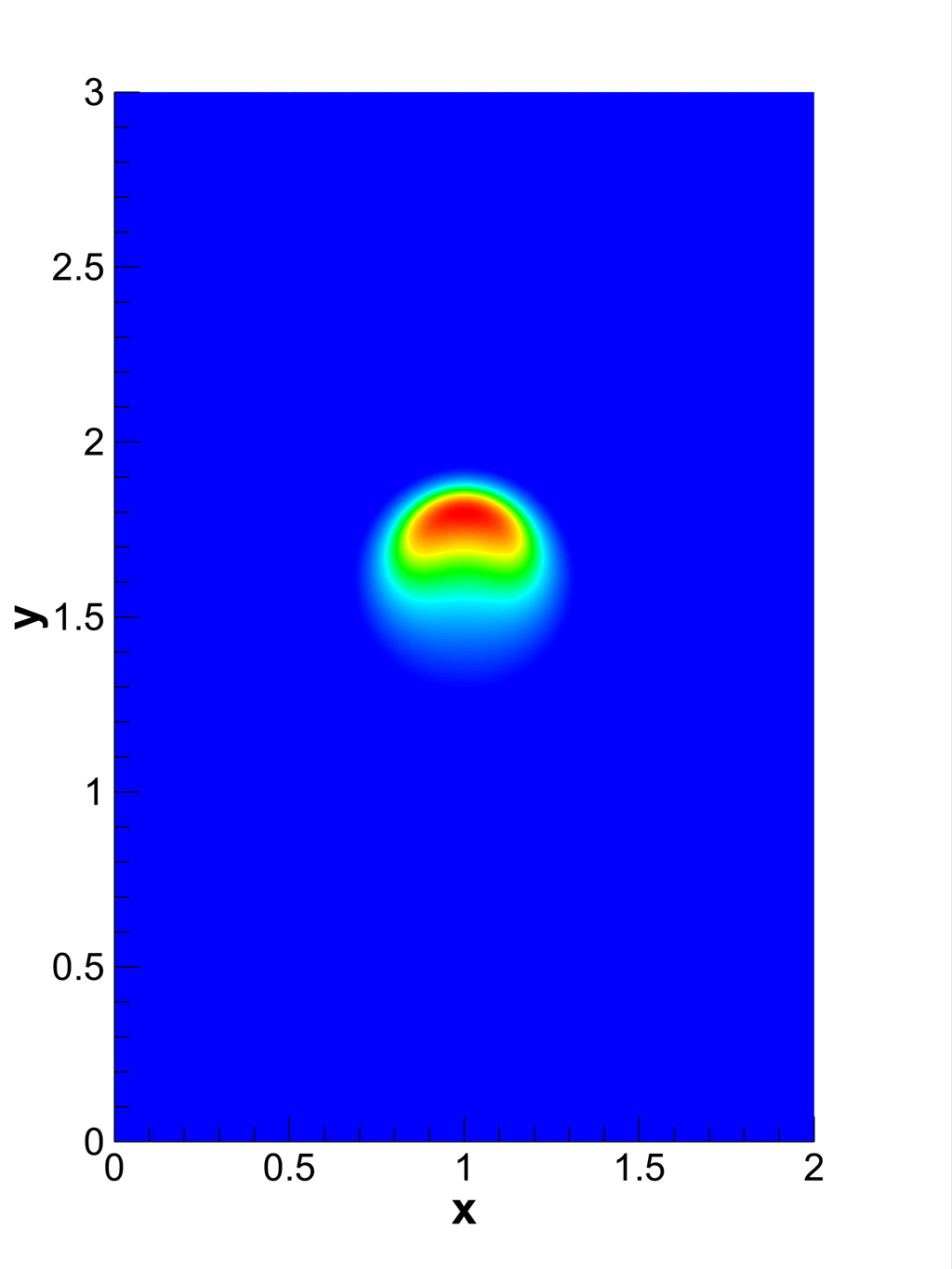}
	\phantom{\begin{minipage}{0.08\linewidth}
			\vspace{-4\linewidth} \includegraphics[width=\linewidth]{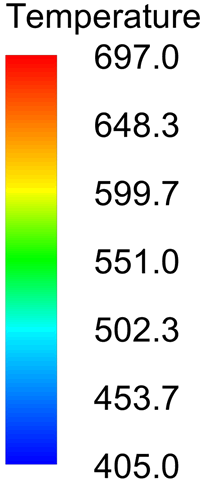}
	\end{minipage}}
	
	\includegraphics[trim= 5 0 5 0,clip,width=0.3\linewidth]{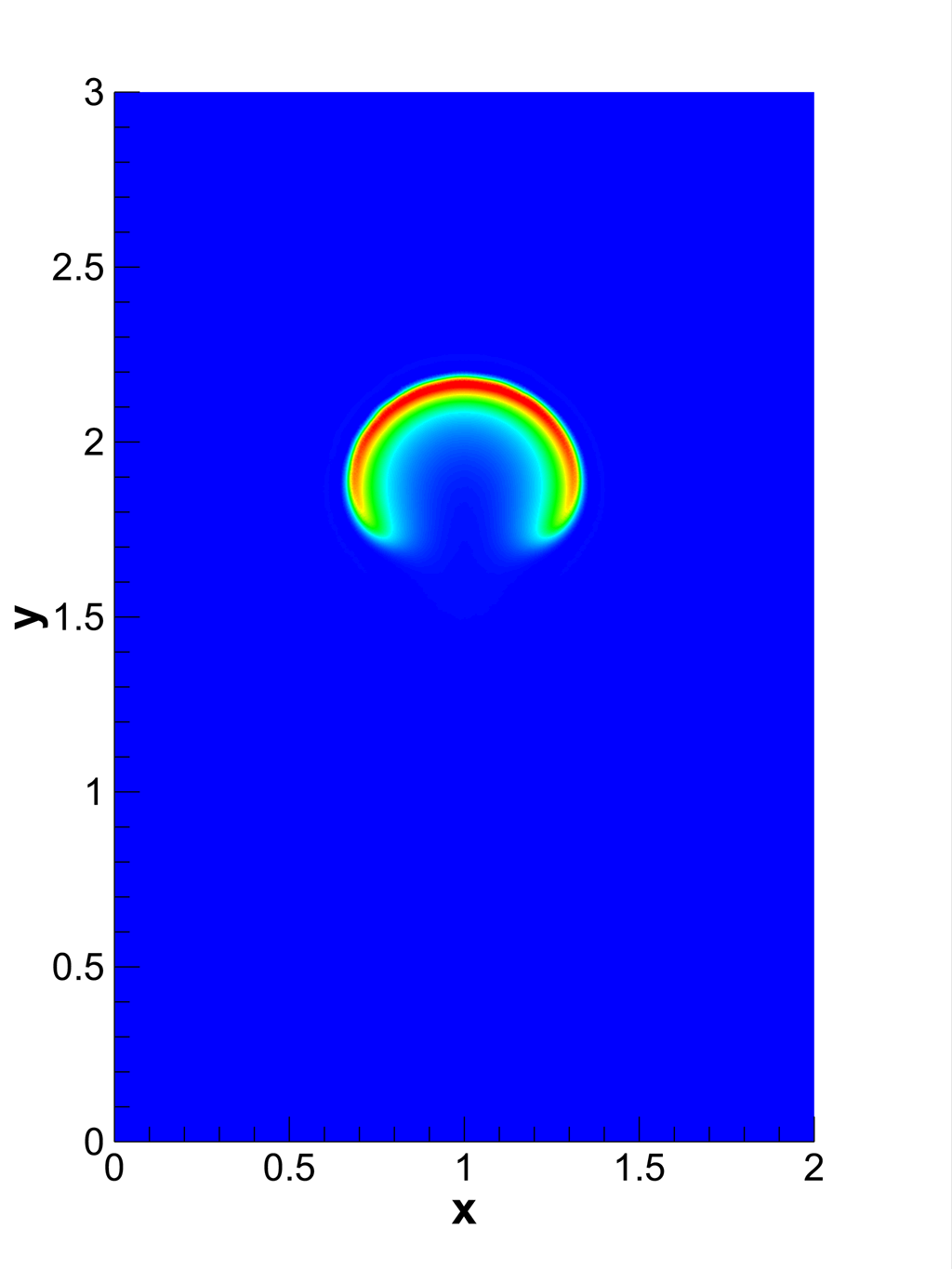}\hspace{0.05\linewidth}
	\includegraphics[trim= 5 0 5 0,clip,width=0.3\linewidth]{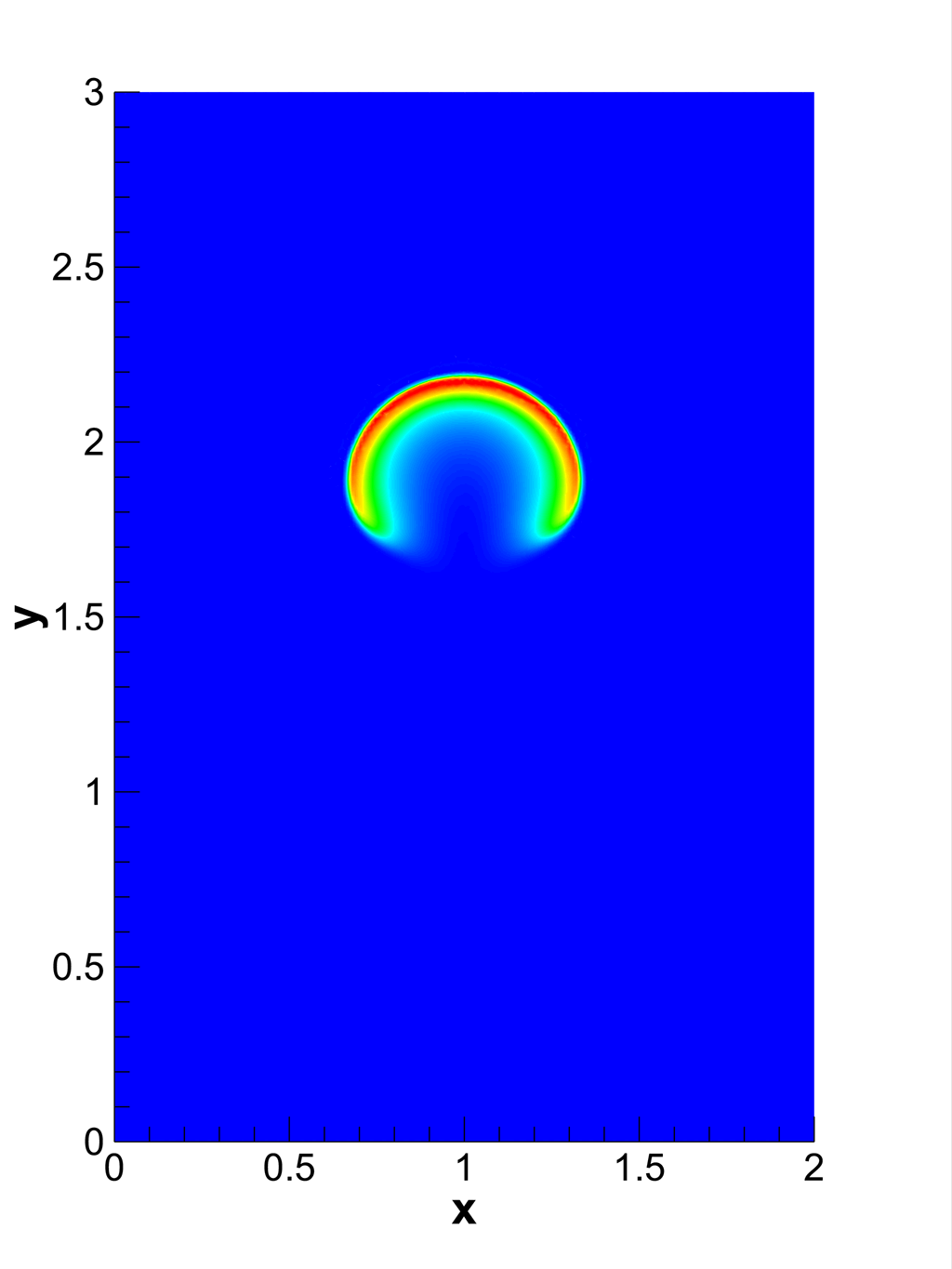}
	\begin{minipage}{0.08\linewidth}
		\vspace{-4\linewidth} \includegraphics[width=\linewidth]{figures/GP_temperature_vlegend.png}
	\end{minipage}
	
	\includegraphics[trim= 5 0 5 0,clip,width=0.3\linewidth]{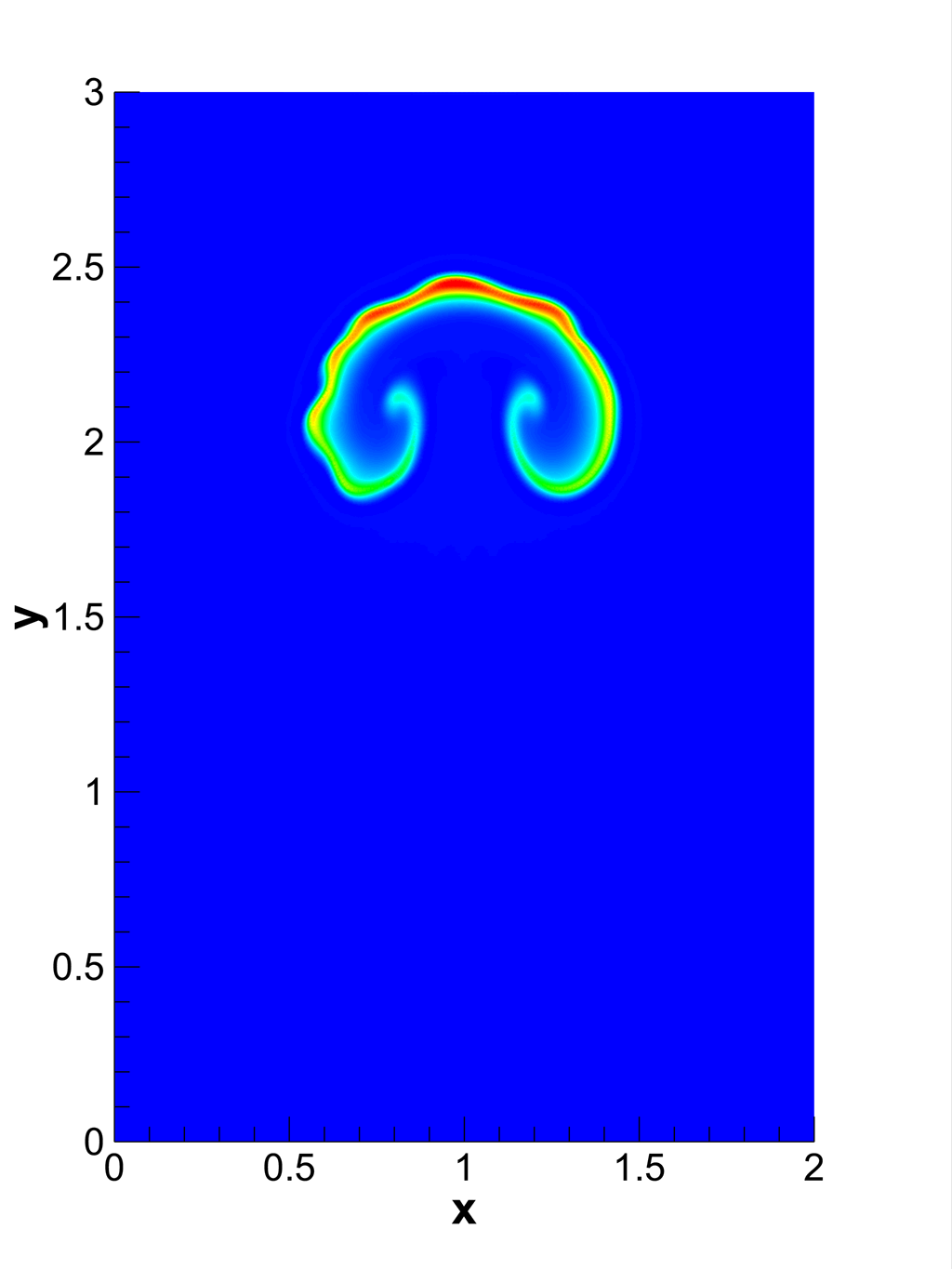}\hspace{0.05\linewidth}
	\includegraphics[trim= 5 0 5 0,clip,width=0.3\linewidth]{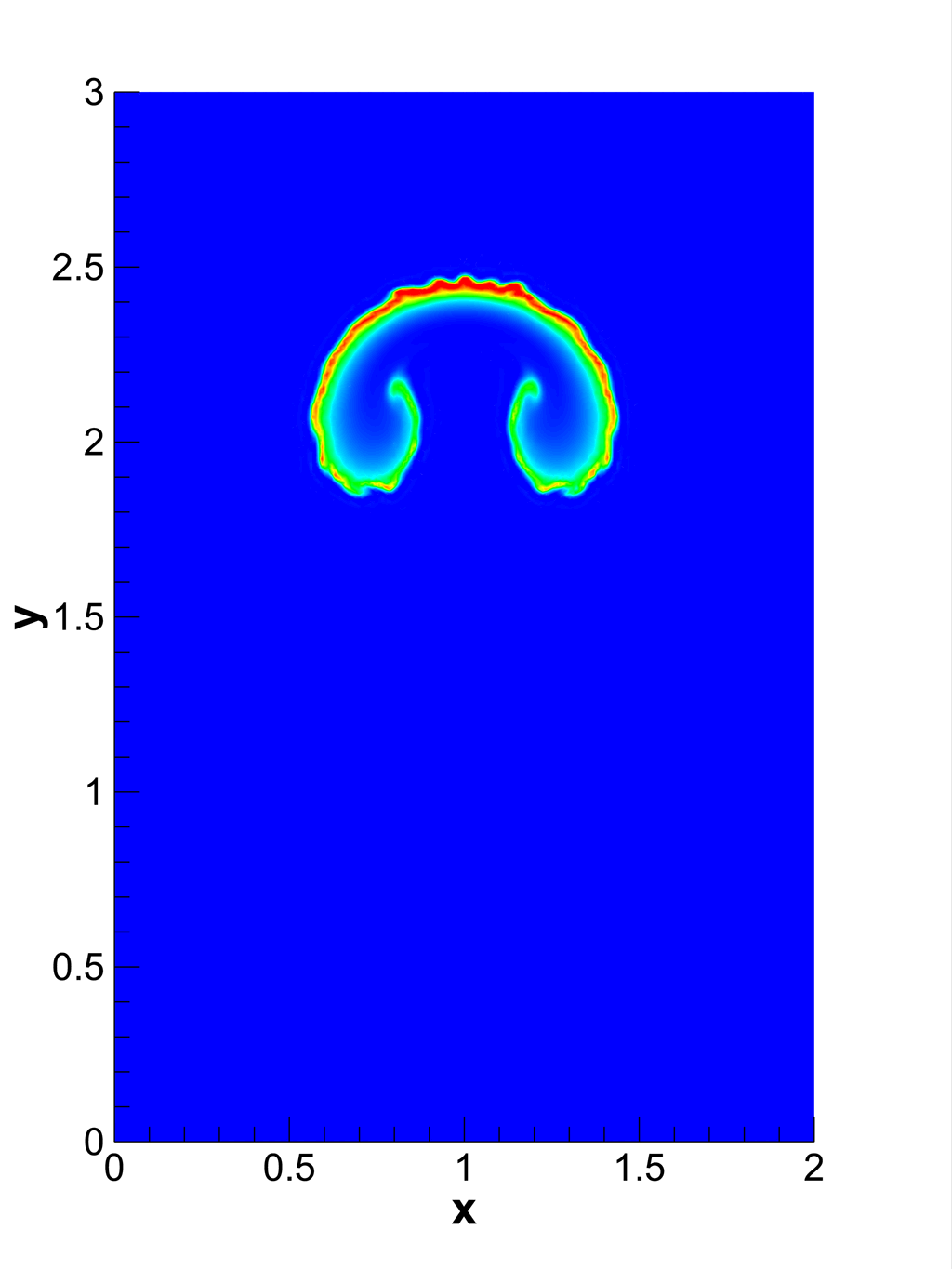}
	\phantom{\begin{minipage}{0.08\linewidth}
			\vspace{-4\linewidth} \includegraphics[width=\linewidth]{figures/GP_temperature_vlegend.png}
	\end{minipage}}
	
	\caption{Rising bubble. Temperature contour plot for $t\in\left\lbrace 0.5,1,1.5\right\rbrace$. Left: Hybrid FV-FE method (LADER-ENO). Right: Semi-implicit DG scheme ($p=3$).}
	\label{GP_temperature}
\end{figure}

\begin{figure}[h]
	\centering
	\includegraphics[trim= 5 0 5 0,clip,width=0.3\linewidth]{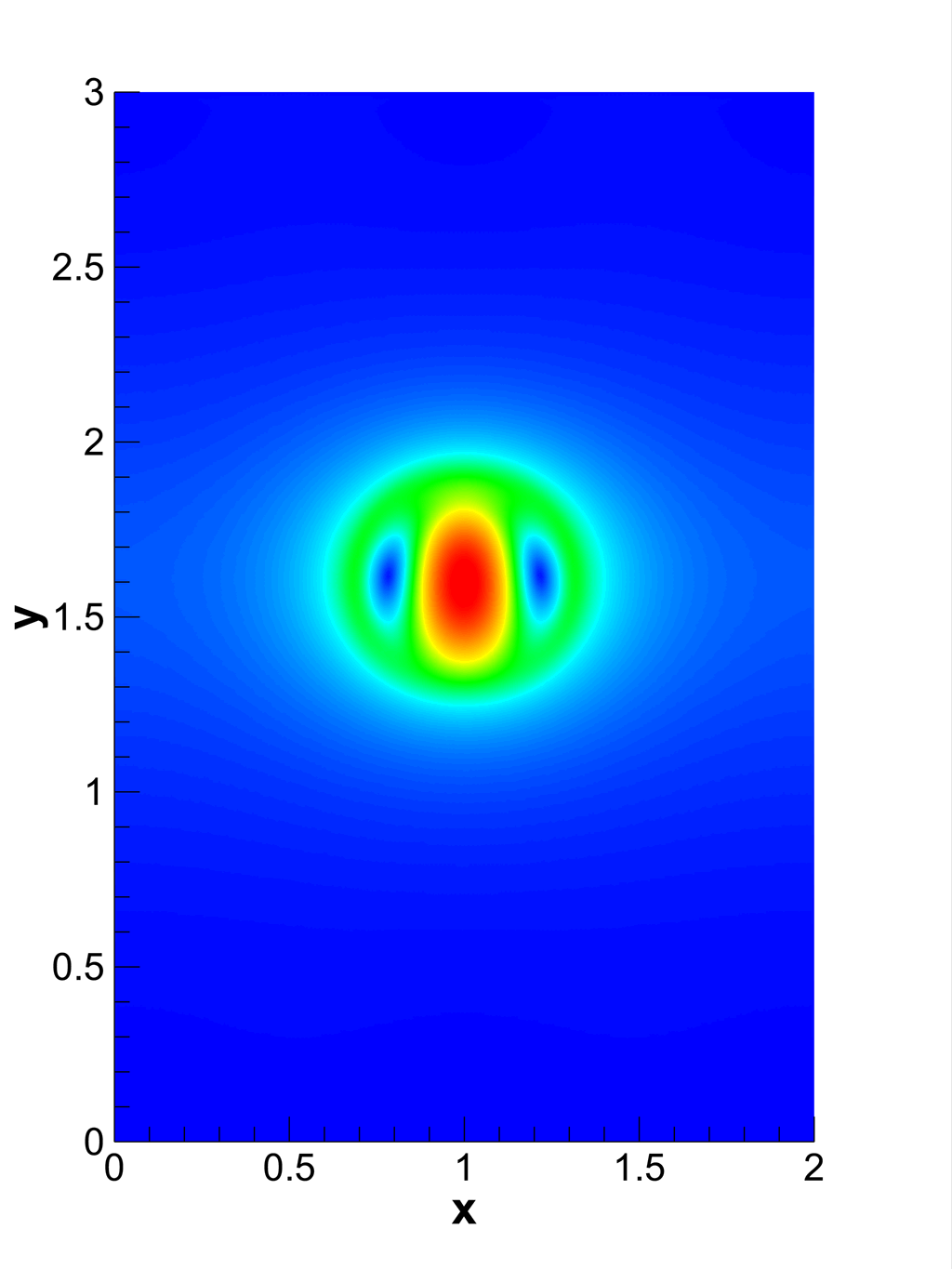}\hfill
	\includegraphics[trim= 5 0 5 0,clip,width=0.3\linewidth]{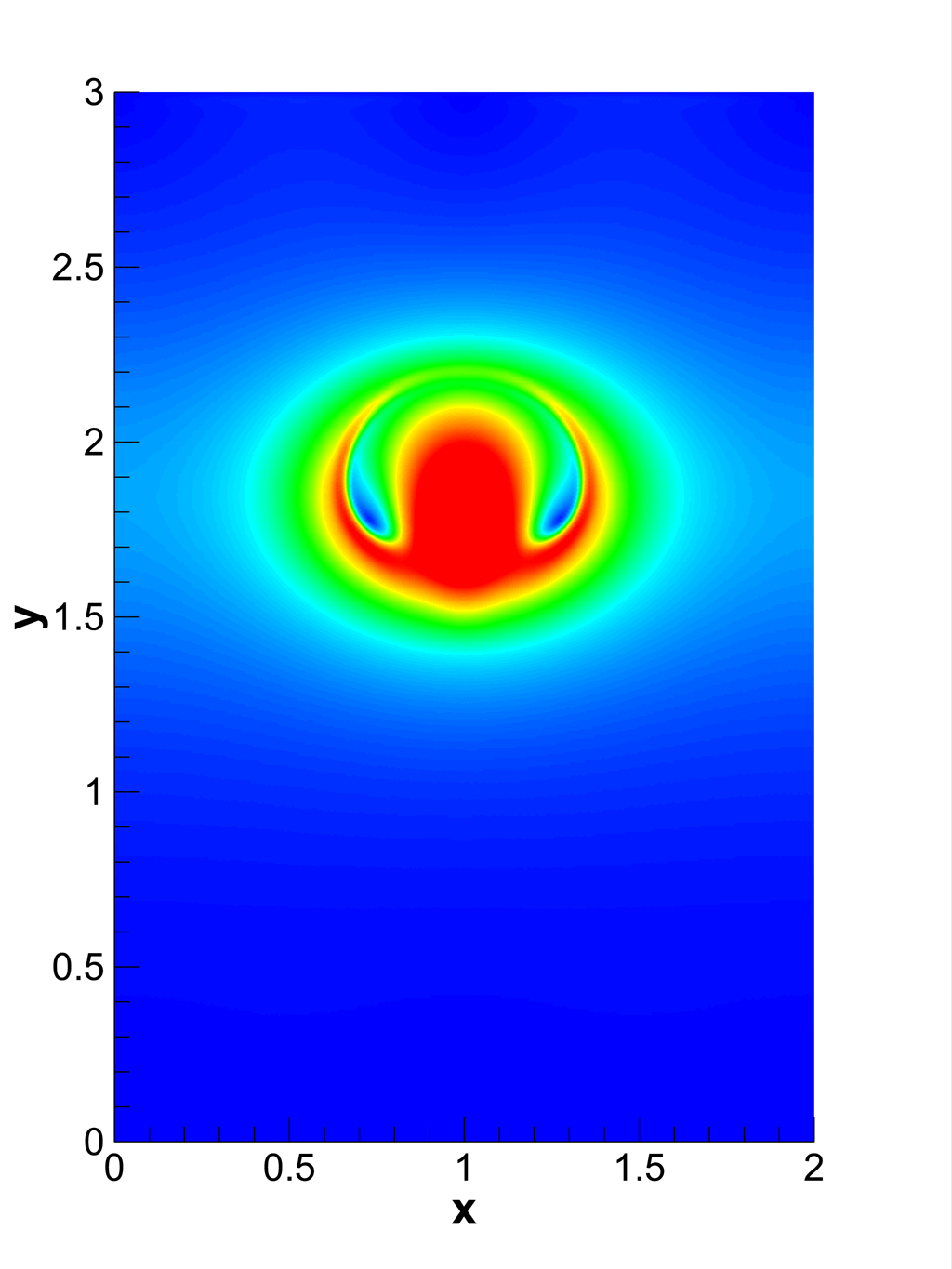}\hfill
	\includegraphics[trim= 5 0 5 0,clip,width=0.3\linewidth]{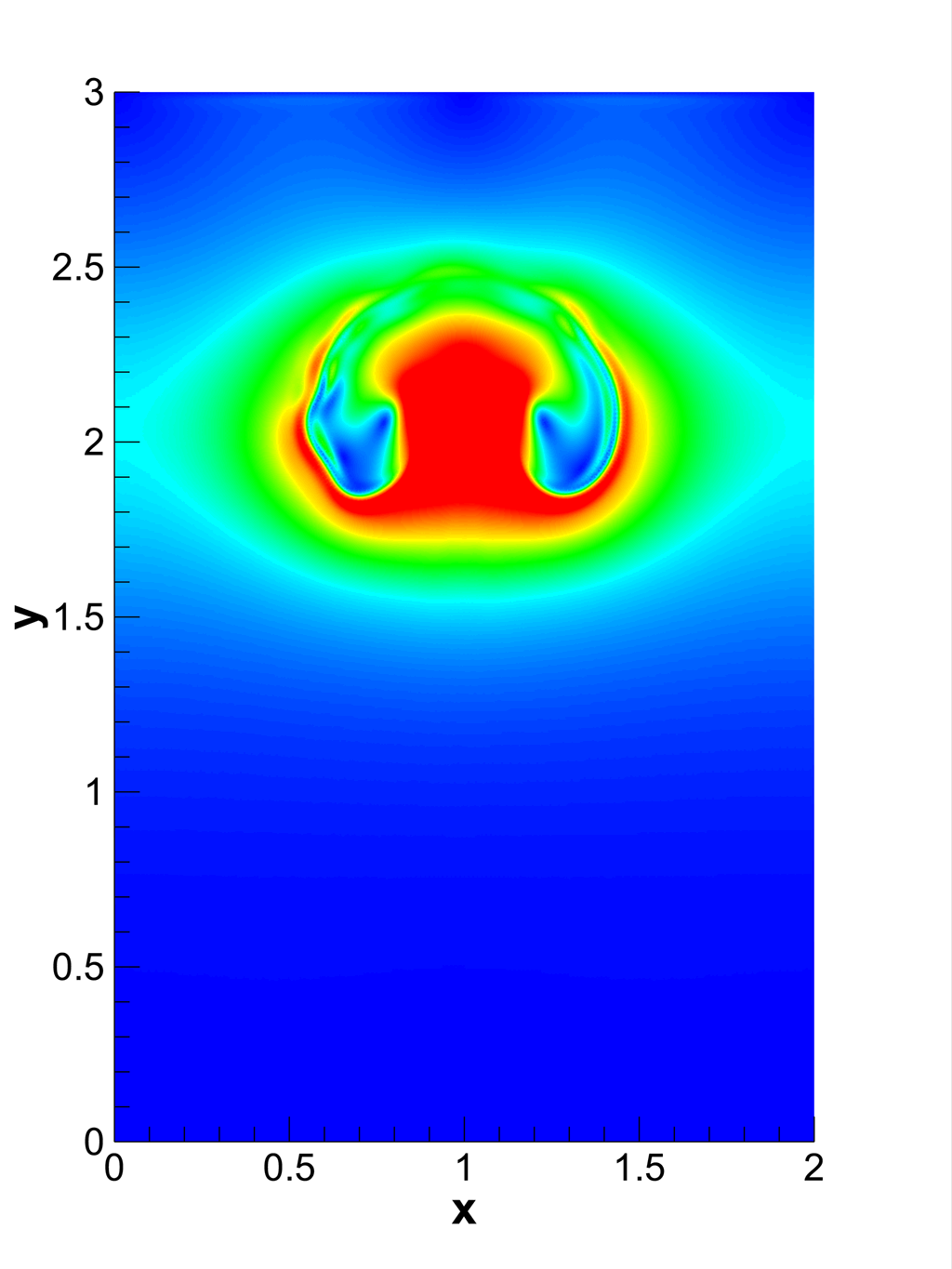}
	
	\includegraphics[trim= 5 0 5 0,clip,width=0.3\linewidth]{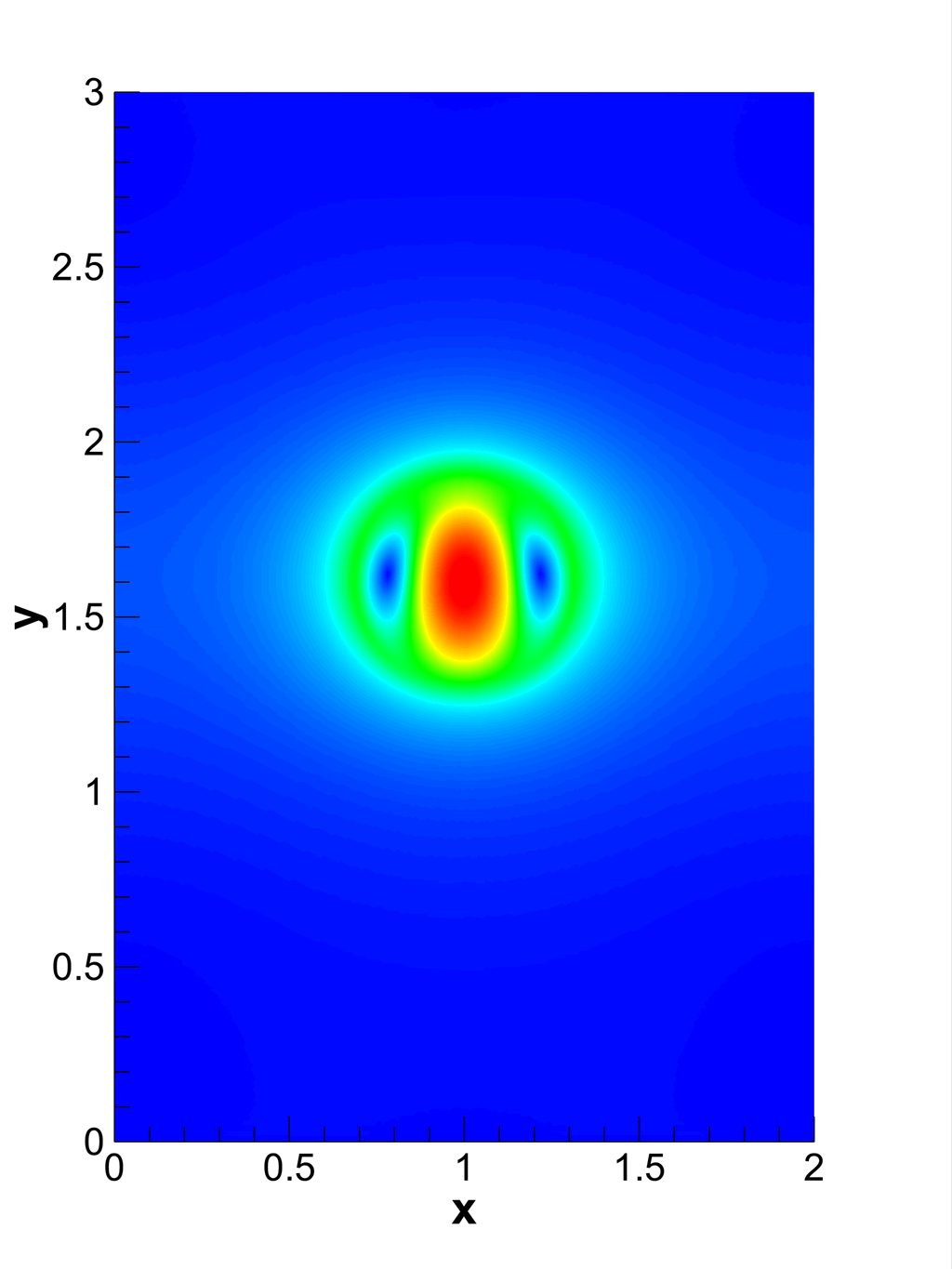}\hfill
	\includegraphics[trim= 5 0 5 0,clip,width=0.3\linewidth]{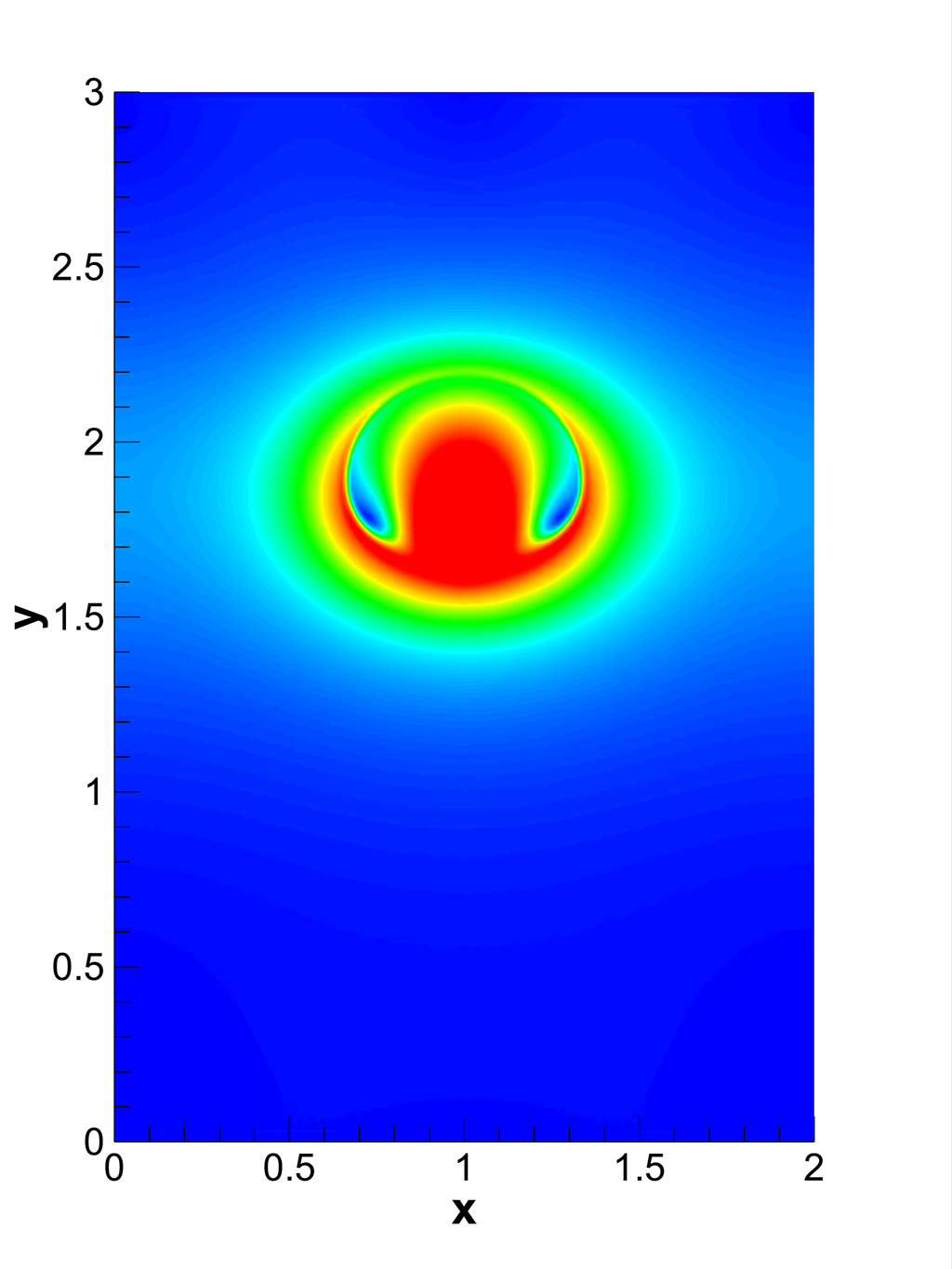}\hfill
	\includegraphics[trim= 5 0 5 0,clip,width=0.3\linewidth]{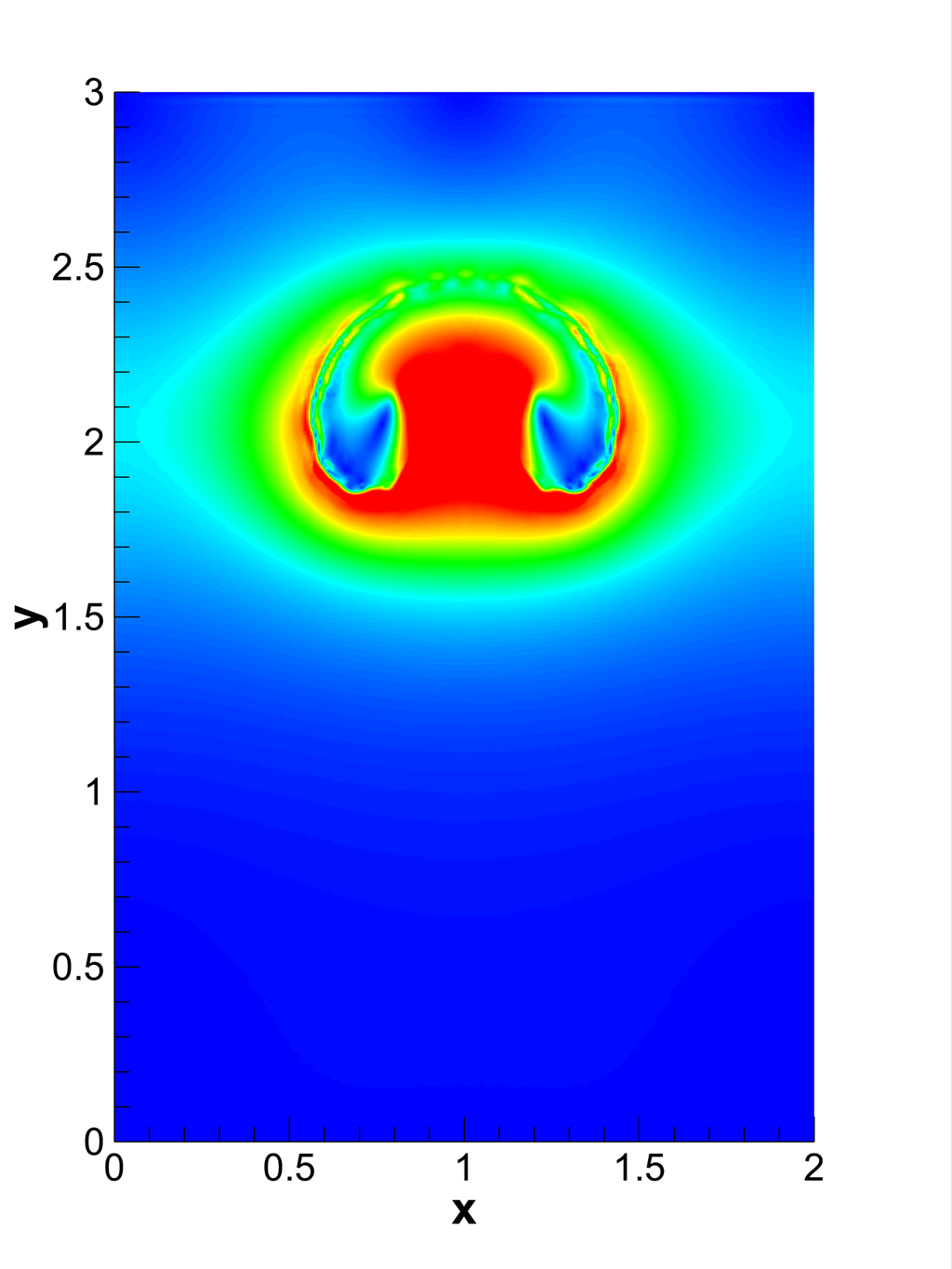}
	
	\includegraphics[width=0.4\linewidth]{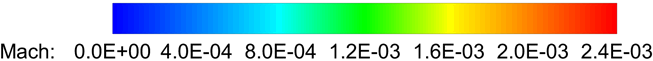}
	\caption{Rising bubble. Mach number contour plot for $t\in\left\lbrace 0.5,1,1.5\right\rbrace$. Top: Hybrid FV-FE method (LADER-ENO). Bottom: Semi-implicit DG scheme ($p=3$).}
	\label{GP_mach}
\end{figure}

Rising bubble benchmarks are typically used for the assessment of thermal convection-driven problems (see \cite{MBGW13,YMLGW14,BLY17,Yi18,BTBD20}). Within this article we have considered an initial Gaussian bubble of the form 

\begin{equation}
	\rho\left(x,y,0\right) =  \left\lbrace \begin{array}{lr}
		1-\dfrac{1}{2} e^{-\frac{1}{2} \left( \frac{r}{0.2}\right)^2} & \mathrm{ if } \; r^{2} \le 0.1,\\[10pt]
		1-\dfrac{1}{2} e^{-\frac{5}{4}} & \mathrm{ otherwise},
	\end{array}\right. \qquad
	\press \left(x,y,0\right) = 10^{5}+ yg_{2} \qquad
	\mathbf{u} \left(x,y,0\right) = 0,
\end{equation} 
in the computational domain $\Omega=\left[0,2\right] \times \left[0,3\right]$, with $r =\sqrt{\left(x-1\right)^2+\left(y-1.5\right)^2}$ the radius with respect to the centre of the bubble, $\mathbf{x}=\left(1,1.5\right)$, $\mathbf{g}=\left(0,-9.81\right)^{T}$ the gravity, $\mu=3.36\cdot 10^{-3}$ the viscosity, $\lambda= 0.0238$ the thermal conductivity and $\gamma=1.4$ the adiabatic index. Periodic boundary conditions are set in the  $x$-direction, whereas the top and bottom boundaries are assumed to be adiabatic walls. The temperature variation, under gravity effects, leads to the movement of the fluid. 
The simulation is carried out on a mesh of $85528$ primal elements up to time $t_{\mathrm{end}}=1.5$ {\color{black} with $\mathrm{CFL}_{c}=2083.3$}. 
To validate the results, we have also run the simulation with the high order staggered semi-implicit discontinuous Galerkin scheme presented in \cite{TD14,TD15,TD17,BTBD20}. More precisely, we have taken a polynomial approximation degree of $p=3$ in space, using a primal mesh composed of $20448$ elements. Second order in time is reached thanks to the theta method. The  temperature and Mach number contour plots are depicted in Figures \ref{GP_temperature} and \ref{GP_mach}, respectively. Let us remark that the unstructured and non-symmetric grids employed may lead to a loss of symmetry at large times. For small times, we observe a good agreement of the results obtained with both numerical schemes. At $t=1.5$, we observe that 
both methodologies are able to capture the small instabilities arising on the top of the main structure of the thermal bubble and that would lead to the so called Kelvin-Helmholtz instabilities for larger times. As expected, the shape observed is smoother with the hybrid FV-FE solver, probably due to the effect of the extra artificial viscosity $c_{\alpha}=0.5$ that was imposed for this test. 

\subsection{3D spherical explosion}\label{sec:CE3D}
\begin{figure}[h]
	\centering 	
	\includegraphics[trim= 5 0 5 0,clip,width=0.45\linewidth]{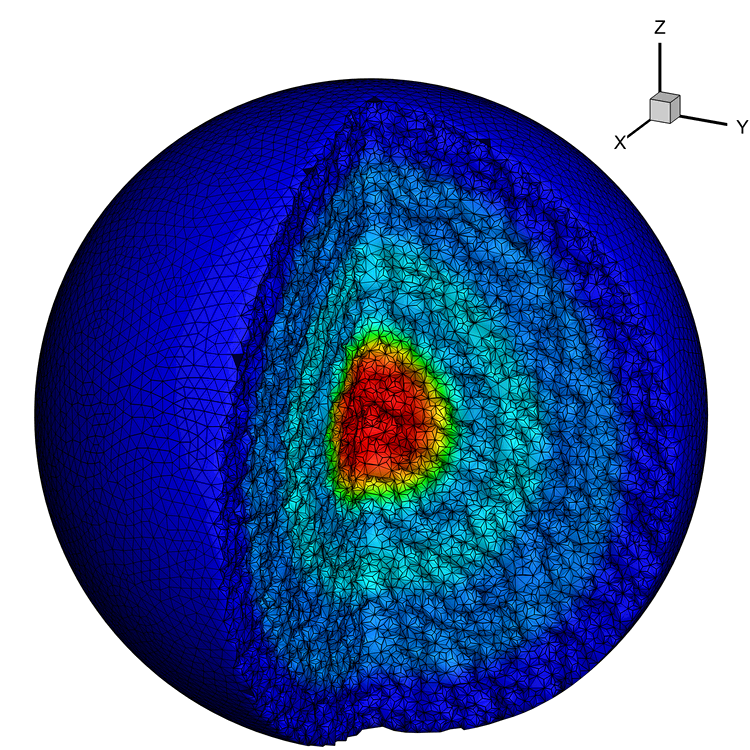}\hspace{0.05\linewidth}
	\includegraphics[trim= 5 0 5 0,clip,width=0.45\linewidth]{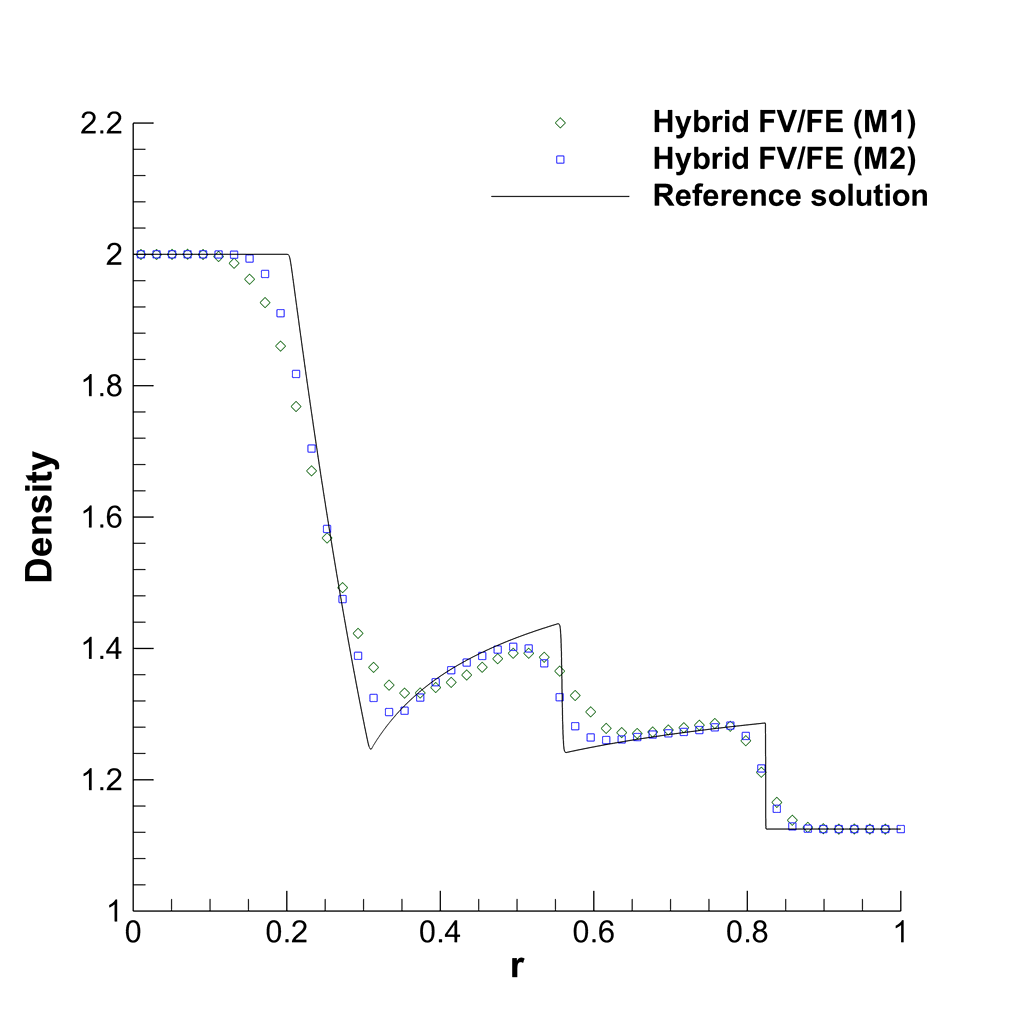}
	
	\vspace{0.05\linewidth}
	\includegraphics[trim= 5 0 5 0,clip,width=0.45\linewidth]{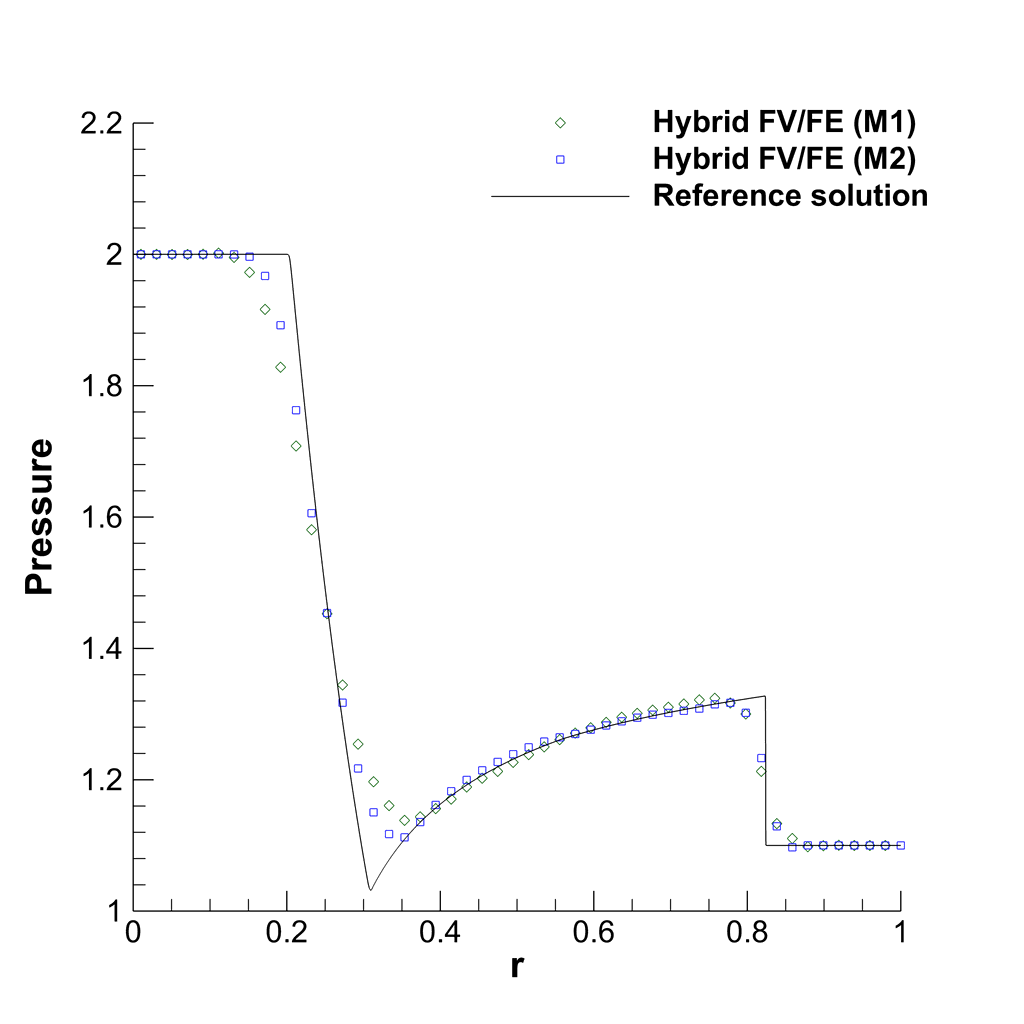}\hspace{0.05\linewidth}
	\includegraphics[trim= 5 0 5 0,clip,width=0.45\linewidth]{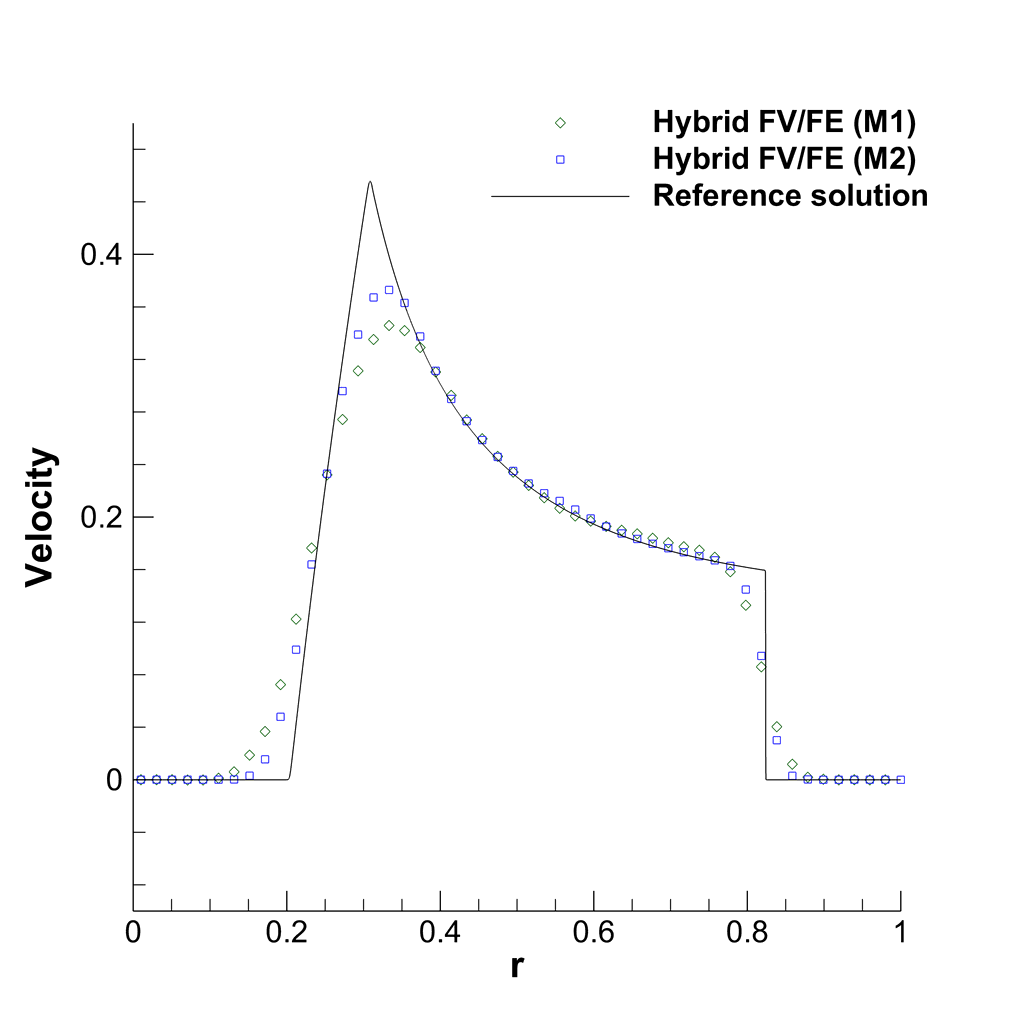}
	\caption{3D spherical explosion. Solution obtained at $t_{\mathrm{end}}=0.25$ using LADER-ENO {\color{black} with $\mathrm{CFL}_{c}=3.3$}. From left top to right bottom: 3D mesh and density contours for M1;
		1D profile for $x\in\left[0,1\right], \, y=0,\, z=0$ (M1 green diamonds, M2 blue squares) and reference solution (black line) of density, pressure and velocity magnitude.}
	\label{fig:CE3d_Le_t025}
\end{figure}

In order to show the capability of the numerical scheme to handle also 3D problems, we present the results obtained for a three-dimensional spherical explosion problem. The sphere of radius $R=1$ centred at the origin is taken as computational domain and initial conditions are given by 

\begin{equation}
	\rho\left(x,y,z,0\right) =  \left\lbrace \begin{array}{lr}
		2 & \mathrm{ if } \; r \le 0.5,\\
		1.125 & \mathrm{ if } \; r > 0.5,
	\end{array}\right. \qquad
	\press \left(x,y,z,0\right) = \left\lbrace \begin{array}{lr}
		2 & \mathrm{ if } \; r \le 0.5,\\
		1.1 & \mathrm{ if } \; r > 0.5,
	\end{array}\right. \qquad
	\mathbf{u} \left(x,y,z,0\right) = 0. \label{eq:IC_CE3D}
\end{equation}
We assume  $\mu = \lambda = 0$ and Dirichlet boundary conditions on the surface of the sphere.
Let us remark that this initial condition differs form the one in Section \ref{sec:CE} where 
density and pressure values were lower yielding to a Mach number greater than one, {\color{black} whereas in this simulation $M\approx 0.3$}. The simulation is run until time $t_{\textrm{end}}=0.25$ on two different meshes, M1 consisting of $557147$ tetrahedra and M2 made of $2280182$ primal elements.
In Figure \ref{fig:CE3d_Le_t025}, we have plotted the solution obtained using the LADER-ENO scheme with auxiliary artificial viscosity $c_{\alpha}=3$. We observe a good agreement with the reference solution that has been obtained using the 1D code which solves compressible Euler equations with appropriate geometrical source terms, see \cite{Toro}. {\color{black} The CPU time employed for a serial simulation in M1 on an Intel$^{\textrm{\textregistered}}$  Xeon$^{\textrm{\textregistered}}$ Gold 6126 is of $125666.58s$, corresponding to a CPU time per dual element and iteration $t_{e}= 82.92\mu s$ whereas for M2 we get CPU time = $773741.77s$ and  $t_{e}= 77.65\mu s$  .}

\section{Conclusions}\label{sec:conclusions}
In this paper we have presented a novel semi-implicit hybrid finite volume-finite element method for the simulation of weakly compressible flows in two and three space dimensions. It corresponds to an extension of the pressure-based solver introduced in \cite{BFSV14,BFTVC17} for incompressible flows. The unstructured staggered meshes considered allow for an easy discretization of complex domains, while avoiding the checker-board phenomena typical for collocated grids. 

Within this paper, the original compressible Navier-Stokes equations have been rewritten by replacing the total energy conservation equation with the governing PDE for the pressure, whose formal derivation has been also included. The time discretization has been then performed on the resulting system of equations and subsequently a projection algorithm has been applied. One of the main advantages of the developed semi-implicit methodology is that it allows for the  decoupling of the density and linear momentum variables and the computation of the pressure. The former ones are obtained at the transport diffusion stage using a finite volume scheme. More precisely, we have extended the LADER methodology introduced in \cite{BFTVC17} to account for density variations. Consequently, the explicit scheme used for the transport diffusion stage is second order accurate in space and time. Moreover, this technique profits from the dual mesh structure, resulting in an efficient numerical scheme with a small stencil for the final pressure system to be solved in each time step. The use of a projection technique splits the original pressure equation into three parts. The first of them includes a non conservative product that is approximated in the finite volume framework using a path conservative scheme. The second part consists in a density dependent term that has been computed on the primal mesh using a finite volume approach.  Finally, the third part is coupled with the second PDE derived from the momentum conservation equation, yielding a Poisson-type system for the pressure unknown.  
In the case when the sound speed tends to infinity, the pressure system reduces identically to the pressure Poisson equation of the incompressible Navier-Stokes equations. The pressure system is solved in the projection stage by employing a classical $\mathbb{P}_{1}$ continuous finite element method. The temperature and heat flux needed at each iteration are obtained, using the EOS, as a postprocessing of the previous time step. Passing data from one mesh to the other is done thanks to a weighted average. 

Numerous tests have been presented, aiming to validate the final algorithm. First, the order of accuracy has been assessed numerically at the aid of the Taylor-Green vortex benchmark problem, for which an exact solution of the incompressible Navier-Stokes equations is available. 
{\color{black} The computational efficiency of the new scheme proposed in this paper has been carefully studied by comparing it against a fully incompressible flow solver and a fully explicit density-based Godunov-type finite volume scheme, implemented in the same code basis and using the same mesh and the same computer.} The behaviour of the scheme in the presence of weak discontinuities has been tested via several Riemann problems. Moreover, a 2D circular explosion problem has been presented and compared with a reference solution that makes use of the angular symmetry of the problem. A second test case with known analytical solution and which was studied in this paper was the first problem of Stokes, in which viscous effects play a dominant role.  
To show the capability of the method to deal with very low Mach number flows, we have also considered the lid driven cavity test and the double shear layer benchmark. The obtained results have been successfully compared with reference solutions available in the literature.  
A smooth acoustic wave test has been run to check the correct propagation of sound waves. Like for the circular explosion test, we have observed an excellent agreement between the results given by the Hybrid FV/FE method and a reference solution that makes use of the angular symmetry of the problem. Regarding heat driven flows, we have included two different tests: a simple one-dimensional heat conduction test, initiated by a temperature jump in the initial condition, as well as rising bubble test, that has been validated against the semi-implicit high order discontinuous Galerkin scheme presented in \cite{TD17,BTBD20}.  Finally, also a three dimensional spherical explosion problem has been included, showing the good performance of the methodology also in three space dimensions.  

As future research, we will extend the above methodology to solve all Mach number flows by using a discrete form of the total energy conservation equation, instead of the pressure equation used in this paper. Moreover, the parallelization of the code will be considered, aiming at decreasing the wall clock time consumption and allowing for the simulation of more complex and realistic problems in three space dimensions.

\section*{Acknowledgements}
This work was financially supported by INdAM (\textit{Istituto Nazionale di Alta Matematica}, Italy) under a Post-doctoral grant of the research project \textit{Progetto premiale FOE 2014-SIES} and by the Spanish MECD under grant FPU13/00279; by Spanish MICINN projects MTM2013-43745-R, MTM2015-68275-R; by Spanish MCIU under project MTM2017-86459-R; by FEDER and Xunta de Galicia funds under the ED431C 2017/60 project. 
S.B. and M.D. acknowledge funding from the Italian Ministry of Education, University 
and Research (MIUR) in the frame of the Departments of Excellence Initiative 2018--2022 
attributed to DICAM of the University of Trento (grant L. 232/2016) and in the frame of the 
PRIN 2017 project \textit{Innovative numerical methods for evolutionary partial differential equations and  applications}. Furthermore, M.D. has also received funding from the University of Trento via the Strategic Initiative \textit{Modeling and Simulation} and acknowledges partial support 
of the European Union's Horizon 2020 Research and Innovation Programme under 
the project \textit{ExaHyPE}, grant no. 671698 (call FETHPC-1-2014).  
S.B. and M.D. are members of the GNCS-INdAM group.

\bibliographystyle{elsarticle-num}
\bibliography{./mibiblio}

\appendix
\section{Energy equation in terms of pressure}\label{sec:app_equations}	
The goal of this appendix is to write the energy equation in terms of pressure as main variable. Firstly, we recall that the speed of sound $c$ is defined by 
$$
c^2:=\Big(\frac{\partial \press}{\partial \rho}\Big)_s,
$$
where $\press$ denotes pressure,  $s$ specific entropy,  $\rho$ density and  $\nu=\rho^{-1}$ specific volume. We notice that $c$ is a function of two thermodynamic variables, for instance, $\press$ and $\nu$ or $\rho$. Let $e$ denote the specific internal energy:
$$
e:=E-\frac{1}{2}|\bf u|^2, 
$$
where $E$ is the specific total energy. We will use the following thermodynamic equalities:
\begin{lemma}
	We have
	\begin{equation}
		\label{eq:1}
		c^2=\Big(\frac{\partial \press}{\partial \rho}\Big)_e+\frac{\press}{\rho^2}\Big(\frac{\partial \press}{\partial e}\Big)_\rho,
	\end{equation}
	\begin{equation}
		\label{eq:energy6}
		\Big(\frac{\partial p}{\partial e}\Big)_\rho=\frac{1}{c_v} \Big(\frac{\partial p}{\partial \theta}\Big)_\rho,
	\end{equation}
	where $c_v$ denotes the specific heat at constant volume defined by
	$$
	c_v:= \Big(\frac{\partial e}{\partial \theta}\Big)_\rho.
	$$
\end{lemma}
In what follows the dot over a field denotes its material derivative with respect to time. Let us recall that $\dot{\varphi}=\displaystyle\frac{\partial\varphi}{\partial t} +{\bf u}\cdot \grae \varphi$.
\begin{proposition}
	The energy conservation equation can be written as
	\begin{equation}
		\label{eq:energy4}
		\dot{\press}+\rho c^2\dive {\bf u}= \frac{1}{\rho c_v}\Big(\frac{\partial \press}{\partial \theta}\Big)_\rho \big( \boldsymbol{\tau}\cdot \gra \mathbf{u}- \dive \mathbf{q} +f\big),
	\end{equation}
	\noindent where  ${\bf q}$ is the heat flux vector and $f$ is the volumetric heat source density.
\end{proposition}
\begin{proof} We start from the following form of the energy conservation equation in terms of the total energy $E$:
	\begin{equation}
		\label{ETLO}
		\rho\dot{E}=\dive (\boldsymbol{\sigma}\mathbf{u})+\mathbf{b}\cdot \mathbf{u}-\dive \mathbf{q}+f,
	\end{equation}
	where  $\boldsymbol{\sigma}$ is the Cauchy stress tensor and $\mathbf{b}$ is the body force.
	By using the momentum equation it is easy to prove that (see, for instance, \cite[Prop. 1.4.6]{Ber05}):
	\begin{equation}
		\label{eq:energy1}
		\rho \dot{e} =\boldsymbol{\sigma}\cdot \gra \mathbf{u} - \dive \mathbf{q} +f,
	\end{equation}
	For fluids,
	$$
	\boldsymbol{\sigma}=-\press\,  \mathbf{I} + \boldsymbol{\tau},
	$$
	being $\boldsymbol{\tau}$ the viscous stress tensor. Replacing in~\eqref{eq:energy1} we obtain
	\begin{equation}
		\label{eq:energy2}
		\rho \dot{e} =-\press \dive \mathbf{u} + \boldsymbol{\tau}\cdot \gra \mathbf{u} - \dive \mathbf{q} +f.
	\end{equation}
	Moreover, writing $\press=\hat{\press}(\rho,e)$ we have,
	\begin{align}
		\dot{\press}= \Big(\frac{\partial \press}{\partial \rho}\Big)_e \dot{\rho} + \Big(\frac{\partial \press}{\partial e}\Big)_\rho \dot{e}
	\end{align}
	and using the mass conservation equation, $\dot{\rho} +\rho \mathrm{div} {\bf u}=0$, \eqref{eq:1} and~\eqref{eq:energy2} we deduce
	\begin{align}
		\label{eq:energy3}
		\nonumber
		\dot{\press}= -\Big(\frac{\partial \press}{\partial \rho}\Big)_e \rho \dive {\bf u}+ \frac{1}{\rho}\Big(\frac{\partial \press}{\partial e}\Big)_\rho \big(-\press\dive  {\bf u} + \boldsymbol{\tau}\cdot \gra \mathbf{u} - \dive {\bf q} +f\big)
		\\
		\nonumber
		=-\rho\dive  \mathbf{u}\Big[\Big(\frac{\partial \press}{\partial \rho}\Big)_e + \frac{\press}{\rho^2}\Big(\frac{\partial \press}{\partial e}\Big)_\rho \Big]+ \frac{1}{\rho}\Big(\frac{\partial \press}{\partial e}\Big)_\rho \big(\boldsymbol{\tau}\cdot \gra\mathbf{u} - \dive {\bf q} +f\big)
		\\
		= -\rho c^2\dive  \mathbf{u} + \frac{1}{\rho}\Big(\frac{\partial \press}{\partial e}\Big)_\rho \big( \boldsymbol{\tau}\cdot \gra \mathbf{u} - \dive {\bf q} +f\big),
	\end{align}
	which leads to~\eqref{eq:energy4} by using~\eqref{eq:energy6}.
\end{proof}
\begin{corollary}
	The following equation holds:
	\begin{equation}
		\label{eq:energy8}
		\frac{\partial \press}{\partial t}+ {\bf u}\cdot(\grae \press-c^2\grae \rho)+c^2\dive  \left(\rho \mathbf{u}\right)= \frac{1}{\rho c_v}\Big(\frac{\partial \press}{\partial \theta}\Big)_\rho \big( \boldsymbol{\tau}\cdot \gra\mathbf{u}- \dive {\bf q} +f\big).
	\end{equation}
\end{corollary}
\begin{proof}. It easily follows from~\eqref{eq:energy4} by using the equalities
	\begin{align*}
		\dot{\press}=\frac{\partial \press}{\partial t}+\mathbf{u}\cdot \grae \press,
		\\
		\rho\dive {\bf u}=\dive \left(\rho\mathbf{u}\right)-\grae \rho\cdot \mathbf{u}.
	\end{align*}
\end{proof}
We notice that equation~\eqref{eq:energy4} is valid for any fluid, independently of its constitutive law. In the case of ideal gases,
$$
\press=\rho R\theta,
$$
being $R$ the specific gas constant and $\theta$ the absolute temperature. Then,
$$
\frac{1}{\rho c_v}\Big(\frac{\partial \press}{\partial \theta}\Big)_\rho=\frac{1}{\rho c_v}{\rho R}=\frac{R}{c_v}=\frac{c_{\press}-c_v}{c_v}=\gamma-1,
$$
where $c_{\press}$ is the specific heat at constant pressure, namely,
$$c_{\press}:= \Big(\frac{\partial h}{\partial \theta}\Big)_{\press},
$$
being $h:= e+\displaystyle \frac{\press}{\rho}$ the specific enthalpy, $\gamma:=\displaystyle\frac{c_{\press}}{c_v}$ is the adiabatic index, and we have used the Mayer relation $R=c_{\press}-c_v$. Thus, we have proved the following form of the energy equation for ideal gases.
\begin{corollary}
	For ideal gases, equation~\eqref{eq:energy8} becomes
	\begin{equation}
		\label{eq:isentropic4}
		\frac{\partial \press}{\partial t}+ \mathbf{u}\cdot(\grae \press-c^2\grae \rho)+c^2\dive  \left(\rho \mathbf{u}\right)= (\gamma-1) \big( \boldsymbol{\tau}\cdot \gra \mathbf{u}- \dive \mathbf{q} +f\big).
	\end{equation}
\end{corollary}
Let us remark that if the volumetric heat source density is zero, then equation \eqref{eq:isentropic4} results
\begin{equation}
	\frac{\partial \press}{\partial t}+ {\bf u}\cdot(\grae \press-c^2\grae\rho)+c^2\dive \left(\rho\mathbf{u}\right)= (\gamma-1) \big( \boldsymbol{\tau}\cdot \gra \mathbf{u}- \dive \mathbf{q} \big).\label{eq:energy_igsf}
\end{equation}
which corresponds with equation \eqref{eq:pressure2}.

\end{document}